\documentclass[11pt]{article}
\usepackage{eurosym}
\usepackage{amsfonts}
\usepackage{amsmath}
\usepackage{geometry}
\usepackage{amssymb}
\usepackage{graphicx}
\usepackage[dvipsnames]{xcolor}
\usepackage{hyperref}
\hypersetup{
    colorlinks=true,
    citecolor=Green,
    linkcolor=blue,
    filecolor=magenta,      
    urlcolor=cyan,
    }
\usepackage[author={Oana}]{pdfcomment}

\usepackage{titling}
\usepackage{mathrsfs}
\usepackage{enumitem}
\usepackage{subcaption}
\usepackage[color=yellow,colorinlistoftodos]{todonotes}
\usepackage{framed}
\usepackage[toc,page]{appendix}

\newtheorem{theorem}{Theorem}

\newtheorem{remark}[theorem]{Remark}

\newcommand{\id}{\,\mathrm{d}}
\newcommand{\Z}{\mathbb{Z}}
\newcommand{\R}{\mathbb{R}}
\newcommand{\C}{\mathbb{C}}
\newcommand \N{\mathbb{N}}

\newcommand{\mb}[1]{\mathbf{#1}}
\newcommand{\mbb}[1]{\mathbb{#1}}

\newcommand{\bs}[1]{\boldsymbol{#1}}

\newcommand{\dt}{\, \mathrm{d}t}

\def\p{{\partial}}

\def\rmd{{\color{black}{\rm d}}}
\def\ba{{\mathbf{a}}}

\def\bb{{\mathbf{b}}}
\def\bc{{\mathbf{c}}}
\def\bd{{\mathbf{d}}}
\def\bh{{\mathbf{h}}}
\def\bk{{\mathbf{k}}}
\def\bp{{\mathbf{p}}}
\def\bq{{\mathbf{q}}}

\def\bu{{\mathbf{u}}}
\def\bv{{\mathbf{v}}}
\def\bx{{\mathbf{x}}}

\newcommand{\mbs}[1]{\ensuremath{\boldsymbol{#1}}}

\begin{document}

\date{}
\title{\textbf{Comparison of Stochastic Parametrization Schemes using Data Assimilation on Triad Models}}
\author{B. Chapron, D. Crisan, D. Holm, O. Lang, A. Lobbe, E. M\'{e}min}

\maketitle


\abstract{ 
In recent years, stochastic parametrizations have been ubiquitous in modelling
uncertainty in fluid dynamics models. One source of 
model uncertainty comes from the coarse graining of the fine-scale data and is in common usage in computational simulations at coarser scales. In this paper, we look at two such stochastic parametrizations: the Stochastic Advection by Lie Transport (SALT) parametrization introduced by Holm in \cite{dh2015} and the Location Uncertainty (LU) parametrization introduced by M\'{e}min in 
\cite{em2014}. Whilst both parametrizations are available for full-scale models, we study their reduced order versions obtained by projecting them on a complex vector Fourier mode triad of eigenfunctions of the curl. Remarkably, these two parametrizations lead to the same reduced order model, which we term the \textit{helicity-preserving stochastic triad} (HST). This reduced order model is then compared with an alternative model which preserves the energy of the system, and which is termed the \textit{energy preserving stochastic triad} (EST). These low-dimensional models are ideal benchmark models for testing new Data Assimilation algorithms: they are easy to implement, exhibit diverse behaviours depending on the choice of the coefficients and come with natural physical properties such as the conservation of energy and helicity. 
}



\section{Introduction} 

The introduction of stochasticity in fluid dynamics has recently been the subject of intense research effort. This approach involves using random processes to model, for example, unresolved scales, or to take into account neglected physical effects. 
A stochastic formulation for the fluid flow introduces a probabilistic basis for modelling unresolved scales. This is different from the deterministic causal modelling which is  difficult to achieve in practice due, for instance, to unknown initial conditions. In addition, some phenomena such as energy backscattering are directly accessible as stochastic processes. Another usage of stochastic modelling is to generate ensembles of realizations of the model. This facilitates the analysis of model uncertainty quantification for different low-resolution computational simulations and their usage to approximate the true state of the fluid, instead of using a single high-resolution numerical simulation. 

Some stochastic schemes have been proposed in the literature  by considering a variety of \textit{ad-hoc} perturbations. However, a principled approach is desirable.
The formulation of stochastic dynamical systems based on physical principles has recently been proposed in various settings. For a review and classification of approaches to stochastic parameterisation based on physical principles, see \cite{berner2017stochastic}. 
The present work treats two additional new approaches. The first one, called  stochastic advection by Lie transport (SALT) relies on the variational principle for fluid dynamics \cite{dh2015}. The second one, called modelling under location uncertainty (LU) is derived from Newton's principle \cite{em2014}. Both frameworks introduce stochasticity into the Lagrangian specification of the flow field, rather than directly into the Eulerian frame.  

{In the deterministic case, it is known that three-dimensional fluid flows may trigger a cascade of dynamics across multiple length and time scales. This multiscale behavior poses considerable challenges in the computational simulation using standard Navier-Stokes equations (see e.g. \cite{bohr1998dynamical}, \cite{novikov}, \cite{Biferale}). }
When modelling turbulence numerically, {specialised} discretisation methods {are needed to} decompose the underlying partial differential equations into a very large number of ordinary differential equations. Alternative approaches have been introduced where the Navier-Stokes dynamics in the Fourier space is mimicked using a finite number of variables, {say} $u^1, u^2, \ldots, u^N$. The Fourier space is divided into $N$ shells, and each shell $\mathrm{s}_i$ comprises the set of wave vectors $\underline{s}$ with magnitude $|\underline{s}|\in(s_02^i, s_02^{i+1})$. Each $u^{i}$ satisfies an ODE and it represents the magnitude of the velocity field on a length scale of $s_i^{-1}$ (\cite{gledzer}, \cite{bohr1998dynamical}). The quadratic nonlinearity in the Navier-Stokes equations produces triads of interacting vector Fourier modes within each shell. Shell models involving multiple triads have had considerable success in modelling energy and helicity cascades, as well as modelling intermittency in chaotic dynamical systems \cite{Biferale-2003,chen2003intermittency,chen2003joint}. Simplified shell models with only a few triads date back to the 1970s and have provided {major} insight into fluid modal interaction. Even the dynamical system representation of Euler's fluid equations on a single triad has been  
quite insightful, see e.g. \cite{waleffe1992nature,waleffe1993inertial}. 

More recently, the problem of correctly parameterizing effects of small-scale physical processes together with the need for probabilistic ensemble forecasting and uncertainty quantification has led to modern stochastic approaches in the study of turbulence using reduced order \textit{shell models}. 
In this work we will explore reduced order models for SALT and LU models obtained by projecting onto helical basis functions \cite{Biferale-2003,waleffe1992nature,waleffe1993inertial}. These helical basis functions, defined as eigenfunctions of the curl operator, enable one to construct reduced order stochastic models of fluid flow with a simplified nonlinear interaction. As we will see, under projection onto the basis of helical triad modes, both LU and SALT result in the \textit{same} reduced order model and this projected model conserves helicity, but it does not conserve energy. Because of this coincidence in projecting the SALT and LU models onto the helical basis, a second reduced order scheme with a strong energy conservation property inspired by \cite{holm2021stochastic} and known as the \textit{energy preserving stochastic triad} (EST) model will be proposed for  comparison. 

While the EST model is not of transport type, it will provide comparison between two different classes of stochastic dynamical systems. The two classes treated here are (1) the \textit{helicity preserving stochastic triad} (HST) (comprising both LU and SALT on the helical basis) and (2) the \textit{energy preserving stochastic triad} (EST) of \cite{holm2021stochastic} projected onto the helical basis. The solution behaviour of the HST model will be compared to that for the EST model for several data assimilation objectives formulated on the helical triad modes. 
%
%
For classical deterministic models one obtains a system of ordinary differential equations. However, for stochastic dynamics a set of  stochastic differential equations (SDEs) is obtained \cite{chapron2018largescale,geurts2020lyapunov,Resseguier-SIAMUQ-22}.




The goal of the data assimilation procedure in this context is twofold: firstly, it is used to calibrate the uncertainty of the model (the amplitude of the noise). Secondly, once the calibration is complete, the particle filtering methodology can be used to reduce the uncertainty. We want the distribution of the fluctuations to be properly approximated. In the absence of stochasticity, all particles would go in the same direction and the initial spread would rapidly disappear because of the hyperbolic character of the model. In the absence of a reasonable spread, the particle filter methodology will eventually collapse. For this reason, we need to introduce stochasticity into the system that correctly  characterises the fluctuation dynamics. In particular, we want to find the type of noise amplitude and the stochastic parameters for which the distribution of the output samples is reasonably \textit{uniform}.

\paragraph{Structure of the Paper}

In Section~\ref{sec:roms} we introduce the triad models for incompressible flows modelled by the Euler equation in its deterministic and stochastic form. To this end, we introduce the stochastic parametrisation paradigms. 
Building upon these models for the 3D Euler equation, we then present reduced order triad models derived from the original equations. The derivation follows the classical approach for triad models from the literature, that has been successfully employed in the deterministic case. Our full derivations, complete also for the stochastically parametrised models, can be found in Appendix~\ref{sec:derivations}. The Data Assimilation experiments are carried out in Section \ref{sect:dataassimilation}. We first briefly explain the standard particle filter methodology, and then in Section~\ref{subsect:numericalstudies} we present the findings of our numerical studies. 
In particular, Section~\ref{subsect:numericalstudies} presents the results of the main numerical studies in this work. These are:
\begin{itemize}
    \item The model realisations of the stochastic models which we propose in this work are presented in Section~\ref{sec:model_realisation} for different realisations of the noise. Further, we numerically confirm here the physical conservation properties of the models that motivated their theoretical design.
    \item The model statistics for a large number of stochastic realisations are presented in Section~\ref{sec:model_stats}. Crucially, the evolution of the model statistics in time reveals differences between the two stochastic models on the level of the individual triad energies which go beyond the conservation properties. Moreover, we show here the (in-)stability of either model with respect to large noise parameters. 
    \item In Section~\ref{sec:Numerical_DA} we exhibit the results of the data assimilation performed using the transition kernels derived from the theoretical models developed in this paper. We show here that using the stochastic transition kernels associated with our proposed stochastic models improves the particle filtering procedure and produces efficient ensemble evolutions that are well-suited for data assimilation purposes.
\end{itemize}
In Section \ref{sect:conclusions} we describe our conclusions on this topic. 
We conclude the paper with a number of appendices: in Appendix \ref{appendix:notations} one can find a list of notations and standard identities, in Appendix \ref{appendix:derivTriademodels} we present a detailed derivation of shell models (deterministic and stochastic), in Appendix \ref{sec:suppl_num} we introduce some supplementary numerics related to the noise amplitude calibration. 

\paragraph{Code Availability}
The code corresponding to the numerical experiments in this paper is archived in~\cite{alobbe_2023_7845270}. The GitHub repository is located at \href{https://github.com/alobbe/stochastic-triads}{https://github.com/alobbe/stochastic-triads}.
\section{Reduced Order Models for Incompressible Fluids}
\label{sec:roms}

\subsection{Reduced Order Models for the 3D Euler Equation}
\label{sec:StochParam}

The 3D Euler equations model incompressible inviscid fluid dynamics. These equations may be written 
 by using the Leray operator $\mathcal{P}$ to project onto the divergence-free part  of its operand as
\begin{align}
\begin{split}
\frac{\p \mb{v}}{\p t} &=  \mathcal{P}\Big( \mb{v} \times {\rm curl}\mb{v}\Big)
\\&= 
\mathcal{P} \bigg( \frac{\delta E}{\delta  \mb{v}} \times \frac{\delta C}{\delta  \mb{v}} \bigg)
\\ \hbox{with conserved Energy } E( \mb{v}) &=   \int_{\mbb{R}^3} \frac12\mathcal{P}\mb{v}\cdot \mb{v}\it  d^3x 
\\ \hbox{and conserved helicity } C( \mb{v}) &=  \frac12\int_{\mbb{R}^3}  \mb{v} \cdot {\rm curl}\,\mb{v}\, d^3x 
\,
\end{split}
\label{eqn: Euler}
\end{align}
where $\delta/\delta \mb{v}$ represents variational derivative with respect to the fluid velocity $\mb{v}$. 

Following \cite{chen2003joint,chen2003intermittency} we use a Galerkin expansion in orthogonal \emph{vector} modes that are eigenfunctions of the curl operator. Assume the fluid is contained in a periodic box $\mathcal{D}\subset\R^3$ of side length $L>0$. Then the velocity $\bv$ and vorticity $\bs{\omega}:=\nabla\times\bv$ may be expanded in circularly polarised or helical modes $\bh_{\pm}(\bk)\exp(i\bk\cdot\bx)$, with the wave vectors $\bk\in\mathcal{K}:=(2\pi/L)\Z^3$. The modes shall be orthogonal, i.e.
\begin{equation}
    \int_{\mathcal{D}}\bh_{s_p}(\bp)\exp(i\bp\cdot\bx)\cdot [\bh_{s_q}(\bq)\exp(i\bq\cdot\bx)]^* \id\bx = C\delta_{\bp,\bq}\delta_{s_p,s_q}\,;\; C>0 \text{ const}.
\end{equation}
The complex vector amplitudes $\bh_{\pm}(\bk)$ should satisfy $\bk\cdot\bh_{\pm}(\bk)=0$ and $i\bk\times\bh_{\pm}(\bk)=\pm |\bk|\bh_{\pm}(\bk)$.
A convenient choice of basis for the $\bh_{\pm}(\bk)$ is then given by
\begin{align} 
    \bh_{\pm}(\bk)&:=\bs{\nu}\times\bs{\kappa}\pm i\bs{\nu},
    \; \text{ with } \;
    \bs{\kappa}:=\bk/k,
    \; \bs{\nu} := \bk\times\bs{\Gamma}/|\bk\times\bs{\Gamma}|,\; \bs{\Gamma}:=\text{ const},
\label{eqn: h-basis}
\end{align}
for which $|\bh_\pm(\bk)|^2:=\bh_\pm(\bk)\cdot\bh_\pm(\bk)^*=2$ and $\bh_{\pm}(\bk)\cdot\bh_{\mp}(\bk)^*=0$.

At this point, one notices the key features of the helical modes $\bh_{\pm}(\bk)\exp(i\bk\cdot\bx)$ which greatly simplifies analysis of modal expansions of the 3D Euler and related equations, such as 3D Navier-Stokes. Namely, the helical modes $\bh_{\pm}(\bk)\exp(i\bk\cdot\bx)$, are divergence-free eigenfunctions of the curl operator. Specifically,
\begin{equation}
   \nabla \cdot  \bh_s(\bk) e^{i\bk\cdot\bx} =  i\bk\cdot\bh_s(\bk) e^{i\bk\cdot\bx} = 0
\end{equation}
and
\begin{equation}
    \nabla \times \bh_s(\bk) e^{i\bk\cdot\bx} = i\bk\times\bh_s(\bk)e^{i\bk\cdot\bx} 
    = s |\bk|\bh_s(\bk)e^{i\bk\cdot\bx}
    .
\label{eqn: eigencurl}
\end{equation}
See \cite{chen2003joint,chen2003intermittency,waleffe1992nature,waleffe1993inertial} for more information about how this Galerkin decomposition into divergence-free eigenfunctions of the curl are used as a standard tool in analysis of 3D solution behaviour of the deterministic Euler fluid equations and Navier-Stokes fluid equations. In particular, the helical mode expansions in equations \eqref{eqn: defs} comprise the source of the popular \emph{shell models} as finite-dimensional expansions of turbulent fluid dynamics. Thus, this expansion provides a useful framework for studying low-dimension stochastic models of 3D Navier-Stokes turbulence.\\

In terms of the basis of helical modes $\bh_{\pm}(\bk)\exp(i\bk\cdot\bx)$, the divergence-free fluid velocity $\bv(\bx,t)$ and vorticity $\bs{\omega}(\bx,t)$ are expressed in \cite{waleffe1992nature,waleffe1993inertial} in terms of complex vector amplitudes $\bu(\bk,t),\bs{\varpi}(\bk,t)\in \mbb{C}^3$, respectively,
\begin{align}
\begin{split}
\bv(\bx,t) &:= \sum_{\bp} \bu(\bp,t)e^{i\bp\cdot\bx} := \sum_{\bp} \sum_{s_p=\pm} a_{s_p}(\bp,t)\bh_{s_p}(\bp) e^{i\bp\cdot\bx}
\,,\\
\bs{\omega}(\bx,t) &:= \sum_{\bq} \bs{\varpi}(\bq,t)e^{i\bq\cdot\bx} :=\sum_{\bq} \sum_{s_q=\pm} s_q|\bq|\, a_{s_q}(\bq,t)\bh_{s_q}(\bq) e^{i\bq\cdot\bx}
\,.
\end{split}
\label{eqn: defs}
\end{align}
Here, the choice 
\begin{equation}
    \bu(\bk,t) := a_+(\bk,t)\bh_+(\bk) + a_-(\bk,t)\bh_-(\bk) = \sum_{s_k=\pm} a_{s_k}(\bk,t)\bh_{s_k}(\bk),
    \label{eqn:u_amplitudes}
\end{equation}
with $a_s^*(\bk) = a_s( - \bk) $ \cite{waleffe1993inertial},
was made, so that \eqref{eqn: eigencurl} implies
\begin{equation}
    \bs{\varpi}(\bk,t):= |\bk| \Big(a_+(\bk,t)\bh_+(\bk) - a_-(\bk,t)\bh_-(\bk)\Big) =  |\bk|\!\!\sum_{s_k=\pm}  s_k\, a_{s_k}(\bk,t)\bh_{s_k}(\bk).
    \label{eqn:omega_amplitudes}
\end{equation}
The conservation laws for the Euler fluid kinetic energy and helicity -- expressed as integrals over the spatially periodic box ${\cal D}$ -- can be evaluated in Fourier space via Parseval's theorem, as follows,
\begin{align}
\begin{split}
\frac12\int_{\cal D} |\bv(\bx,t)|^2 \, d^3x 
&=
\frac12 \sum_{\bk}  \bu(\bk,t)\cdot \bu^*(\bk,t) = \sum_{\bk} \sum_{s_k=\pm} a_{s_k}(\bk,t)\cdot a^*_{s_k}(\bk,t)
\,,\\
\int_{\cal D} \bv(\bx,t)\cdot {\rm curl} \bv(\bx,t) \, d^3x
&=
 \sum_{\bk}  \bu(\bk,t) \cdot \bs{\varpi}^*(\bk,t)
\\
&=
 \sum_{\bk}\sum_{s_k=\pm}  ks_k \,a_{s_k}(\bk,t)a^*_{s_k}(\bk,t) 
\,\bh_{s_k}(\bk)  \cdot \bh^*_{s_k}(\bk)
\\
&=
2 \sum_{\bk}\sum_{s_k=\pm}  ks_k \,a_{s_k}(\bk,t)a^*_{s_k}(\bk,t)  
\,.
\end{split}
\label{eqn: helicity}
\end{align}
Expanding the terms of the Euler equation in curl form \eqref{eqn: Eulcurl}, we obtain the Euler equations for the coefficients $a_{s_k}(\bk,t)$. For all $\bk\in\mathcal{K}$, $s_k\in \{+,-\}$ we have
\begin{equation}
    \partial_t a_{s_k}(\bk, t) = -\frac{1}{4} \sum_{\bp+\bq+\bk=0}\sum_{s_p, s_q} (s_p|\bp|-s_q|\bq|)a_{s_p}^*(\bp, t)a_{s_q}^*(\bq, t)\bh_{s_p}^*(\bp)\times\bh_{s_q}^*(\bq)\cdot \bh_{s_k}^*(\bk).
    \label{eqn:Euler_coefs}
\end{equation}
For the explicit derivation of equation \eqref{eqn:Euler_coefs} see Appendix~\ref{adx:detEuler}.

The elementary interactions in Fourier space take place
between triads of wave vectors such that $\bk + \bp + \bq = 0$, as exemplified in equation \eqref{eqn:Euler_coefs} above. There are two degrees of freedom per wave vector, $(a_+ ,a_-)$, so eight different types of interaction are allowed according to the value of the triplet $(s_k ,s_p ,s_q ) = ( \pm 1, \pm 1, \pm 1)$.
Consider a fixed triple of wave vectors $\bk, \bp, \bq\in \mathcal{K}$ such that $\bk+\bp+\bq=0$ and a fixed triple $s_k,s_p,s_q\in \{+,-\}$. This gives rise to three coefficients $a_{s_k}(\bk,t), a_{s_p}(\bp,t), a_{s_q}(\bq,t)$, which we compactly summarise into the complex vector \[\ba = (a_{s_k}, a_{s_p}, a_{s_q})\in\C^3\,.\] The dynamics of $\ba$ is determined by the three equations obtained from \eqref{eqn:Euler_coefs}
\begin{align}
    \frac{\rmd a_{s_k}}{\rmd t} &= g(s_p|\bp| - s_q|\bq|) a_{s_p}^*a_{s_q}^*+ R,\\
    \frac{\rmd a_{s_p}}{\rmd t} &= g(s_q|\bq| - s_k|\bk|) a_{s_q}^*a_{s_k}^*+R,\\
    \frac{\rmd a_{s_q}}{\rmd t} &= g(s_k|\bk| - s_p|\bp|) a_{s_k}^*a_{s_p}^*+R,
\end{align}
 We pick out a summand and cycle through $k,p,q$
\begin{align}
    \frac{\rmd a_{s_k}}{\rmd t} &= g(s_p|\bp| - s_q|\bq|) a_{s_p}^*a_{s_q}^*,\\
    \frac{\rmd a_{s_p}}{\rmd t} &= g(s_q|\bq| - s_k|\bk|) a_{s_q}^*a_{s_k}^*,\\
    \frac{\rmd a_{s_q}}{\rmd t} &= g(s_k|\bk| - s_p|\bp|) a_{s_k}^*a_{s_p}^*,
\end{align}
with the constant complex scalar
\begin{equation}
    g := -\frac{1}{4}\bh_{s_p}^*(\bp)\times\bh_{s_q}^*(\bq)\cdot \bh_{s_k}^*(\bk).
    \label{eq:g-factor}
\end{equation}
The equations corresponding to the single triad interaction of type  $(s_k, s_p, s_q)$ with $\bk + \bp+\bq=0$ thus have the complex vector form also derived in \cite{waleffe1992nature},
\begin{equation}
\frac{\rmd\ba}{\rmd t} = g \ba^* \times \mathbb{D}\ba^* = g (\ba \times \mbb{D}\ba)^*, 
\label{eqn: complexRB}
\end{equation}
with the constant diagonal matrix \begin{equation}
    \mbb{D}:={\rm diag}\big(s_k|\bk|, s_p|\bp|, s_q|\bq| \big).
\end{equation}
The form of the factor $g$ defined in \eqref{eq:g-factor} above can be calculated from \eqref{eqn: h-basis} to show that it depends on the shape and the orientation of the wave-vector triad, but not on its scale; since the real and imaginary parts of the complex helical vector amplitudes $\bh_{s_k}$ are unit vectors. Moreover, $\mbb{D}\ba$ can be seen to represent the $(s_k, s_p, s_q)$ components of the vorticity vector amplitude $\bs{\varpi}$ through equation \eqref{eqn:omega_amplitudes} above.
Two conservation laws for real-valued triad energy and helicity follow immediately from equation \eqref{eqn: complexRB}, as 
\begin{equation}
\frac{d}{dt}(\ba\cdot\ba^*) =  0
\quad\hbox{and}\quad
\frac{d}{dt}(\ba\cdot\mbb{D}\ba^*) =  0
\,.
\label{eqn: erghelic}
\end{equation}
The dynamical system in equation \eqref{eqn: complexRB} is similar to rigid body dynamics, but replaced by complex angular momentum $(\ba)$, complex angular velocity $(\mbb{D}\ba)$ and real moment of inertia $\mbb{I}=\mbb{D}^{-1}$ relating the two complex quantities.\footnote{Rigid body dynamics with complex angular momentum has  also been discussed previously in \cite{bender2007complexified}.}


\subsection{Stochastic parametrizations for the 3D Euler equation}

In the following we introduce  a reduced order model
for the stochastic parametrizations introduced through
the Stochastic Advection by Lie Transport paradigm as well as the Location Uncertainty paradigm. We explain below the rationale of these parametrizations:

\subsubsection{Modelling under the Stochastic Advection by Lie Transport Principle}

The SALT equations were derived in \cite{dh2015} using a Stratonovich stochastic version of Hamilton's variational principle \cite{holm1998euler} in combination with Kraichnan's scalar turbulence model based on Stratonvich stochastic Lagrangian paths \cite{kraichnan1968small}. 
The application of Hamilton's principle with an imposed stochastic Lie transport constraint implied an Euler-Poincar\'e equation for the fluid motion \cite{holm1998euler}. The 3D SALT Euler equations for divergence-free fluid velocity $\bv(\bx,t)$ are given by
\begin{align}
\begin{split}
    &\id \bv + (\id \bx_t \cdot \nabla) \bv + v_j\nabla \id \bx_t^j = - \nabla \id p \,,\\
    &\quad\hbox{with}\quad\id \bx_t = \bv{dt} + \sum_i \bs{\xi}_i(\bx)\circ \id W_t^i\,.
\end{split}
\label{eqn: velformSALT}
\end{align}
As discussed in \cite{dh2015}, this motion equation yields a Kelvin-Noether circulation theorem for the stochastic system 
\begin{align}
\id\oint_{c(\bx_t)} \bv\cdot d\bx 
= -\oint_{c(\bx_t)}\nabla \id p \cdot d\bx
= 0\,. 
\label{eqn: KelvinSALT}
\end{align}
This stochastic Kelvin circulation theorem is has the same form as that for the deterministic system, except that each line element of the material loop in Kelvin's theorem follows the Stratonovich stochastic Lagrangian path, $\bx_t$.

The real vectors $\bs{\xi}_i$ comprise the time-independent noise amplitudes which are to be determined from data assimilation. The $W^i$ are independent (uni-dimensional) standard Brownian motions and $\circ$ denotes stochastic integration in the Stratonovich sense\footnote{An exposition of Brownian motion, stochastic calculus and the Stratonovich integral is to be found, for example, in the monograph \cite{karatzas2012brownian}.}.
The curl form of the SALT Euler motion equation in \eqref{eqn: velformSALT} is obtained from \eqref{eq:vec_calc_curl_form} and given by
\begin{equation}
    \id \bv - \id \bx_t \times \mathrm{curl}\bv + \nabla (\bv \cdot \id \bx_t) = -\nabla \id p.
\label{eqn: curlform}
\end{equation}
The motion equation \eqref{eqn: curlform} and its curl yielding the SALT vorticity equation implies a formula for the evolution of the helicity of the flow, $\Lambda$, defined as
\begin{equation}
    \Lambda:=\int_{\cal D} \bv \cdot \mathrm{curl}\bv \,d^3x\,.
\label{eqn: helicity}
\end{equation}
Upon applying the divergence theorem, one finds
\begin{equation}
    \id \Lambda = - \int_{\p\cal D} \mathbf{\widehat{n}}\,\cdot
    \Big((\bv \cdot \mathrm{curl}\bv) \id \bx_t + \mathrm{curl}\bv \, \id p
    \Big)\,dS\,.
\label{eqn: helicons}
\end{equation}
Thus, in a periodic 3D domain, or in an infinite 3D domain with asymptotically vanishing boundary conditions, the SALT motion equation in \eqref{eqn: velformSALT} or \eqref{eqn: curlform} preserves the helicity, $\Lambda$, defined in \eqref{eqn: helicity}. However, a glance at the SALT motion equation in \eqref{eqn: curlform} informs us that it will not preserve the kinetic energy, since even with the usual fluid boundary conditions $\mathrm{div}\bv=0$ implies
\begin{equation}
    \tfrac12\id \|\bv\|_{L^2}^2:=\int_{\cal D} \bv \cdot \id \bx_t \times\mathrm{curl}\bv \,d^3x\ne0\,.
\label{eqn: energy}
\end{equation}

\subsubsection{Modeling under the Location Uncertainty Principle}


 The Location Uncertainty principle consists in decomposing  the flow trajectory $\mathbf{x}\colon\Omega\times \R^{+}\rightarrow\Omega$ over a bounded domain, $\Omega\subset\R^{3}$
\begin{equation}\label{eq:Displacement_def}
\mathrm{d}\mathbf{x}_{t} = \mathbf{v}\left(\mathbf{x}_{t},t\right)\,\mathrm{d} t + \mbs{\sigma}\left(\mathbf{x}_{t},t\right)\mathrm{d}\mathbf{W}_{t}
\end{equation}
in terms of $\bv\left(\mathbf{x}_{t},t\right)$, a smooth-in-time component of the (Lagrangian) velocity and a noise $\mbs{\sigma}\left(\mathbf{x}_{t},t\right)\mathrm{d}\mathbf{W}_{t}$, which has here to be understood in the It\^{o} sense and that accounts for the unresolved processes. The Wiener process, $\mathbf{W}_{t}$  is a $(H-$valued (cylindrical) Brownian motion, where H is the Hilbert space of square integrable functions. The noise is then properly defined as the application of an Hilbert-Schmidt symmetric integral kernel $\mbs{\sigma}_t\mbs{f}\left(\mbs{x}\right)=\int_{{\cal S}}\breve{\sigma}\left(\mbs{x},\mbs{y},t\right)\mbs{f}\left(\mbs{y}\right)\,\mathrm{d}\mbs{y}$ to the $H-$valued cylindrical Wiener process $\mathbf{W}$ as
\begin{equation}\label{eq:def_sigma}
    \left(\mbs{\sigma}_t\mathrm{d}\mathbf{W}_{t}\right)^{i} \left(\mbs x\right)= \int_{{\cal S}} \breve{\sigma}_{ik}\left(\mbs x,\mbs{y},t\right) \mathrm{d} W_{t}^{k}\left(\mbs{y}\right) \mathrm{d}\mbs{y},
\end{equation}
The role of the integrable kernel $\breve{\sigma}$ is to impose a spatial correlation on the small-scale component. It leads to the covariance tensor $\mbs{Q}$ 
\begin{eqnarray}
    Q_{ij}\left(\mbs{x},\mbs{y},t,s\right) &=& {\rm I\!E}\left[ \left(\mbs{\sigma}_t \mathrm{d}\mathbf{W}_{t}\left(\mbs{x}\right)\right)^{i}  \left(\mbs{\sigma}_t\mathrm{d}\mathbf{W}_{s}\left(\mbs{y}\right) \right) ^{j} \right] \nonumber \\
    &=& \delta\left(t-s\right)\mathrm{d}t \int_{{\cal S}} \breve{\sigma}_{ik}\left(\mbs{x},\mbs{z},t\right) \breve{\sigma}_{kj}\left(\mbs{z},\mbs{y},s\right) \,\mathrm{d}\mbs{z}, \nonumber 
\end{eqnarray}
 of the centered Gaussian process $\mbs{\sigma}_{t}\mathrm{d}\mathbf{W}_{t} \sim \mathcal{N}\left( 0,\mathbf{Q}\mathrm{d}t \right)$. The diagonal components of the covariance tensor per unit of time, referred to as the variance tensor, $\mathbf{a}$, is a positive definite matrix defined as $\mathbf{a} (\mathbf{x}, t)\delta(t-t')\mathrm{d} t = \mathbf{Q}(\mathbf{x},\mathbf{x},t,t')$, that quantifies the strength of the noise and has the dimension of a viscosity in $\rm{m}^2\rm{s}^{-1}$. The operator $Q$ being compact auto-adjoint positive definite operator on $H$, it admits eigenfunctions $\mbs{\xi}_{n}\left(\mbs{\cdot},t\right)$ with (strictly) positive eigenvalues $\lambda_{n}\left(t\right)$ satisfying $\sum_{n\in\N}\lambda_{n}\left(t\right)<+\infty$. As a consequence, the noise and the variance tensor $\mbs{a}$ can be expressed  through the spectral representation
\begin{eqnarray}
    \mbs{\sigma}_t\mathrm{d}\mathbf{W}_{t} \left(\mbs{x} \right) && = \sum_{n\in\N}\lambda_n^{1/2}\left(t\right)\mbs{\xi}_{n}\left(\mbs{x},t\right)\mathrm{d}\beta_{n} \\
    \mbs{a}\left(\mbs{x},t \right) && = \sum_{n\in\N}\lambda_n\left(t\right) \mbs{\xi}_{n}\left(\mbs{x},t\right)\mbs{\xi}_{n}^{\dag}\left(\mbs{x},t\right).
\end{eqnarray}
 The rate of change of a volume $V_{t}$ of the scalar $q$ is given by  the stochastic Reynolds transport theorem, introduced in \cite{em2014}
\begin{equation}
\label{eq:stoch_REYNOLDS}
\mathrm{d}\int_{V_{t}}q\left( \mbs{x},t \right)\,\mathrm{d}\mbs{x} = \int_{V_{t}} \left\lbrace\mathrm{D}_{t}q+ q\nabla\mbs{\cdot}\left[ \mathbf{v}^{\star}\,\mathrm{d}t + \mbs{\sigma}_t\mathrm{d}\mathbf{W}_{t} \right]\right\rbrace\left( \mbs{x},t \right) \,\mathrm{d}\mbs{x},
\end{equation}
with the transport operator
\begin{equation}
\label{eq:transport-oper}
	\mathrm{D}_{t}q = \mathrm{d}_{t}q + \left[ \mathbf{v}^{\star} \,\mathrm{d}t 
 	+ \mbs{\sigma}_t\,\mathrm{d}\mathbf{W}_{t}
	\right] \mbs{\cdot}\nabla q
	- \frac{1}{2}\nabla\mbs{\cdot}\left(\mathbf{a}\nabla q\right) \,\mathrm{d}t.
\end{equation}
In this formula, the first component of the right-hand side is the increment in time at a fixed location of the process $q$, that is $\mathrm{d}_{t}q=q\left(\mathbf{x}_{t},t+\mathrm{d}t\right) -q\left(\mathbf{x}_{t},t\right)$, playing the role of a derivative in time for a non differentiable process. The effective velocity $\mathbf{v}^{\star}$ is defined as  
\begin{equation}\label{eq:modified_adv}
    \mathbf{v}^{\star} = \mathbf{v}
 	- \frac{1}{2}\nabla\mbs{\cdot}\,\mathbf{a}
 	+ \mbs{\sigma}_t^{\ast}\left(\nabla\mbs{\cdot}\,\mbs{\sigma}_t\right),
\end{equation}
where the  velocity component $\mathbf{v}_{s} = \frac{1}{2}\nabla\mbs{\cdot}\,\mathbf{a}$ results from the noise inhomogeneities. For incompressible homogeneous noise as considered in this work $\mathbf{v}^{\star} = \mathbf{v}$. Besides, the diffusion term exactly balances the noise brought by the noise. With Stratonovich convention and a homogeneous noise the transport operator takes a simplified form similar to the material derivative:
\begin{equation}
\label{eq:transport-oper-S}
	\mathrm{D}_{t}q = \mathrm{d}_{t}q + (\bv \,\mathrm{d}t 
 	+ \mbs{\sigma}_t\,\circ\mathrm{d}\mathbf{W}_{t}) \mbs{\cdot}\nabla q.
\end{equation}
For a divergence free homogeneous noise, the Euler equation, in the LU framework, can then be defined as:
\begin{equation}
    \mathrm{d}_{t}\bv + (\bv \,\mathrm{d}t 
 	+ \mbs{\sigma}_t\,\circ\mathrm{d}\mathbf{W}_{t}) \mbs{\cdot}\nabla \bv = -\nabla dp_t, \;\;\nabla\mbs{\cdot} \bv =0,
\end{equation}
where $dp_t$ denotes the pressure composed of a finite variation term and a martingale pressure term. With the Leray projection, $\mathbb{P}$, this pressure term can be removed and we obtain, the inertial form of the Euler equation:
\begin{equation}
    \mathrm{d}_{t}\bv + \mathbb{P}\bigl((\mathrm{d} \mathbf{x}_{t}\mbs{\cdot}\nabla) \bv\bigr) = 0, \quad \;\nabla\mbs{\cdot} \bv =0.
\end{equation}

\subsection{Triad Model Comparison}

The reduced order model for the full-scale 3D SALT Euler and 3D LU Euler  for a single triad interaction equation is  obtained by projecting the continuous stochastic Euler models onto the helical modes, in the same fashion as for the deterministic equation \eqref{eqn: complexRB}.
Therefore, we introduce an additional Stratonovich stochastic term into the transport velocity in \eqref{eqn:u_amplitudes} as
\begin{equation}
\id \bx_t(\bk,t) := \big(a_+(\bk,t)\bh_+ + a_-(\bk,t)\bh_-\big)\dt  + \sum_i \big(b_+^i(\bk)\bh_+ + b_-^i(\bk)\bh_-\big)\circ dW_t^i\,,
\label{eqn: SALTamp}
\end{equation}
where the $\bk$-dependent complex vector $\bb(\bk):=(b_{s_k},b_{s_p},b_{s_q})^T\in \mbb{C}^3$ represents the time-independent 
noise amplitude which is to be determined from data assimilation, similar to the continuous stochastic models~\eqref{eqn: velformSALT}.
Enumerating the equation for a single triad then yields, after
rearranging using exchange symmetry in $(\bk,\bp,\bq)$, the matrix equation
\begin{align}
\begin{split}
\id
\begin{bmatrix}
{a}_{s_k} \\ 	{a}_{s_p}  \\	 {a}_{s_q} 	
\end{bmatrix}
= g
\begin{bmatrix}
0  & 	-qs_q{a}_{s_q}        &	ps_p {a}_{s_p} 	\\
qs_q {a}_{s_q}     &	  0   & -ks_k {a}_{s_k}  \\
-ps_p {a}_{s_p}    &         ks_k {a}_{s_k}  &	 0
\end{bmatrix}^*
\begin{bmatrix}
{a}_{s_k}\dt \ +\  b_{s_k} \circ dW_t\\ 	{a}_{s_p}\dt \ +\ b_{s_p} \circ dW_t  \\	 {a}_{s_q} 	\dt \ +\ b_{s_q} \circ dW_t
\end{bmatrix}^*
\,.
\end{split}
\label{eqn: 3waves-scalarSALT}
\end{align}
Upon applying the previous steps for the deterministic case to the stochastic velocity in \eqref{eqn: SALTamp}, the single triad interaction 
dynamics for the SALT case would emerge as, cf. equation \eqref{eqn: complexRB},
\begin{equation}
\id\ba =  g \Big(\ba(\bk,t)dt + \bb(\bk)\circ dW_t\Big)^*  \times \mbb{D} \ba^*\,.
\label{eqn: SALTcomplexRB}
\end{equation}
The details of this computation can be found in the Appendix~\ref{adx:luEuler} for the 3D LU Euler model and in Appendix~\ref{adx:saltEuler} for the 3D SALT Euler model.

Remarkably, the HST equation for triad interaction~\eqref{eqn: SALTcomplexRB} still preserves the triad helicity $\ba\cdot\mbb{D}\ba^*$. Hence we name this model the \emph{helicity preserving stochastic triad} (HST) model.
Note that in both equations we use as single source of noise  (One Brownian motion drives the entire system).

It is readily checked that the HST triad evolution~\eqref{eqn: SALTcomplexRB} preserves the helicity. Let's have a look at the diffusion coefficients.
\begin{equation}
    \ba^* \times \mbb{D} \bb = 
    \begin{bmatrix}
        a_p^* s_q q b_q - a_q^* s_p p b_p \\
        a_q^* s_k k b_k - a_k^* s_q q b_q \\
        a_k^* s_p p b_p - a_p^* s_k k b_k
    \end{bmatrix}
    \,,\;\;
    \bb  \times \mbb{D}\ba^* =
    \begin{bmatrix}
        b_p s_q q a_q^* - b_q s_p p a_p^* \\
        b_q s_k k a_k^* - b_k s_q q a_q^* \\
        b_k s_p p a_p^* - b_p s_k k a_k^*
    \end{bmatrix}
    .
\end{equation}
Taking the difference
\begin{equation}
\ba^* \times \mbb{D} \bb - \bb  \times \mbb{D}\ba^* =
    \begin{bmatrix}
        a_p^* (s_q q + s_p p) b_q - a_q^* (s_p p + s_q q) b_p \\
        a_q^* (s_k k + s_q q) b_k - a_k^* (s_q q + s_k k) b_q \\
        a_k^* (s_p p + s_k k) b_p - a_p^* (s_k k + s_p p) b_k
    \end{bmatrix}.
\end{equation}
Writing $\rho:=\operatorname{Tr}{\mbb{D}}$ we get
\begin{equation}
\ba^* \times \mbb{D} \bb - \bb  \times \mbb{D}\ba^* =
    \begin{bmatrix}
        a_p^* (\rho-s_k k) b_q - a_q^* (\rho-s_k k) b_p \\
        a_q^* (\rho - s_p p) b_k - a_k^* (\rho- s_p p) b_q \\
        a_k^* (\rho - s_q q) b_p - a_p^* (\rho-s_q q) b_k
    \end{bmatrix}.
\end{equation}
So that the difference term becomes
\begin{equation}
    \ba^* \times \mbb{D} \bb - \bb  \times \mbb{D}\ba^* = (\rho \operatorname{Id} - \mbb{D})(\ba \times \bb)^*.
\end{equation}

Since the projections of the LU and SALT models onto a single triad are indistinguishable, we introduce a different model that  conserves energy on a single triad to enable a comparison between energy conserving and helicity conserving models. 

\paragraph{Energy-preserving stochastic triad (EST) model.}
We introduce below a modified version of the HST triad equation \eqref{eqn: SALTcomplexRB} that introduces stochasticity into the vorticity instead of into the transport velocity and thereby conserves the energy. This is inspired by the full-scale model introduced in  \cite{holm2021stochastic}. The reduced model is as follows
\begin{equation}
\id\ba =  -\,{g}   \ba^* \times \mbb{D}\Big(\ba(\bk,t)\dt + \bb(\bk)\circ dW_t\Big)^*\,.
\label{eqn: erg}
\end{equation}
We call this model the energy preserving stochastic triad (EST). Written in matrix form the equation 
\eqref{eqn: erg} becomes

\begin{align}
\begin{split}
\id
\begin{bmatrix}
{a}_{s_k} \\ 	{a}_{s_p}  \\	 {a}_{s_q} 	
\end{bmatrix}
= {g}
\begin{bmatrix}
0  & 	-{a}_{s_q}        &	{a}_{s_p} 	\\
{a}_{s_q}     &	  0   & -{a}_{s_k}  \\
-{a}_{s_p}    &         {a}_{s_k}  &	 0
\end{bmatrix}^*
\begin{bmatrix}
ka_k\big({a}_{s_k}\dt \ +\  b_{s_k}(\bk) \circ dW_t\big) \\  ps_p\big({a}_{s_p}\dt \ +\ b_{s_p}(\bp) \circ dW_t \big)
\\	 qs_q\big({a}_{s_q} \dt \ +\ b_{s_q}(\bq) \circ dW_t\big)
\end{bmatrix}^*
\,.
\end{split}
\label{eq:LU_Triad}
\end{align}

The exchange symmetry between the two models HST and EST in the placement of the noise in equations \eqref{eqn: SALTcomplexRB} and \eqref{eqn: erg} is apparent already in the exchange symmetry between velocity and vorticity in Euler's fluid equations \eqref{eqn: Euler}.

\paragraph{Deviation from the conservation laws.}
We can write the equations for the deviation from the conservation laws, which is present in both models.
The SALT model deviates from the energy conservation by
\begin{equation}
    \id_t E_{\text{HST}} = g \bb\cdot(\mbb{D}\ba^*\times \ba^*) \circ \id W_t
\end{equation}
whereas the LU model deviates from the helicity conservation by
\begin{equation}
    \id_t {H}_{\text{EST}} = g \mbb{D}\bb\cdot (\mbb{D}\ba^*\times \ba^*) \circ \id W_t.
\end{equation}
This is seen, since, to get the energy we dot the HST equation with $\ba^*$ and to get helicity we dot the EST equation with $\mbb{D}\ba^*$ and use the standard identities
\begin{align}
    &\ba^* \cdot (\bb  \times \mbb{D}\ba^*) = \bb \cdot (\mbb{D}\ba^*\times \ba^*) \\ 
    &\mbb{D}\ba^* \cdot ( \ba^* \times \mbb{D} \bb ) = \mbb{D}\bb\cdot(\mbb{D}\ba^*\times \ba^*).
\end{align}
Therefore, $\bb$ respects the right scaling and no further scale adjustments between the SALT and LU noise scaling need to be performed in order to compare the models.

\section{Data Assimilation Comparison}\label{sect:dataassimilation}

In this section, we perform a comparative study of the two reduced order models (HST and EST) introduced above by using data assimilation tools. The particular methodology that we make use of is that of particle filters. We will first briefly  explain the particle filtering methodology in a generic framework:

Let $X$ and $Z$ be two processes defined on a given  probability space $(\Omega, \mathcal{F}, \mathbb{P})$. The process $X$ is usually called the \textit{signal process} or the \textit{truth} and $Z$ is the \textit{observation process}. In this paper, $X$ is the pathwise solution of a (deterministic) shell model. The pair of processes $(X,Z)$ forms the basis of the nonlinear filtering problem which consists in finding the best approximation of the posterior distribution of the signal $X_t$ given the observations $Z_1, Z_2, \ldots, Z_{t_{n}}$ \footnote{For a mathematical introduction on the subject, see e.g. \cite{baincrisan}}. The posterior distribution of the signal at time $t$ is denoted by $\pi_t$. We let $d_X$ be the dimension of the state space and $d_Z$ be the dimension of the observation space.  This mixed continuous-discrete time framework can be embedded
into a fully discrete framework, whereby one is interested in computing the
conditional probability law of the signal at the time corresponding to
the observation time. in other words one wants to compute the conditional
distribution $\pi _{{t_n}}$ of $X(t_{n})$ given the data $Z(t_{1}),Z(t_{2}),\ldots,Z(t_{n})$. The process $X$ is assumed to be a Markov process, and we will denote by $\mathcal{K}_n$ its transition kernel, that is  
 \begin{equation} \mathcal{K}_n:\mathbb{R}^{d_X}\times\mathcal{B}(\mathbb{R}^{d_X}) \rightarrow [0,1], \ \mathcal{K}_n(x, B) = \mathbb{P}(X_{t_n} \in B|X_{t_{n-1}}=x)
 \end{equation} 
 for any Borel measurable set $B \in \mathcal{B}(\mathbb{R}^{d_X})$ and $x\in \mathbb{R}^{d_X}$. The process $Z$ models noisy measurements of the truth, using the so-called \textit{observation operator} $\mathscr{H}: \mathbb{R}^{d_X} \rightarrow \mathbb{R}^{d_Z}$:
 \begin{equation} 
 Z_{n} = \mathscr{H}(X_{t_n})+ V_n
 \end{equation}
  where {$\{V_n\}_{n\geq 0}$} are independent identically distributed random variables that represent the measurement noise and $\mathscr{H}$ is a Borel-measurable function. In this paper we will assume that {$\{V_n\}_{n\geq 0}$} have standard normal distributions, but the same methodology can be applied to more general distributions. Observations are incorporated into the system at \textit{assimilation times}. The following recursion formula holds (see \cite{baincrisan})
 \begin{equation}\label{recursionformula}
 \pi_n = g_n \star \pi_{n-1}\mathcal{K}_n
 \end{equation}
  where by '$\star$' we denoted the projective product (see e.g. Definition 10.4 in \cite{baincrisan}).
    
  In the following, we compare approximations of the posterior distribution of the signal using \textit{particle filters}. These are sequential Monte Carlo methods which generate approximations of the posterior distribution $\pi_t$ using sets of \textit{particles}. That is, they generate approximations that are (random) measures of the form
\begin{equation*}
\pi_n \approx \displaystyle\sum_{\ell} \mathrm{w}_{n}^{\ell}\delta({x_{n}^{\ell}}),
\end{equation*}
where $\delta$ is the Dirac delta distribution, $\mathrm{w}_{t}^1, \mathrm{w}_{t}^2, \ldots $ are the \textit{weights} of the  particles and $x_{t}^1, x_{t}^2, \ldots$ are their corresponding positions. Particle filters are used to make inferences about the signal process by using Bayes' theorem, the time-evolution induced by the signal $X$, and the observation process $Z$.

In a standard particle filter, the particles evolve between assimilation times according to the law of the signal.  As we explain below, at each assimilation time the observation is incorporated into the system through the \textit{likelihood function}:
 \begin{equation}
 g_t^{z_t}:\mathbb{R}^{d_X} \rightarrow \mathbb{R}_+, \ g_t^{z_t}(x) = g_t(z_t - \mathscr{H}(x)) \ \ 
 \text{ such that }\ \ \mathbb{P}(Z_t \in dz_t | X_t = x)=g_t^{z_t}(x)dz_t
 \end{equation} 
and all particles are weighted depending on the likelihood of their corresponding position, given the observation. More precisely, the particle $\ell$ is given the weight $w_n^{l} = g_n^{Z_n}(x_\ell)$. Heuristically, the particle weight measures how close the particle trajectory is to the signal trajectory. A selection procedure is then applied to the set of weighted particles.  
Particles with higher conditional likelihood (guided by the observation) have higher weights and will be multiplied, while those which have small likelihoods will be eliminated. For the basic particle filter, this is done by sampling with replacement from the population of particles, with corresponding probabilities proportional to their weights.
 
A Monte Carlo implementation of the transition kernel of the signal  may not always yield good approximations. In many situations one replaces the original transition kernel with {likelihood informed importance proposals}, leading to much better approximations. One situation when this is necessary is when the original process is actually deterministic (aside for the initial position which is assumed to be random). This is the case in our paper. 

To overcome the collapse of the particle filter when using deterministic transition kernels, one can use a Markov Chain Monte Carlo procedure that leaves the deterministic dynamics invariant. This procedure can be costly and might not always introduce enough spread into the sample. In this paper, we propose a different approach, which we illustrate numerically in Section~\ref{sec:DA_det} below. In particular, we propose two different transition kernels based on the physical conservation properties: 
\begin{itemize}
\item The transition kernel associated with the HST model equation (\ref{eqn: 3waves-scalarSALT}). We will denote this transition kernel by $\mathcal{K}^{\text{HST}}$. As we have explained above, this transition kernel preserves the helicity of the system.  

\item The transition kernel associated with the triad model equation (\ref{eq:LU_Triad}) we will denote this transition kernel by  $\mathcal{K}^{\text{EST}}$. 
As we have explained above, this transition kernel preserves the energy of the system.  

\end{itemize}

\subsection{Numerical studies}\label{subsect:numericalstudies}

\subsubsection{Numerical implementation} 

The models are discretised using the stochastic SSPRK3 scheme which is documented, for example, in~\cite{cotter2018numerically}. In our specific case, for instance, the HST model is discretised as
\begin{align*}
    \bq_1^{n} &= \ba_{n} +  g (\ba_n^*\times \mbb{D}\ba_n^*) \Delta t + g(\bb\times\mbb{D}\ba_n^*)\Delta W \\
    \bq_2^{n} &= (3/4) \ba_{n} + (1/4)(\bq_1^{n}+  g ((\bq_1^{n})^*\times \mbb{D}(\bq_1^{n})^*)  \Delta t + g(\bb\times\mbb{D}(\bq_1^{n})^*)\Delta W) \\
     \ba_{n+1} &= (1/3) \ba_{n} + (2/3)(\bq_2^{n}+  g ((\bq_2^{n})^*\times \mbb{D}(\bq_2^{n})^*)  \Delta t + g(\bb\times\mbb{D}(\bq_2^{n})^*)\Delta W)
\end{align*}
where $\Delta t$ denotes the timestep and $\Delta W$ the increment of the driving Brownian motion. Further, $\ba_n$ is the approximate complex vector amplitude at time $t = n\Delta t$. The EST model is discretised completely analogously. For the numerical simulations we chose the following triad throughout. We set 
\begin{equation}
    \bk=[1,0,0],\;\; \bp=[0,-1,1],\;\; \bq=[-1,1,-1],
\end{equation}
with parities $s_k=1, s_p=-1, s_q=-1$ and the initial value $\ba_0 = \frac{1}{\sqrt{3}}[1,1,1]$. We set the parameter $\boldsymbol{\Gamma}=[1,1,1]$ and used a time stepsize of $\Delta t=0.0005$.

\subsubsection{Data Assimilation for the Deterministic Model}
\label{sec:DA_det}
\begin{figure}

\centering
\begin{subfigure}{0.45\textwidth}
    \centering
    \includegraphics[width=0.75\textwidth]{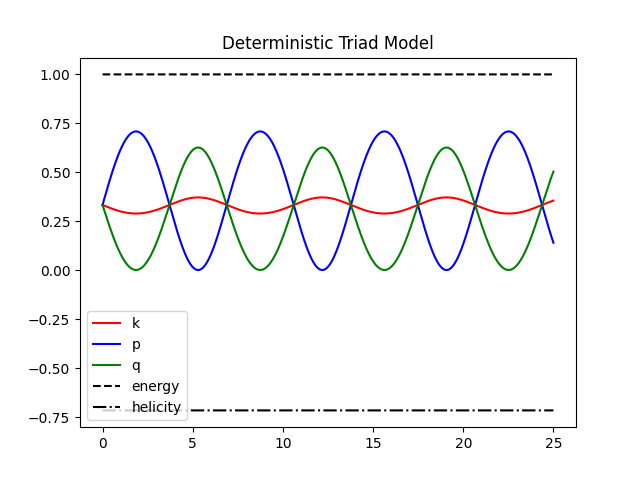}
    \caption{Deterministic Model.}
    \label{fig:det_triad}
\end{subfigure}
\begin{subfigure}{0.45\textwidth}
    \centering
    \includegraphics[width=0.75\textwidth]{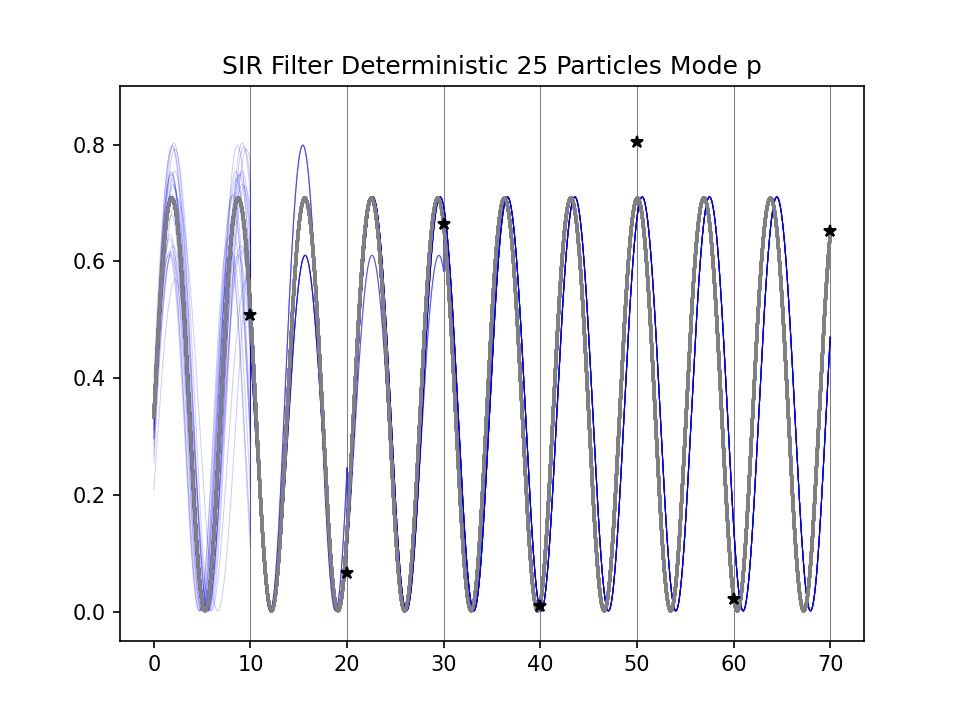}
    \caption{Ensemble for Mode $\bp$.}
    \label{fig:deterministic_filter_p}
\end{subfigure}
\hfill
\begin{subfigure}{0.32\textwidth}
    \includegraphics[width=\textwidth]{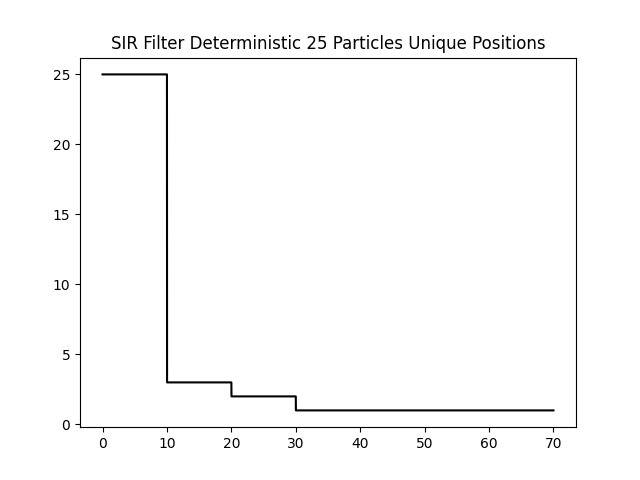}
    \caption{Unique Particles.}
    \label{fig:deterministic_filter_unique}
\end{subfigure}
\hfill
\begin{subfigure}{0.32\textwidth}
    \includegraphics[width=\textwidth]{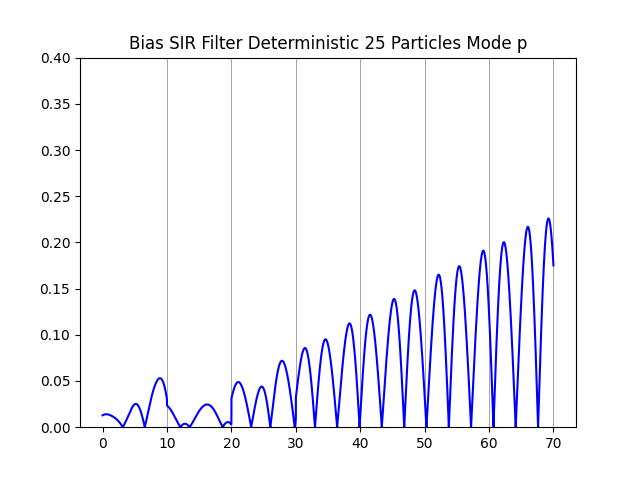}
    \caption{Bias.}
    \label{fig:deterministic_filtering_bias}
\end{subfigure}
\hfill
\begin{subfigure}{0.32\textwidth}
    \includegraphics[width=\textwidth]{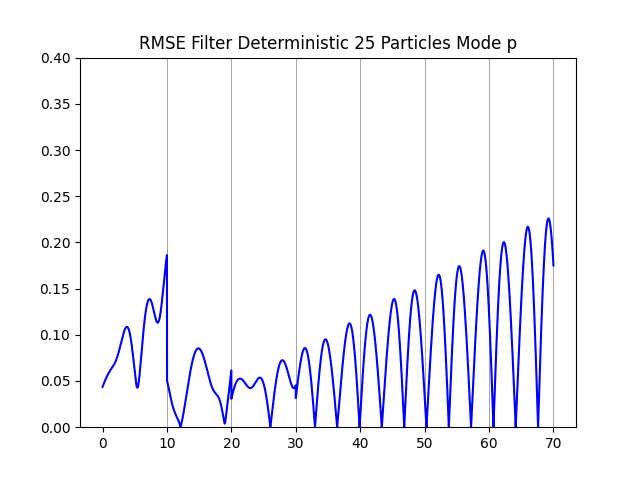}
    \caption{RMSE.}
    \label{fig:deterministic_filtering_rmse}
\end{subfigure}

\caption{Deterministic Triad Model. The horizontal axis shows time. (a) Evolution of the modal energies (colored lines) as well as the total energy (dashed line) and helicity (dash-dotted line). The deterministic model exhibits continually oscillating modal energies. The simulation also confirms the conservation of energy and helicity. (b--e) Data assimilation for the deterministic model using particle filter. (b) Evolution of the energy of mode $\bp$ of the signal (grey line) and evolution of the energy of mode $\bp$ for the particle ensemble (blue lines). Noisy observations (black stars) are made and assimilated every $10$ time units. (c) The number of unique particle positions in the filtering ensemble. (d) The bias of the particle ensemble wrt.~the observations. (e) The RMSE of the particle ensemble wrt.~the observations.} 
\label{fig:deterministic_filter}

\end{figure}

We illustrate the failure of the particle filter with deterministic transition kernel in Figure \ref{fig:deterministic_filter}. In this case, the particle filtering is performed for an ensemble of $n=25$ particles evolving according to the deterministic triad dynamics. The initial ensemble is spread around the initial value $\ba_0$ of the signal according to a Gaussian distribution with standard deviation $1/\sqrt{600}$ and, in particular, does not contain the true initial point. Data assimilation is performed every $10$ time units and the observations are taken from the modal energies of the truth with an observation error $\eta$ distributed as
$
    \eta \sim \mathcal{N}(\mathbf{0}, \mathrm{C})
$
with covariance $\mathrm{C}=\operatorname{diag}(0.005^2, 0.05^2, 0.05^2)$.
We observe that both the bias and the RMSE keep increasing with time to values much larger than the observation error. Moreover, the number of \emph{distinct} particles decreases rapidly: after 30 steps, a single particle remains that is not the true particle since the true particle was not part of the initial cloud. The particle filter does not work.

\subsubsection{Reduced Order Model Realisations}
\label{sec:model_realisation}
The deterministic model in Figure~\ref{fig:det_triad} exhibits continually oscillating triad amplitudes. Plotted are the modal energies. Writing $\ba=[a_k, a_p, a_q]$ we call the real value $a_k^{\phantom{*}}a_k^*$ the \emph{energy of mode $k$}. Similarly for the two other modes.

\begin{figure}

\centering
\begin{subfigure}{0.24\textwidth}
    \includegraphics[width=\textwidth]{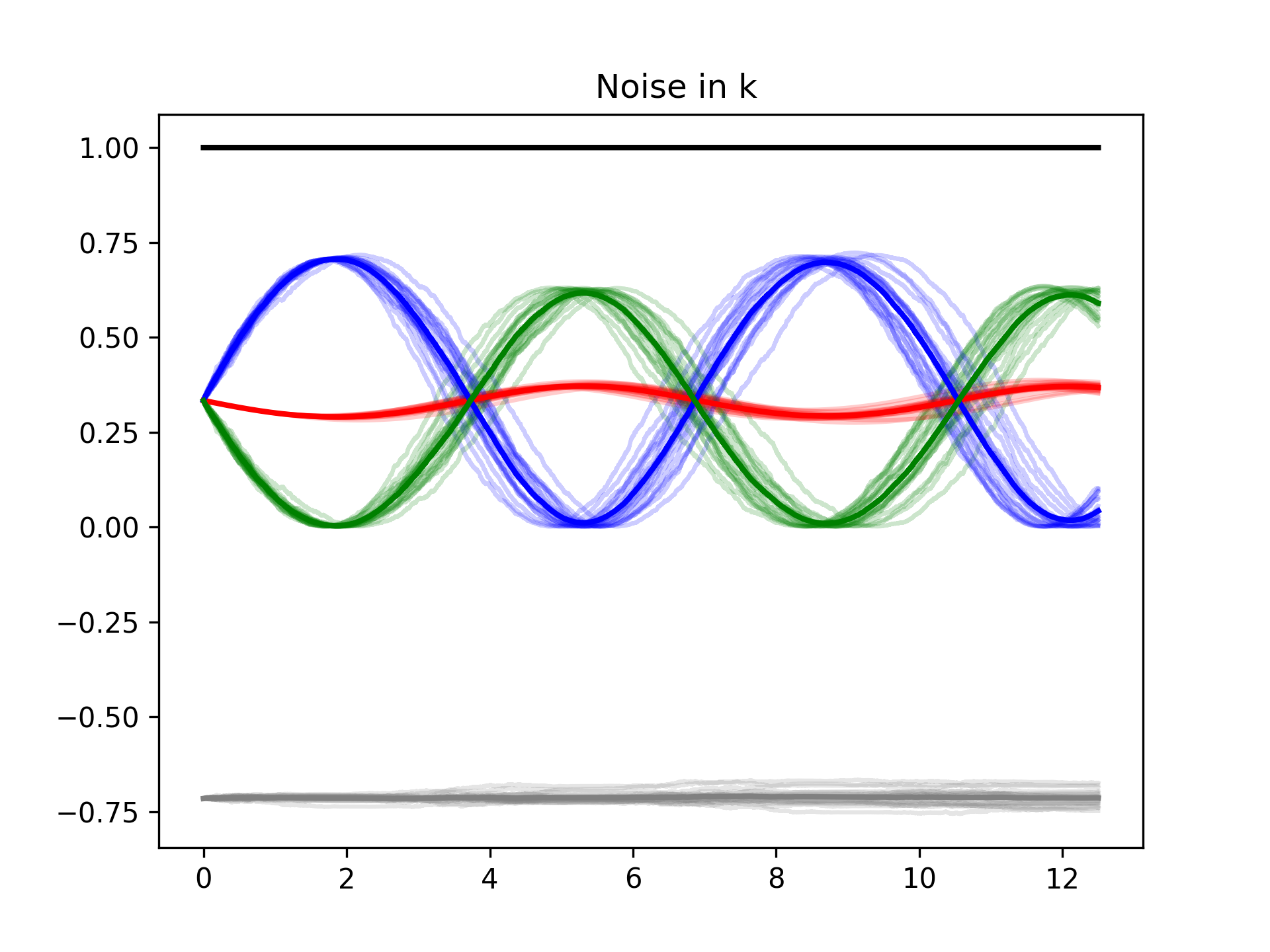}
    \caption{EST $b_k=0.1$.}
    \label{fig:LU_noise_k}
\end{subfigure}
\hfill
\begin{subfigure}{0.24\textwidth}
    \includegraphics[width=\textwidth]{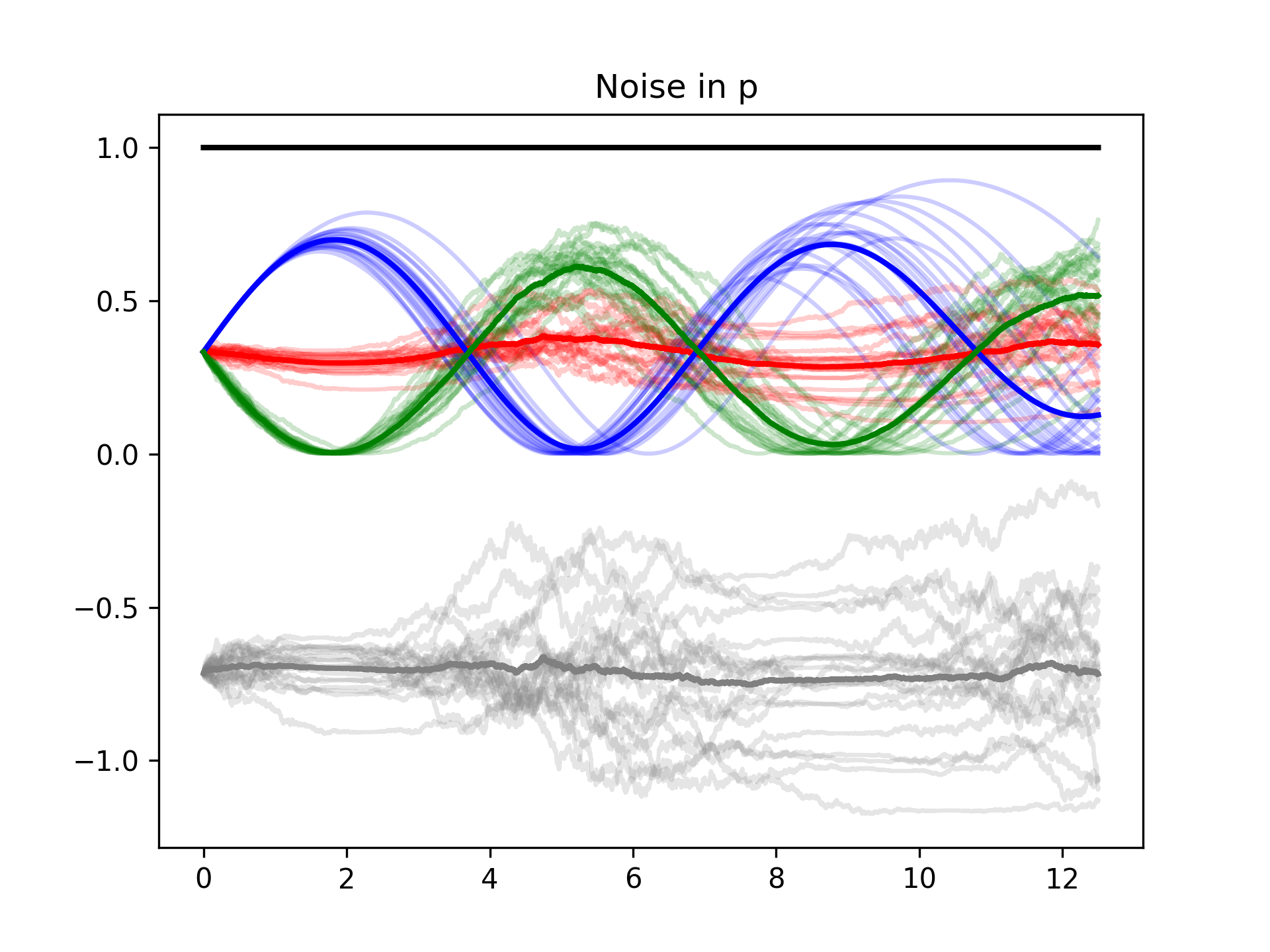}
    \caption{EST $b_p=0.1$.}
    \label{fig:LU_noise_p}
\end{subfigure}
\hfill
\begin{subfigure}{0.24\textwidth}
    \includegraphics[width=\textwidth]{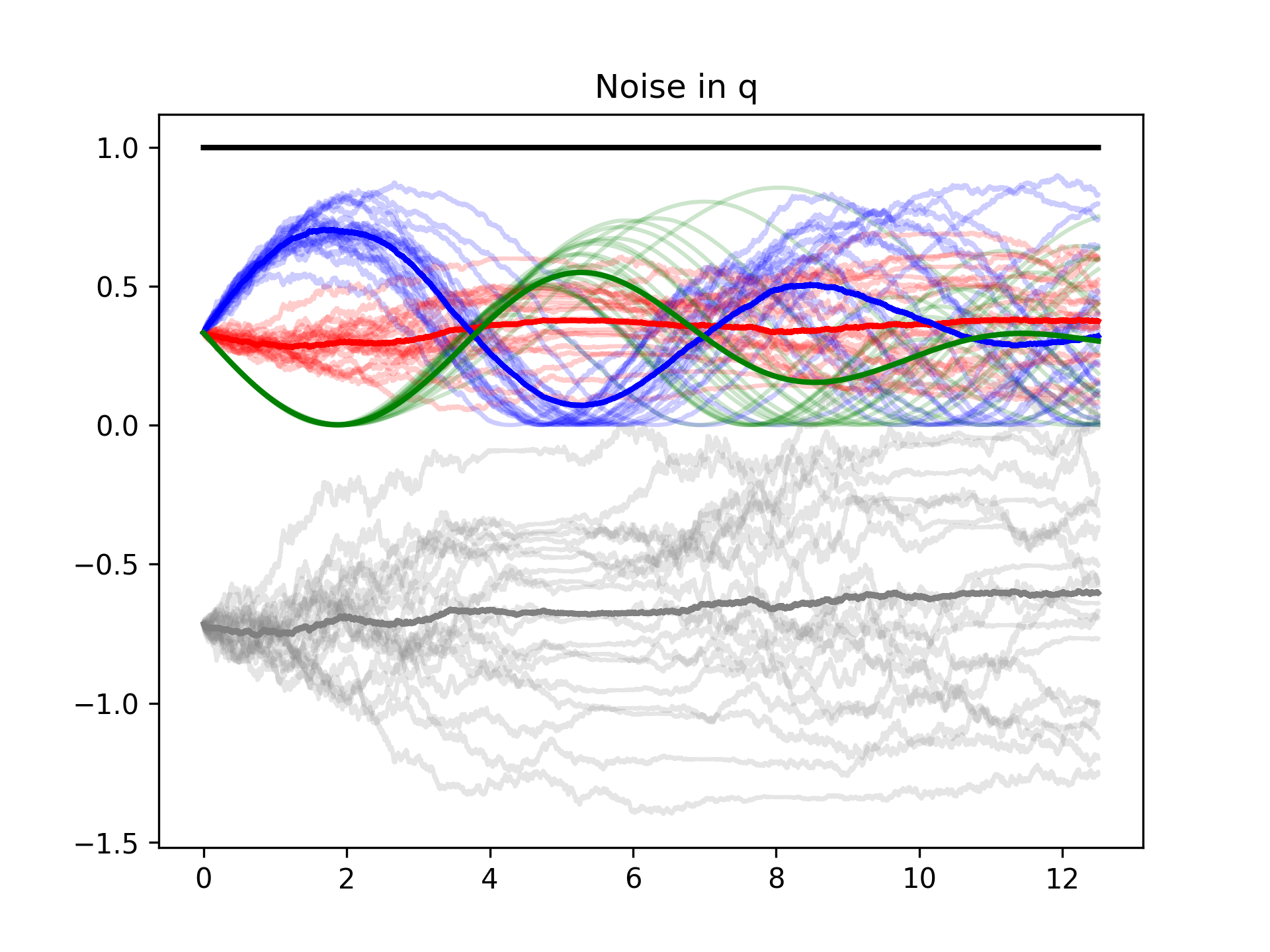}
    \caption{EST $b_q=0.1$.}
    \label{fig:LU_noise_q}
\end{subfigure}
\hfill
\begin{subfigure}{0.24\textwidth}
    \includegraphics[width=\textwidth]{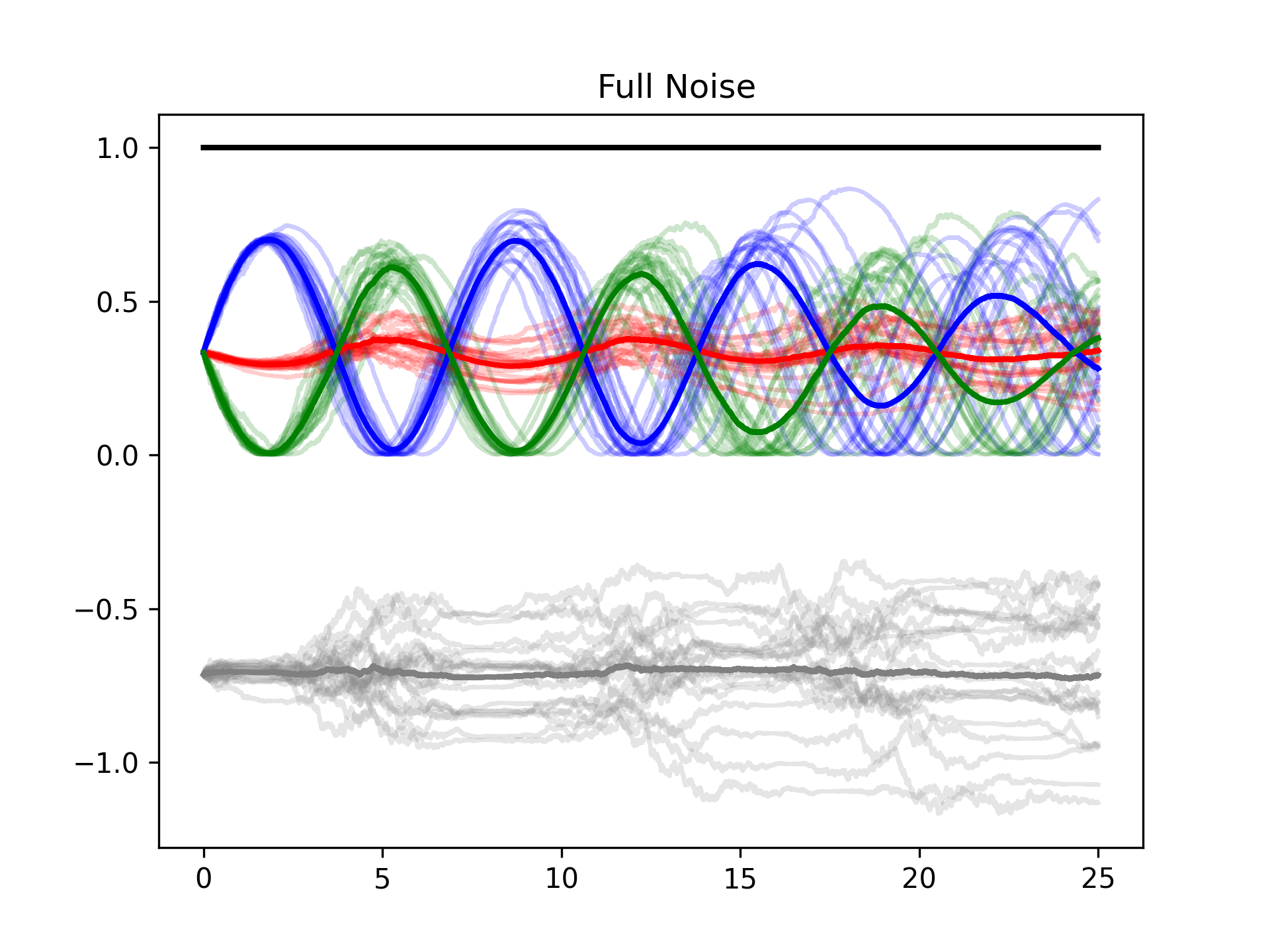}
    \caption{EST full noise.}
    \label{fig:LU_full_noise}
\end{subfigure}
\hfill
\begin{subfigure}{0.24\textwidth}
    \includegraphics[width=\textwidth]{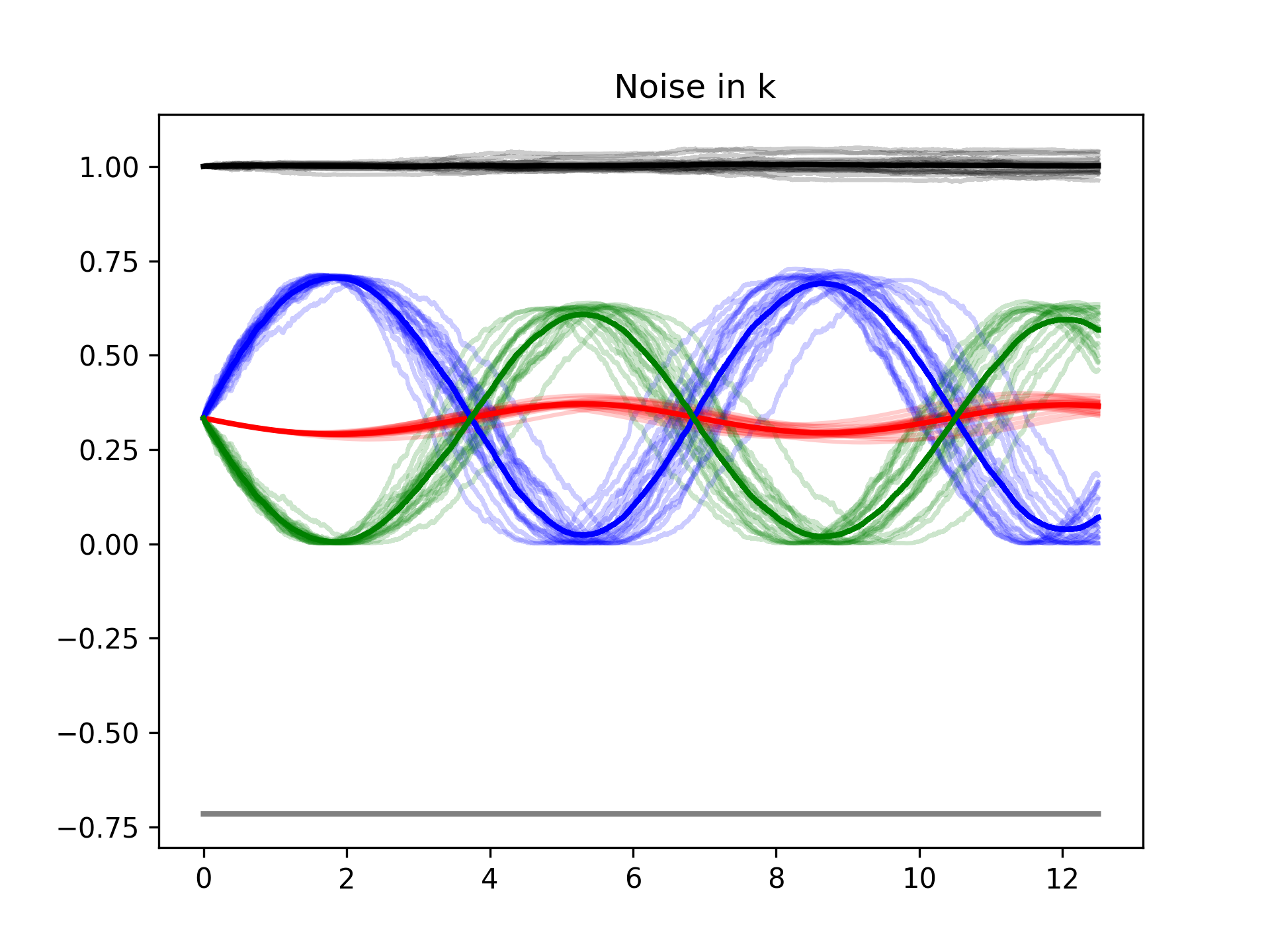}
    \caption{HST $b_k=0.1$.}
    \label{fig:SALT_noise_k}
\end{subfigure}
\hfill
\begin{subfigure}{0.24\textwidth}
    \includegraphics[width=\textwidth]{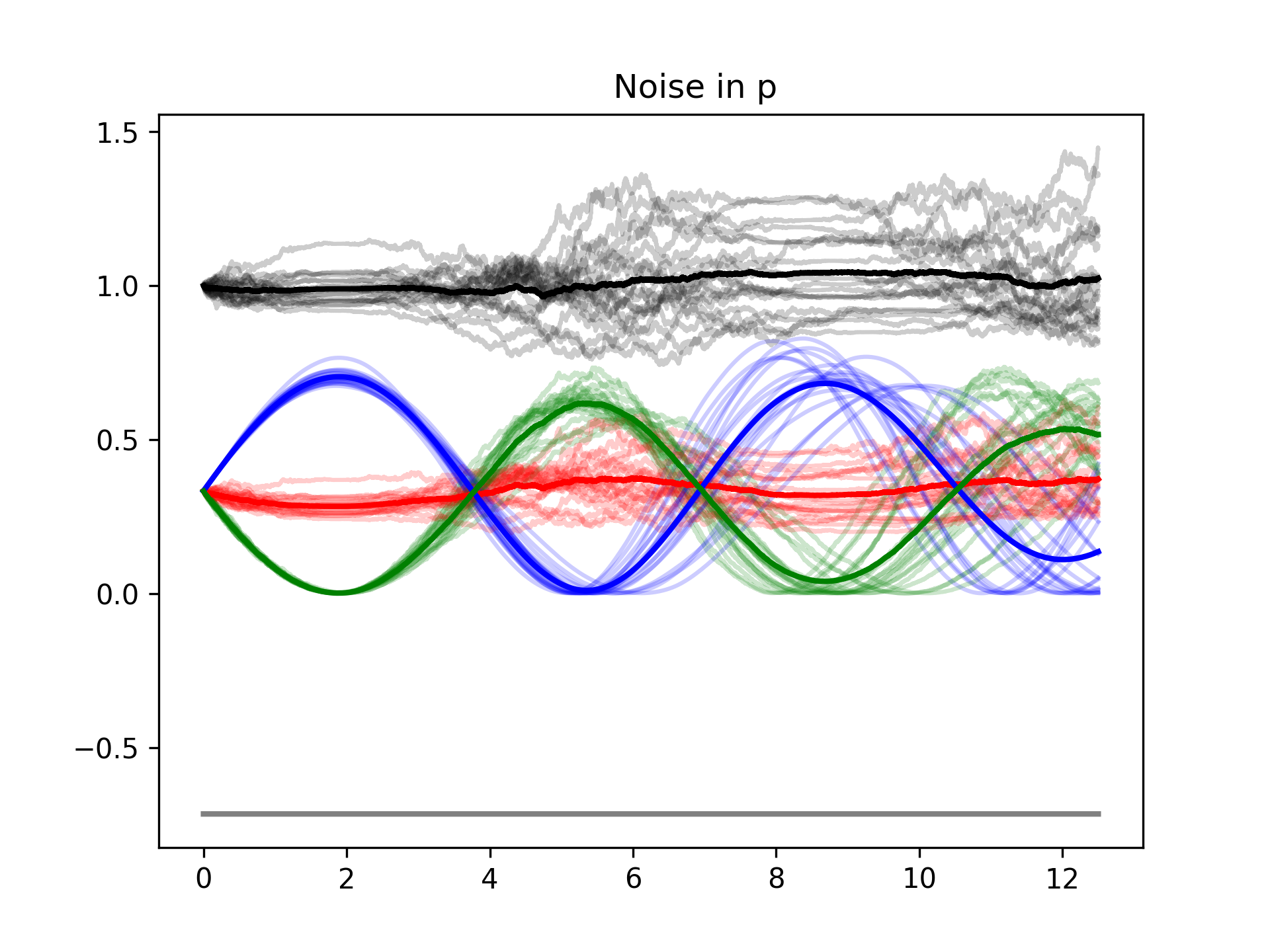}
    \caption{HST $b_p=0.1$.}
    \label{fig:SALT_noise_p}
\end{subfigure}
\hfill
\begin{subfigure}{0.24\textwidth}
    \includegraphics[width=\textwidth]{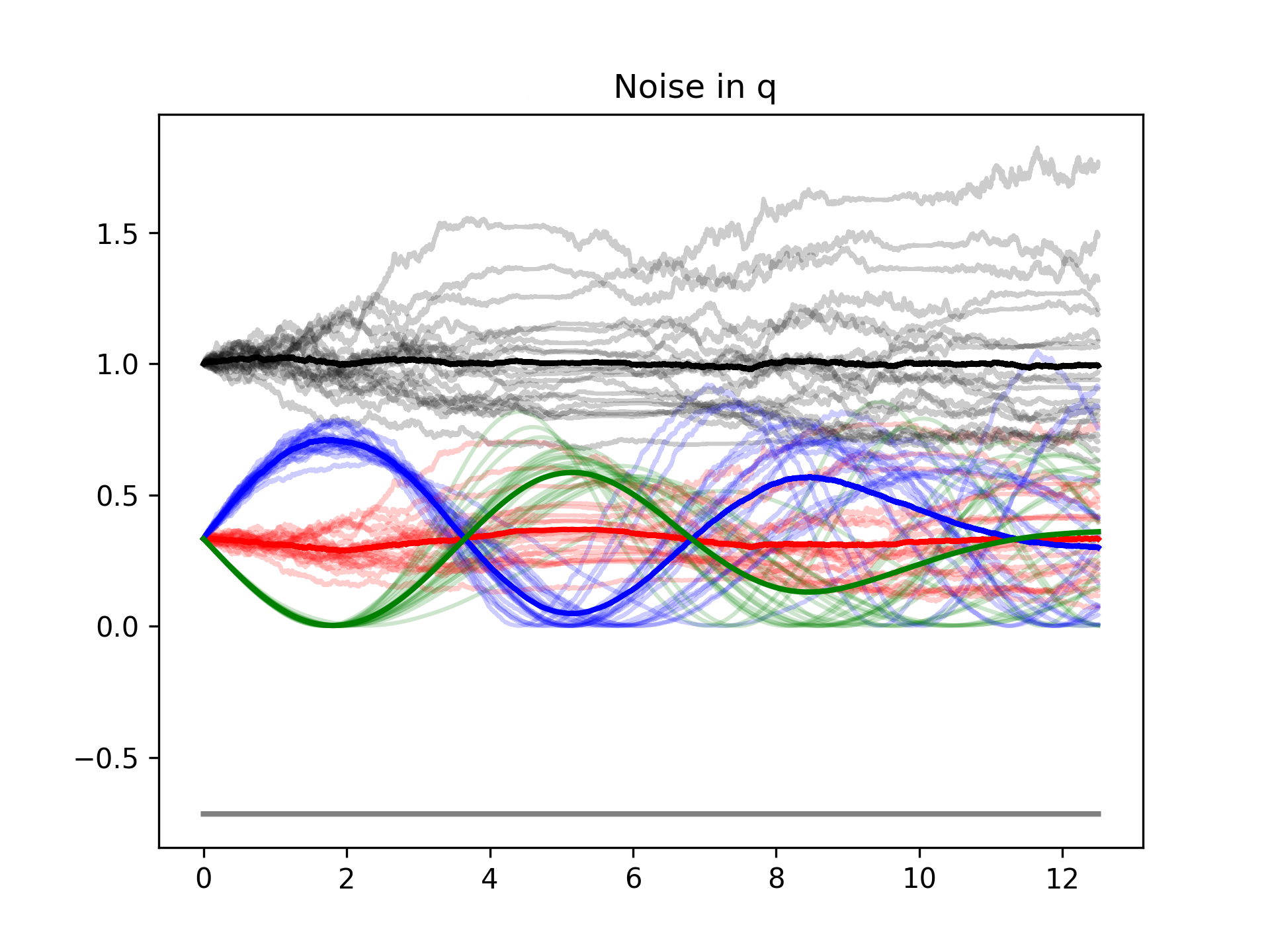}
    \caption{HST $b_q=0.1$.}
    \label{fig:SALT_noise_q}
\end{subfigure}
\hfill
\begin{subfigure}{0.24\textwidth}
    \includegraphics[width=\textwidth]{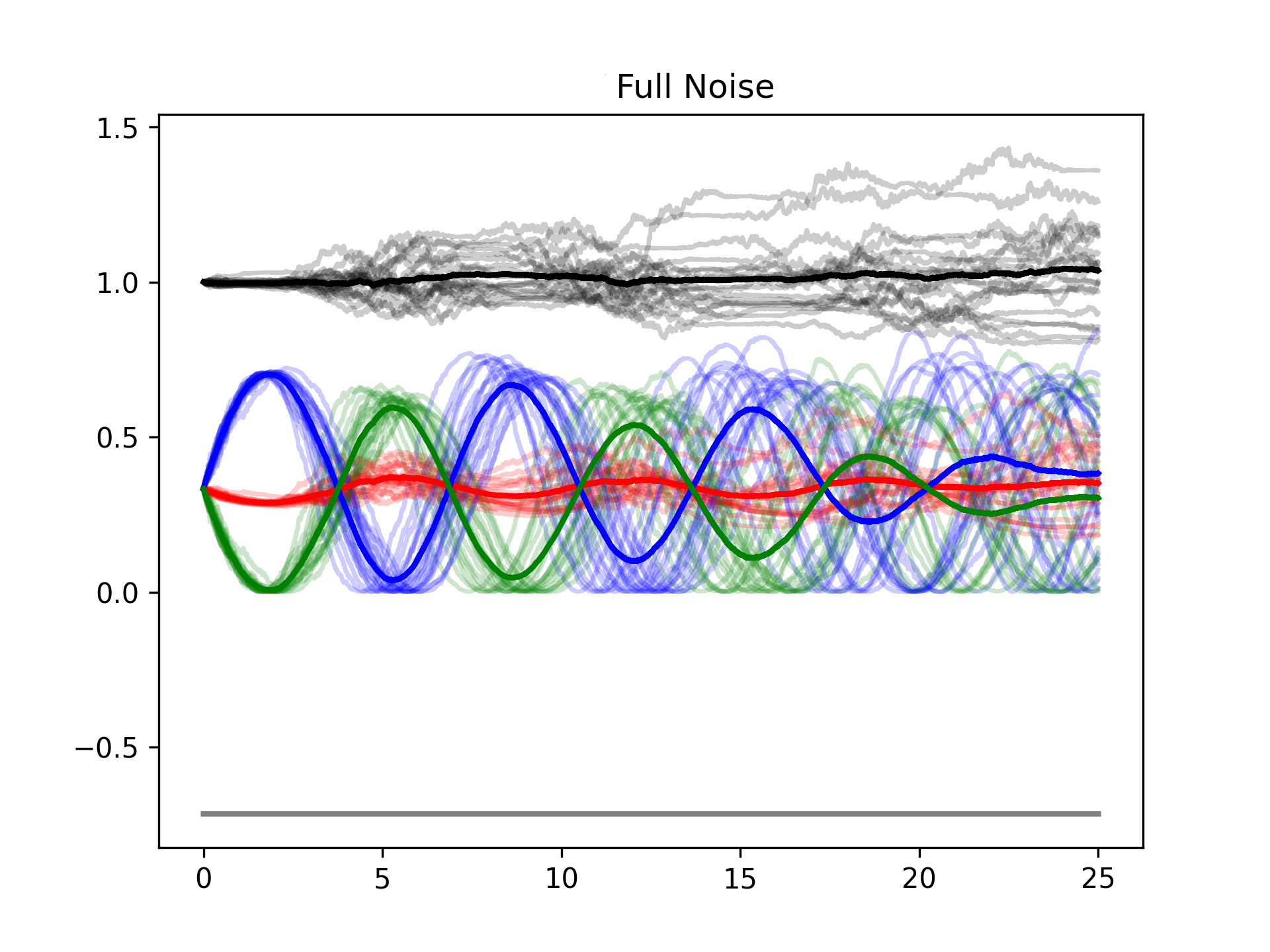}
    \caption{HST full noise.}
    \label{fig:SALT_full_noise}
\end{subfigure}

\caption{Model realisations for both stochastic triad models. Plotted are the modal energies (colored lines), total energy (black line) and helicity (grey line). The respective thick lines are the ensemble means, and the the thin lines represent the different stochastic realisations. (a+e) The noise coefficient $\bb = [0.1,0,0]$. (b+f) The noise coefficient $\bb = [0,0.1,0]$. (c+g) The noise coefficient $\bb = [0,0,0.1]$. (d+h) The noise coefficient $\bb = [0.1,0.05,0.01]$}
\label{fig:single_mode_noise}

\end{figure}

We simulate the model realisations for different noise scenarios. We simulate the effect of noise in each single mode. Let the noise amplitude vector be $\bb=[b_k, b_p, b_q]$. Then we simulate the two models for $\bb=[b_k=0.1, b_p=0, b_q=0]$, $\bb=[b_k=0, b_p=0.1, b_q=0]$, and $\bb=[b_k=0, b_p=0, b_q=0.1]$. The trajectories of the modal energies for $n=20$ realisations of the driving noise for each scenario are shown in Figure~\ref{fig:single_mode_noise}.
We also simulate the case of full noise for the noise amplitude vector $\bb=[b_k=0.1, b_p=0.05, b_q=0.01]$ which was calibrated to the data assimilation objective using the procedure explained in Appendix~\ref{sec:noise_calib}. The ensemble of $n=20$ realisations of the driving noise in the full noise case is depicted in Figures~\ref{fig:LU_full_noise} and~\ref{fig:SALT_full_noise}. In all cases, it can be observed that the mean energy amplitudes of the modes are dampened in both stochastic models. Furthermore, we can experimentally verify the conservation of triad energy for the EST model and the conservation of triad helicity for the HST model.

\subsubsection{Model Statistics}
\label{sec:model_stats}
\begin{figure}

\centering
\begin{subfigure}{0.23\textwidth}
    \includegraphics[width=\textwidth]{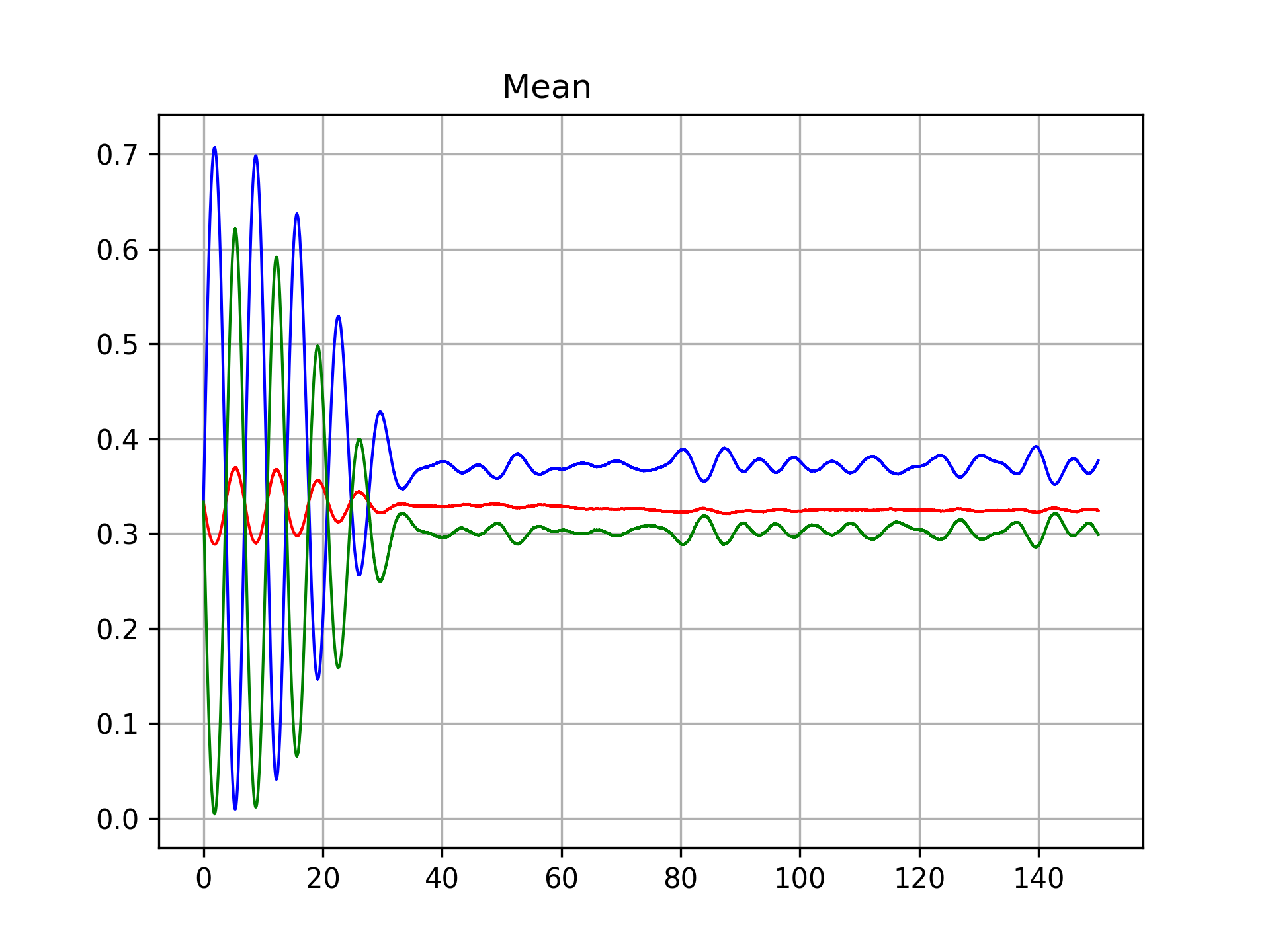}
    \caption{EST mean.}
    \label{fig:LU_mean}
\end{subfigure}
\hfill
\begin{subfigure}{0.23\textwidth}
    \includegraphics[width=\textwidth]{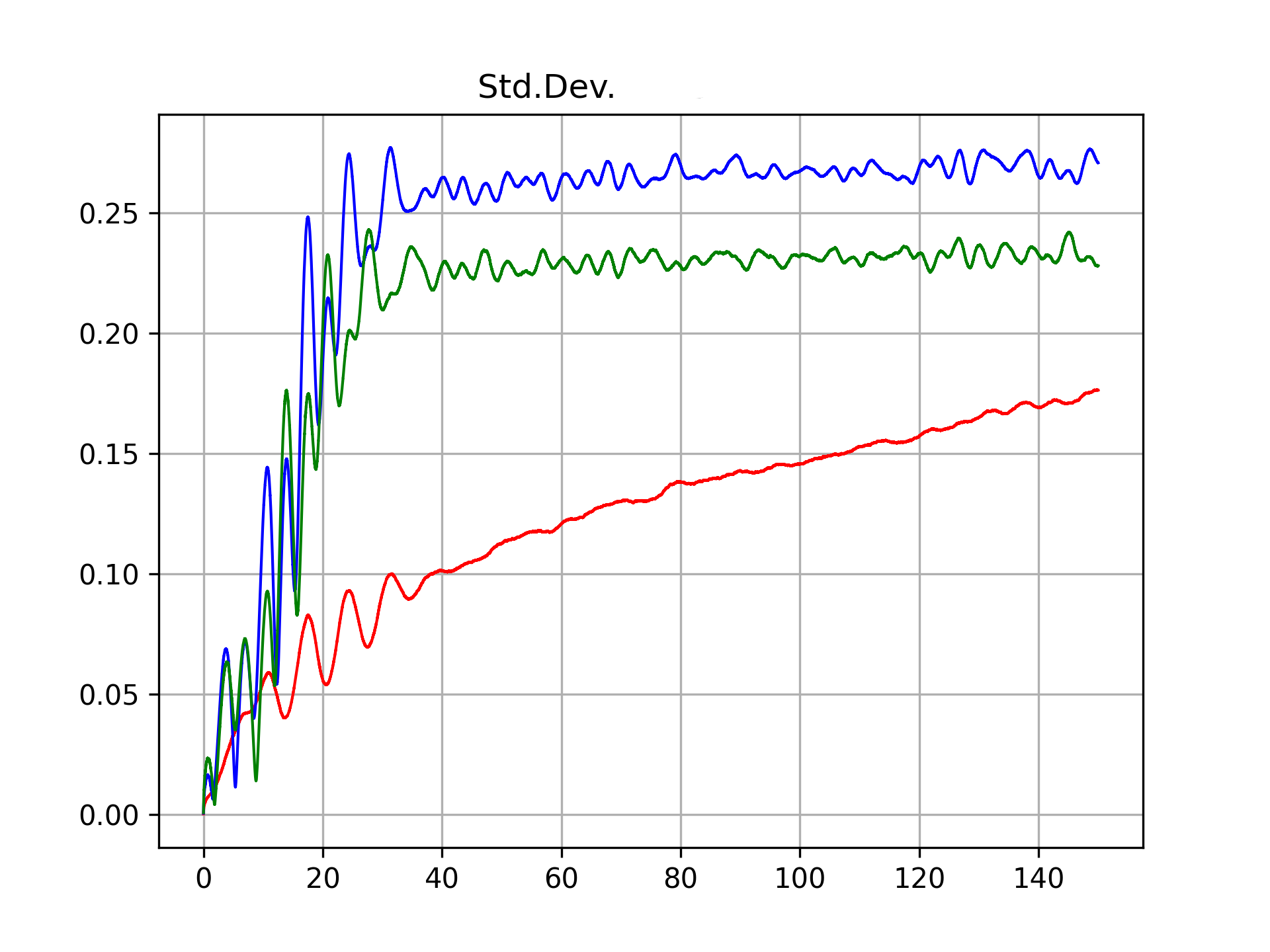}
    \caption{EST std.~dev.}
    \label{fig:LU_std}
\end{subfigure}
\hfill
\begin{subfigure}{0.23\textwidth}
    \includegraphics[width=\textwidth]{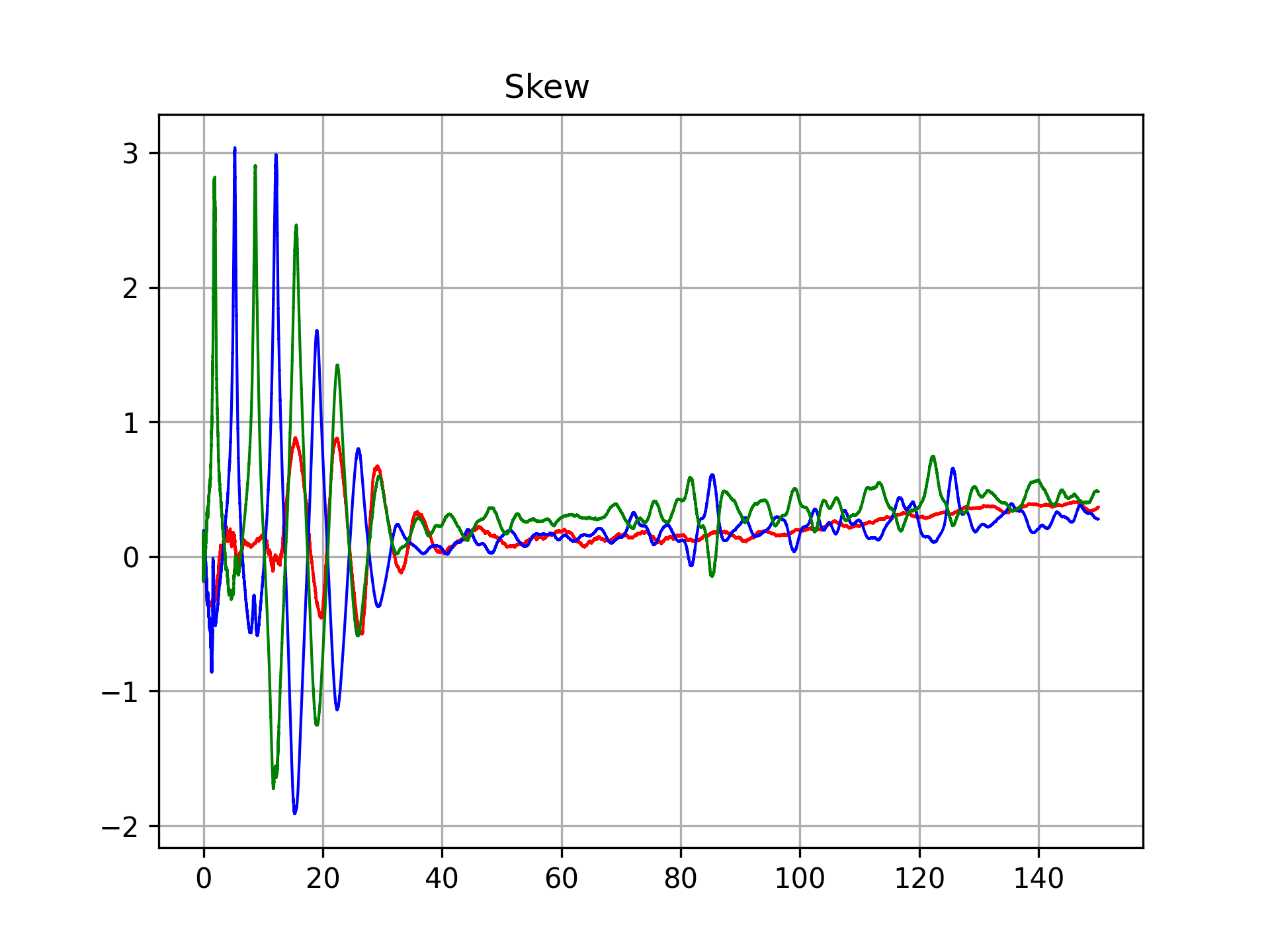}
    \caption{EST skew.}
    \label{fig:LU_skew}
\end{subfigure}
\hfill
\begin{subfigure}{0.23\textwidth}
    \includegraphics[width=\textwidth]{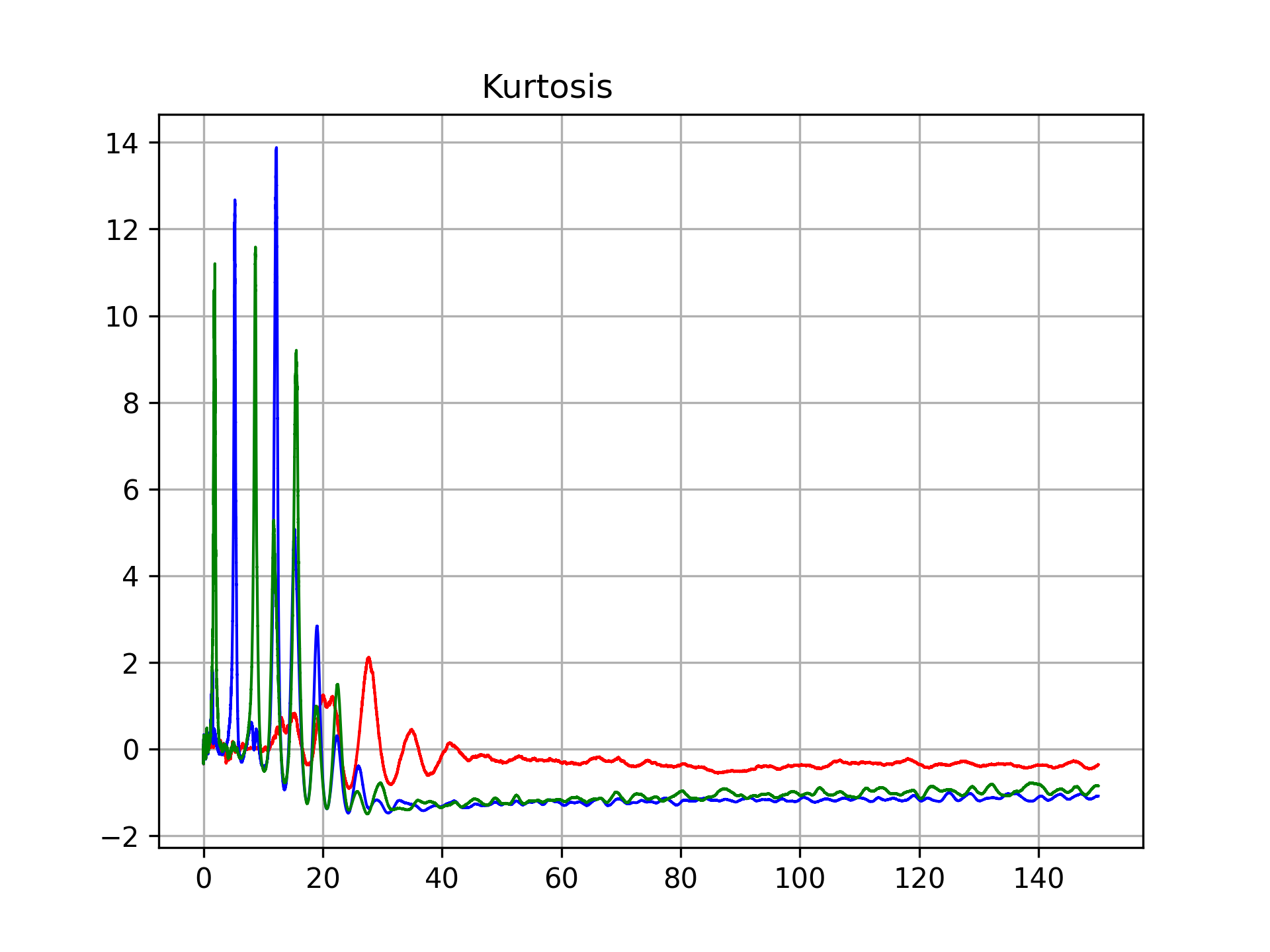}
    \caption{EST kurtosis.}
    \label{fig:LU_kurtosis}
\end{subfigure}
\begin{subfigure}{0.23\textwidth}
    \includegraphics[width=\textwidth]{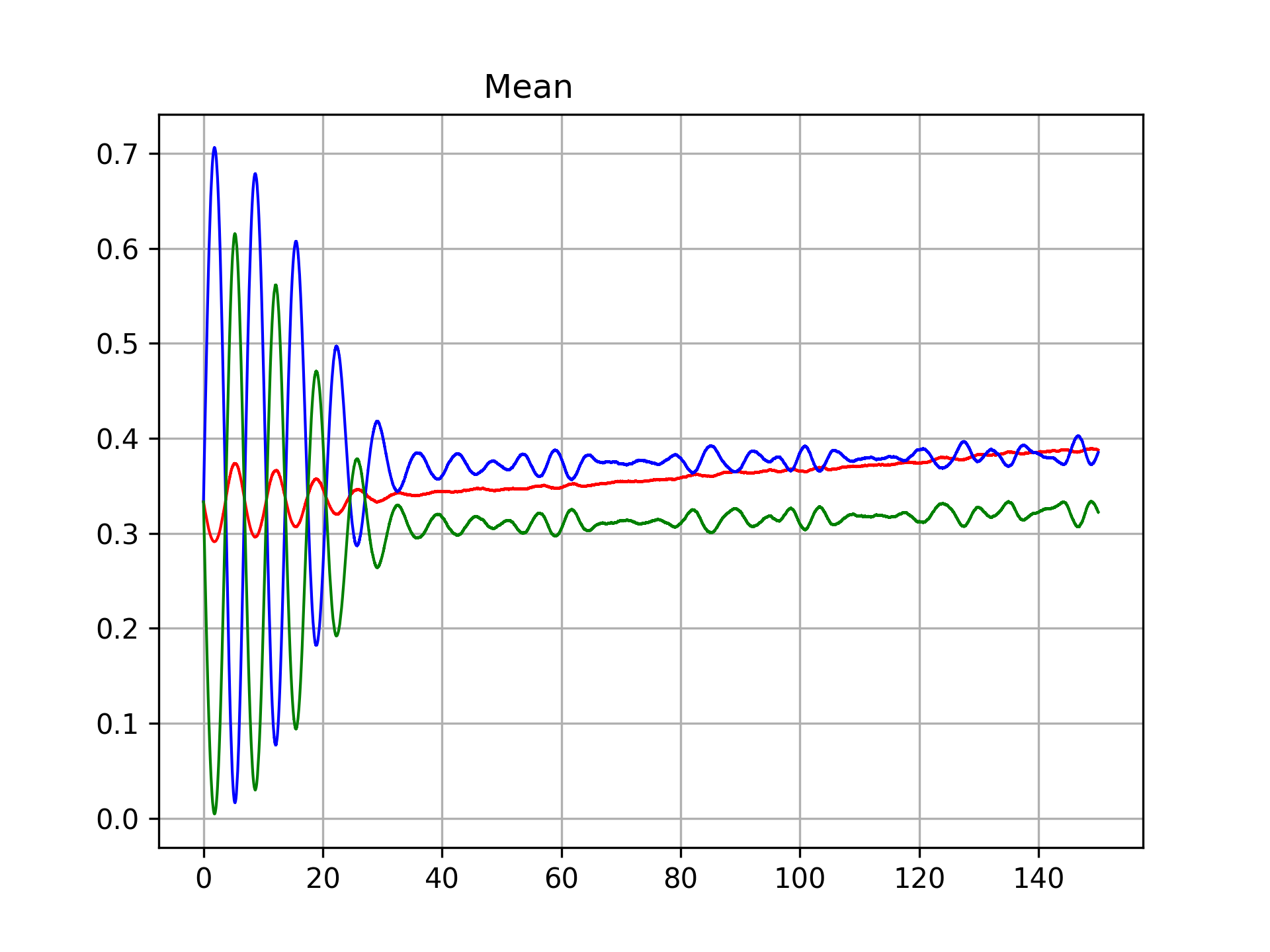}
    \caption{HST mean.}
    \label{fig:SALT_mean}
\end{subfigure}
\hfill
\begin{subfigure}{0.23\textwidth}
    \includegraphics[width=\textwidth]{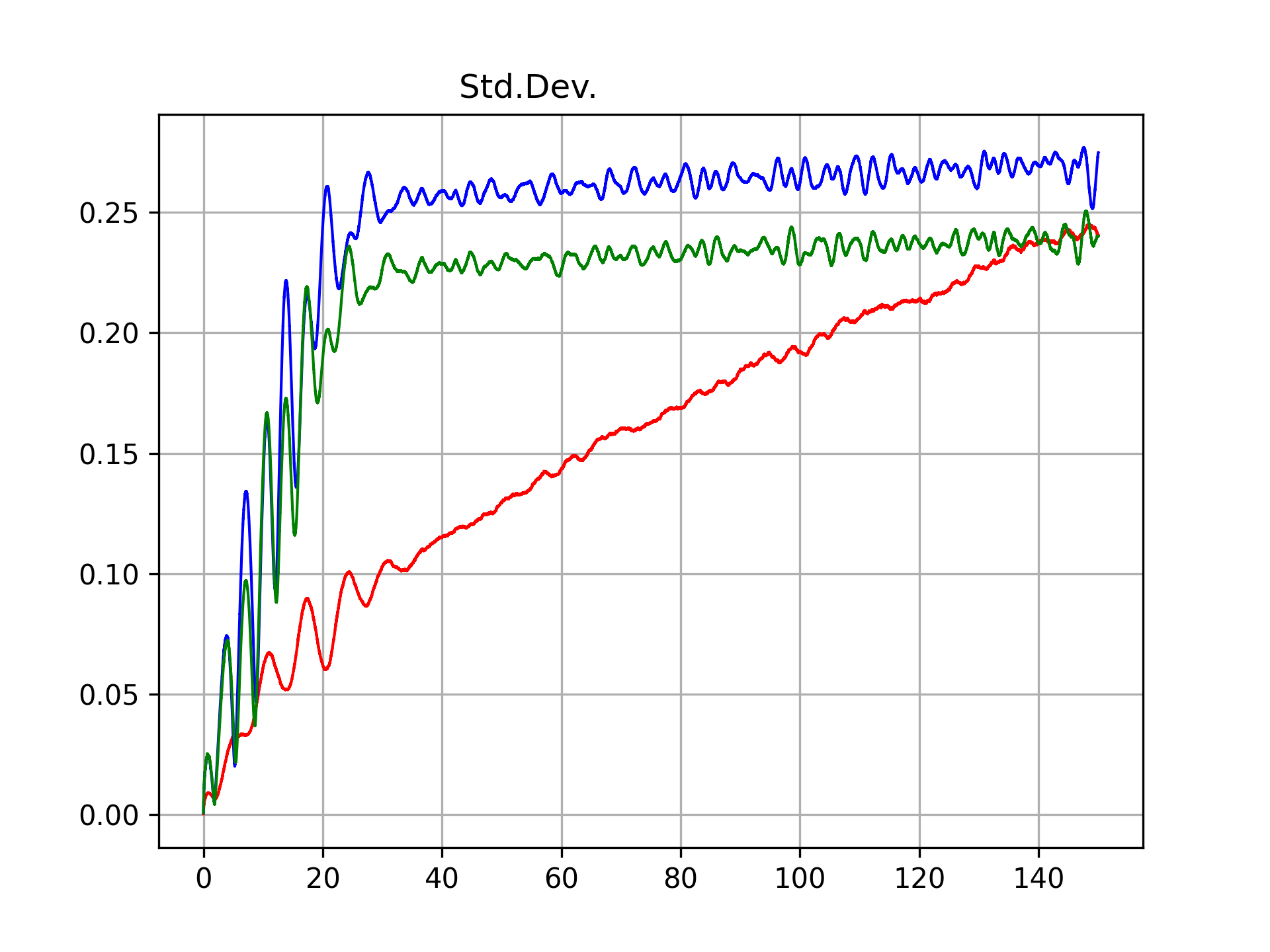}
    \caption{HST std.~dev.}
    \label{fig:SALT_std}
\end{subfigure}
\hfill
\begin{subfigure}{0.23\textwidth}
    \includegraphics[width=\textwidth]{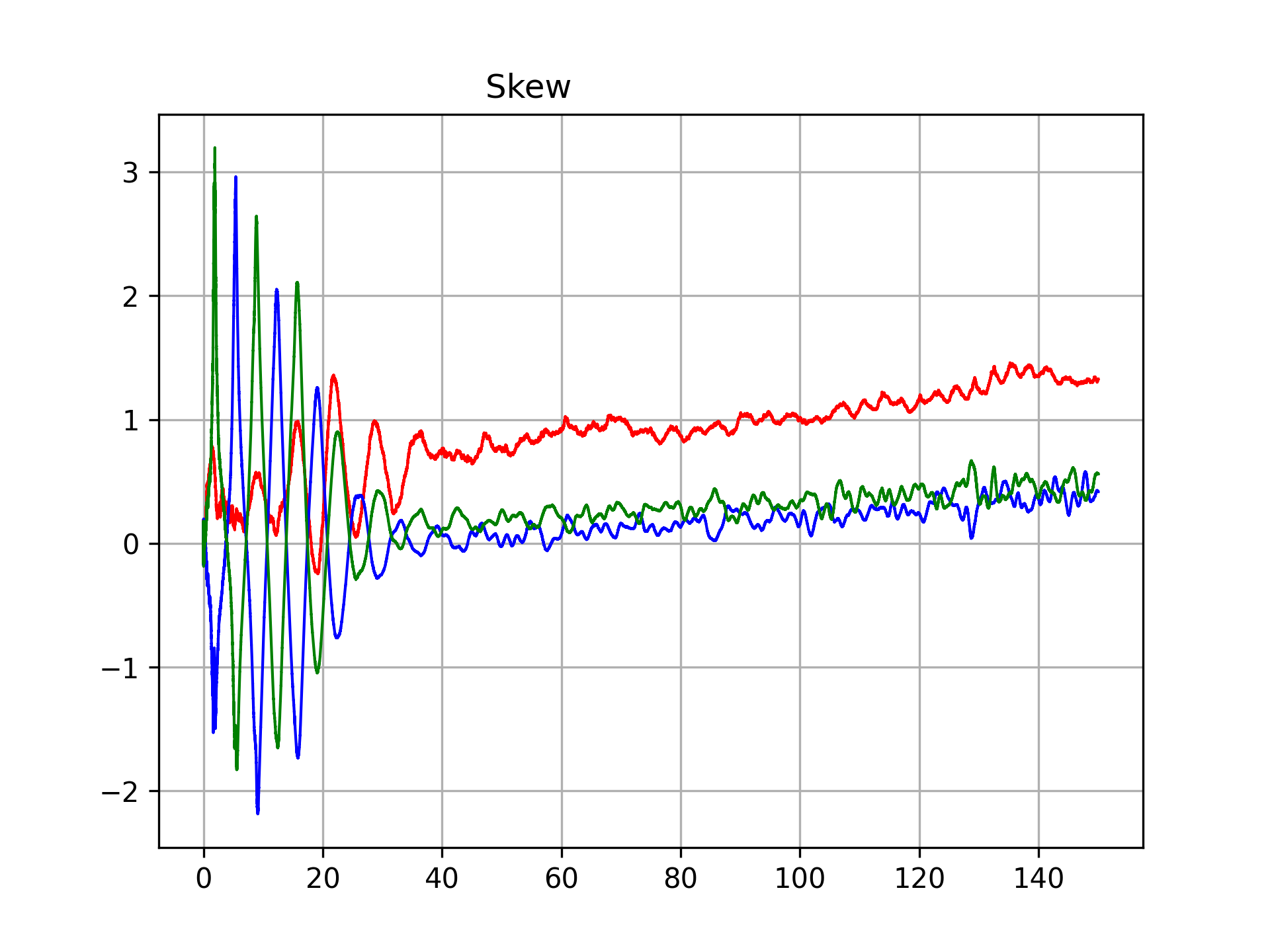}
    \caption{HST skew.}
    \label{fig:SALT_skew}
\end{subfigure}
\hfill
\begin{subfigure}{0.23\textwidth}
    \includegraphics[width=\textwidth]{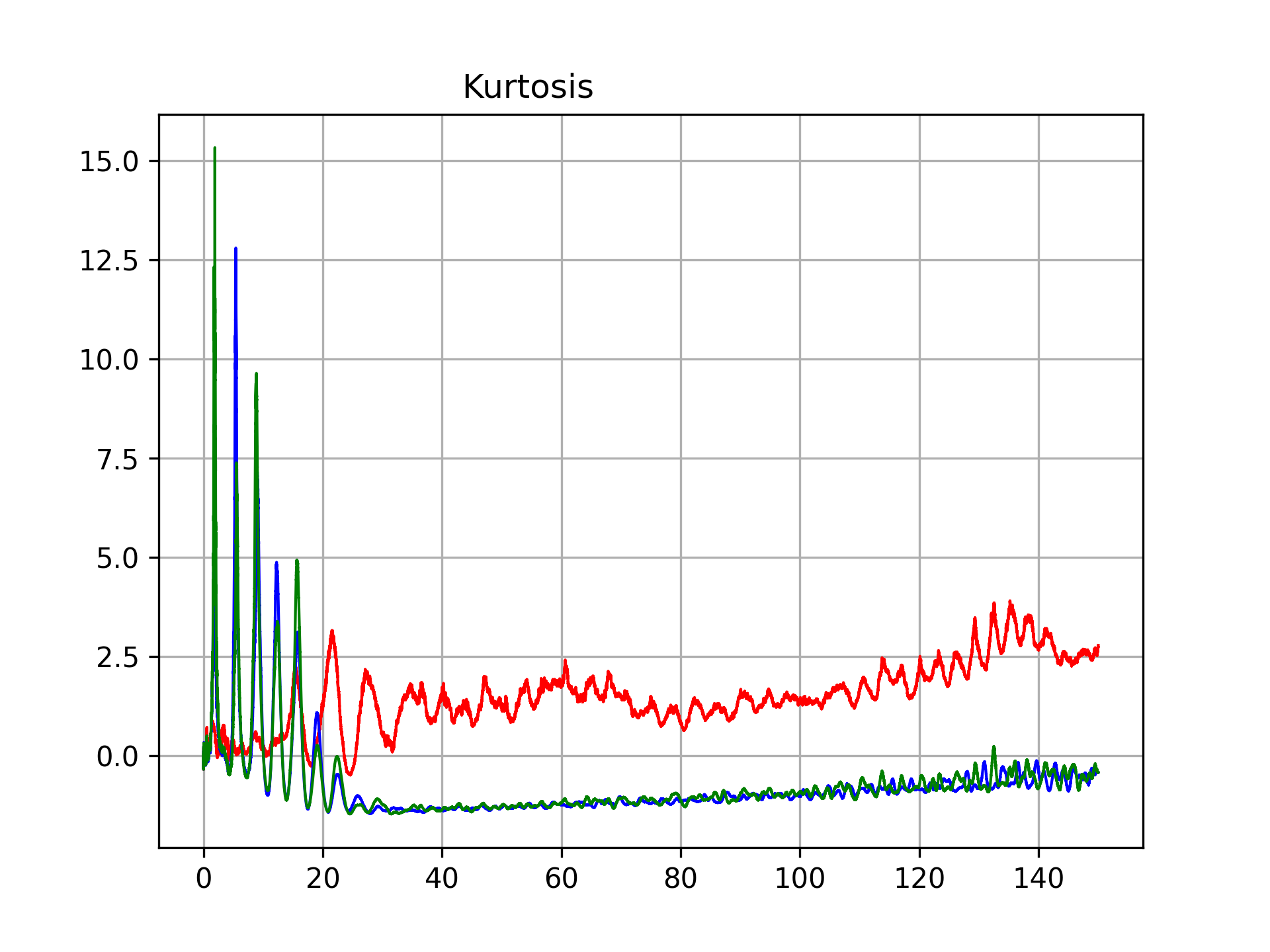}
    \caption{HST kurtosis.}
    \label{fig:SALT_kurtosis}
\end{subfigure}

\caption{Evolution of statistical moments of the modal energies for both stochastic models in the full noise case. The statistics are computed pointwise in time from an ensemble of $1000$ realisations up to a final time of $150$. (a+e) Ensemble mean. (b+f) Ensemble standard deviation. (c+g) Ensemble skew. (d+h) Ensemble kurtosis.}
\label{fig:model_stats}

\end{figure}

Figure~\ref{fig:model_stats} shows various statistics of the generated ensembles of $n=1000$ particles for the HST and EST triad models in the full noise case introduced above. We plot the ensemble mean, standard deviation, skew, and kurtosis.

\begin{figure}

\centering
\begin{subfigure}{0.32\textwidth}
    \includegraphics[width=\textwidth]{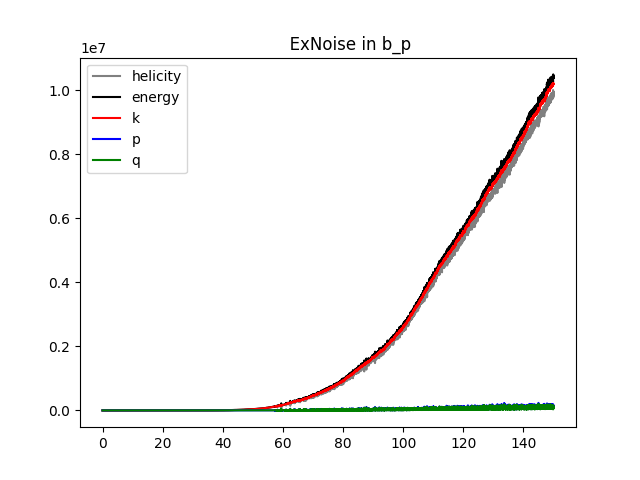}
    \caption{HST.}
    \label{fig:SALT_ex_noise_p}
\end{subfigure}
\hfill
\begin{subfigure}{0.32\textwidth}
    \includegraphics[width=\textwidth]{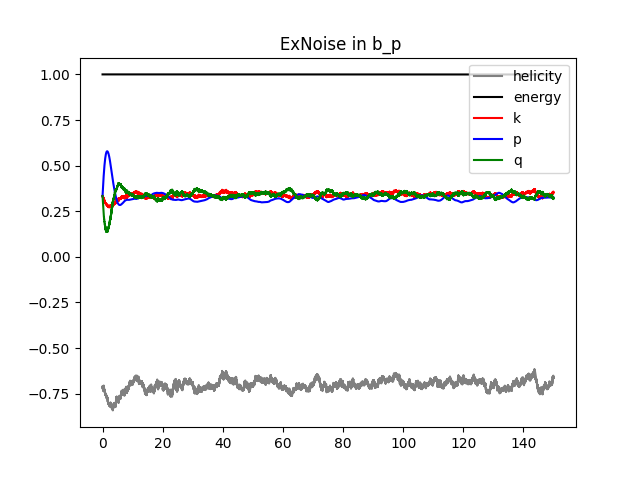}
    \caption{EST.}
    \label{fig:LU_ex_noise_p}
\end{subfigure}
\hfill
\begin{subfigure}{0.32\textwidth}
    \includegraphics[width=\textwidth]{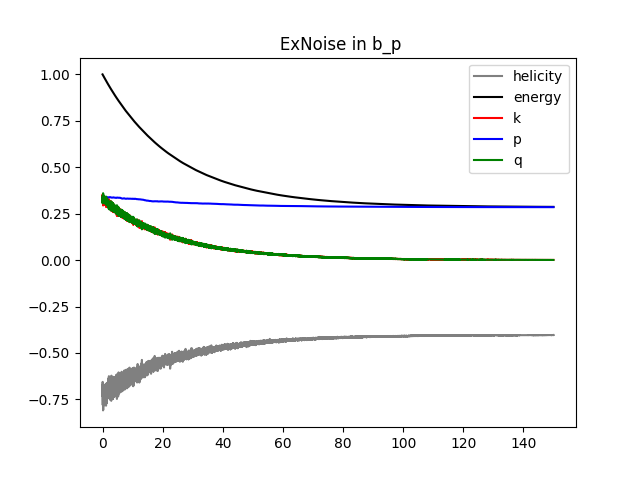}
    \caption{EST strong noise.}
    \label{fig:LU_Ex_noise_extreme_p}
\end{subfigure}
\caption{The effect of large noise coefficients on the mean. Evolution of the mean modal energies (colored lines), mean total energy (black), and mean helicity (grey) for the HST (a) and EST (b) model with noise coefficient $\bb = [0,1,0]$. (c) Evolution of the mean modal energies, mean total energy, and mean helicity for the EST model with the strong noise coefficient $\bb = [0,10,0]$.}
\label{fig:ex_noise}

\end{figure}

The effect of large noise coefficients is exemplified in Figure~\ref{fig:ex_noise}.
We observe that the HST model explodes whereas  the EST model is more tolerant to large noise coefficients, and even in the extreme case, does not become unstable in the mean. The ensemble means are computed from $n=500$ realisations, using the noise coefficient $\bb = [0.0,1.0,0.0]$. Moreover, to stress the EST model, we also ran the same experiment with a noise coefficient of $\bb = [0.0,10.0,0.0]$ for the EST model alone.

\begin{figure}

\centering
\begin{subfigure}{0.48\textwidth}
    \includegraphics[width=\textwidth]{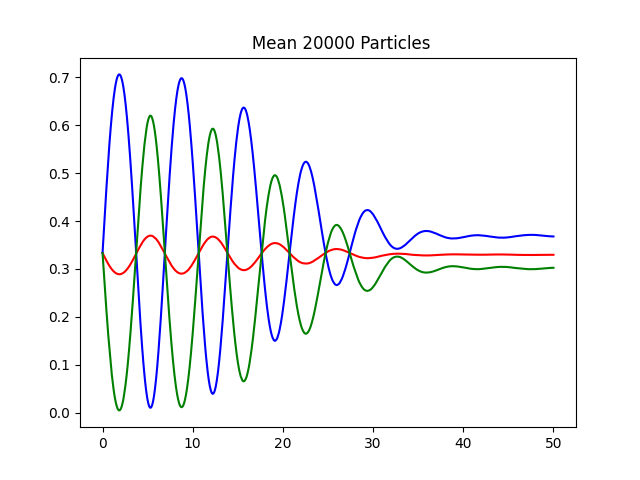}
    \caption{EST.}
    \label{fig:LU_many_particles}
\end{subfigure}
\hfill
\begin{subfigure}{0.48\textwidth}
    \includegraphics[width=\textwidth]{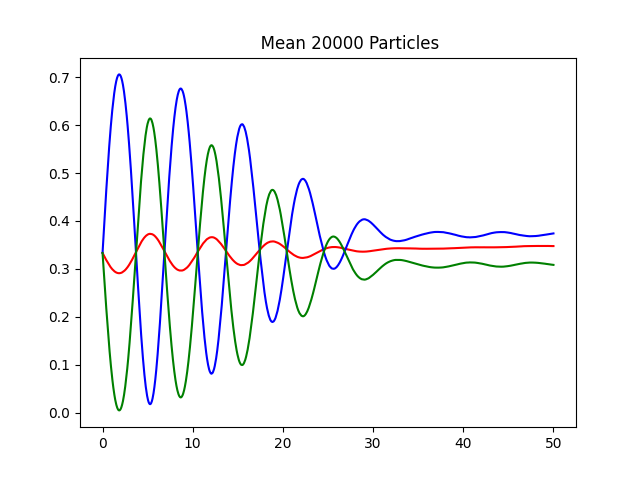}
    \caption{HST.}
    \label{fig:SALT_many_particles}
\end{subfigure}
\caption{Evolution of mean modal energies for a very large number of realisations for the EST (a) and HST (b) models. The mean is computed from $20,\!000$ particles in the full noise case.}
\label{fig:many_particles}
\end{figure}

The mean ensemble for a large number of particles, $n=20,\!000$, is shown in Figure~\ref{fig:many_particles}. We can observe that, compared to Figure~\ref{fig:LU_mean} and \ref{fig:SALT_mean} the oscillations after time $40$ are reduced for a very large number of particles. Hence, we believe that the system stabilizes in the mean to stationary modal energies as the limiting effect of the noise.

\subsubsection{Data Assimilation}
\label{sec:Numerical_DA}
Using the two stochastic models in the full noise case described above, we perform the data assimilation tests using the following framework:

The signal process (the \emph{truth}) is given by the deterministic triad model. The observations are the modal energies of the deterministic model, observed every $10$ time units, and perturbed by noise of the form
\begin{equation}
    \eta\sim\mathcal{N}(\mathbf{0}, \mathrm{C}),
\end{equation}
where the covariance matrix $\mathrm{C}\in\mathbb{R}^3$ is chosen to be the diagonal matrix $\mathrm{C}=\operatorname{diag}(0.005^2, 0.05^2, 0.05^2)$.

\begin{figure}

\centering
\begin{subfigure}{0.3\textwidth}
    \includegraphics[width=\textwidth]{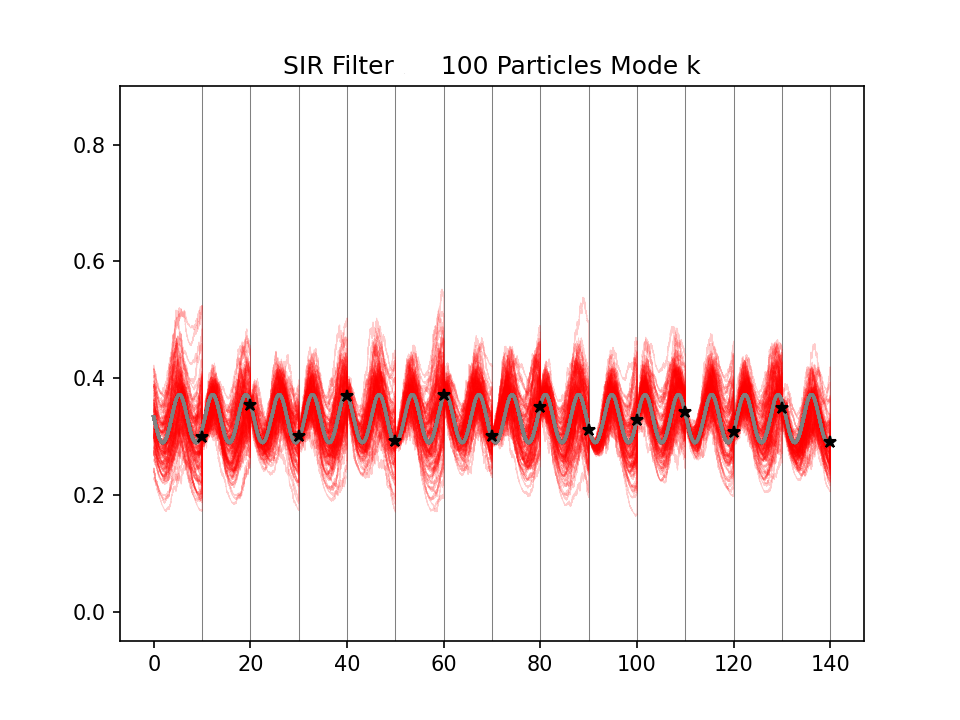}
    \caption{EST filter k.}
    \label{fig:LU_filter_k}
\end{subfigure}
\hfill
\begin{subfigure}{0.3\textwidth}
    \includegraphics[width=\textwidth]{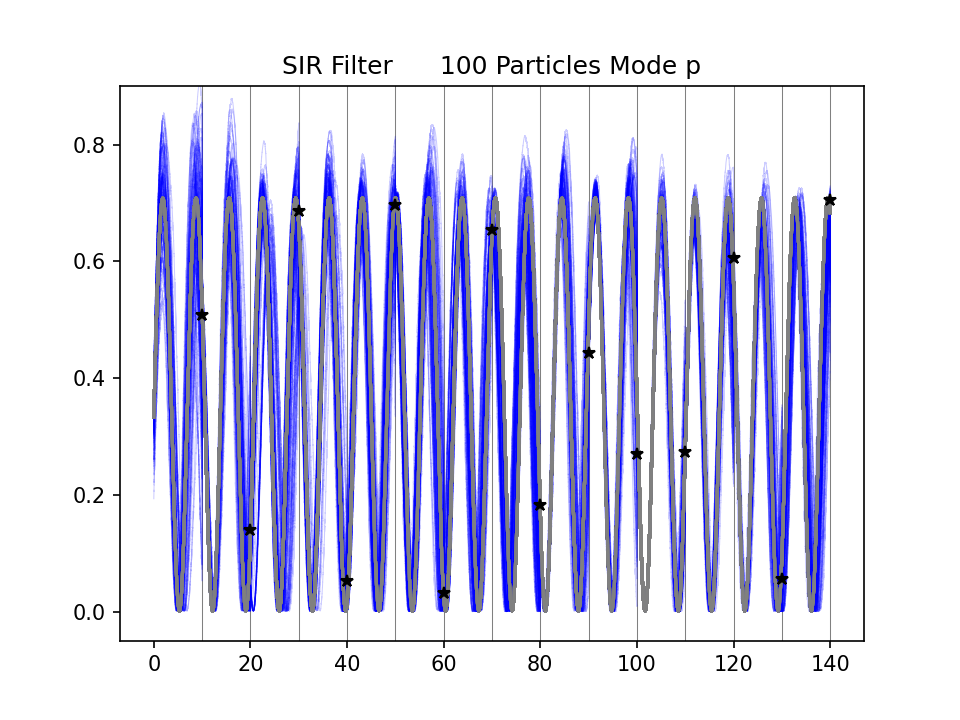}
    \caption{EST filter p.}
    \label{fig:LU_filter_p}
\end{subfigure}
\hfill
\begin{subfigure}{0.3\textwidth}
    \includegraphics[width=\textwidth]{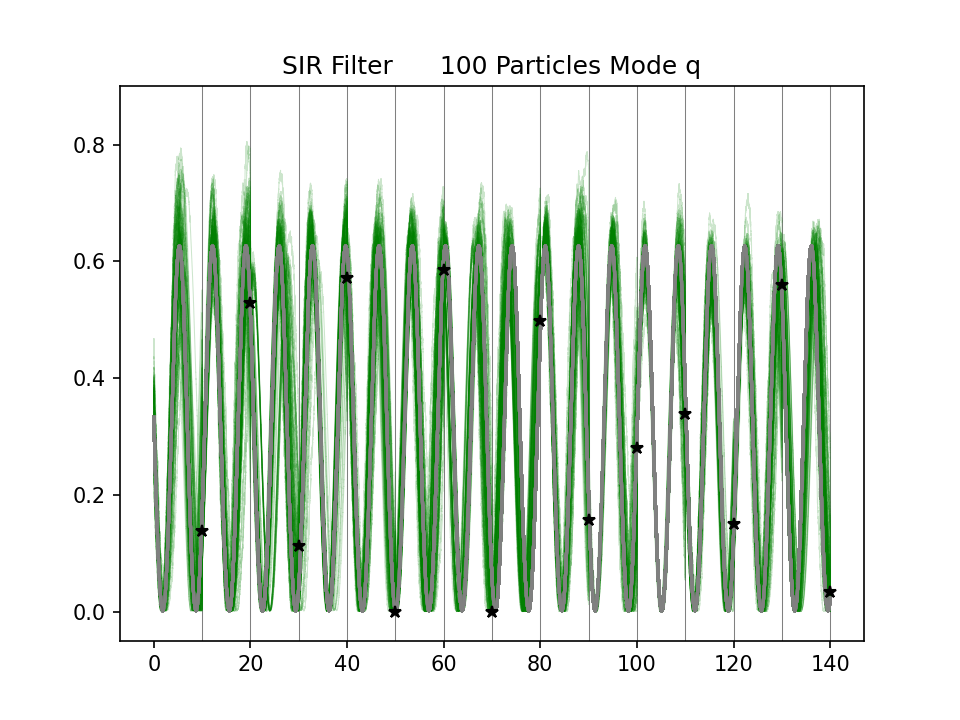}
    \caption{EST filter q.}
    \label{fig:LU_filter_q}
\end{subfigure}

\begin{subfigure}{0.3\textwidth}
    \includegraphics[width=\textwidth]{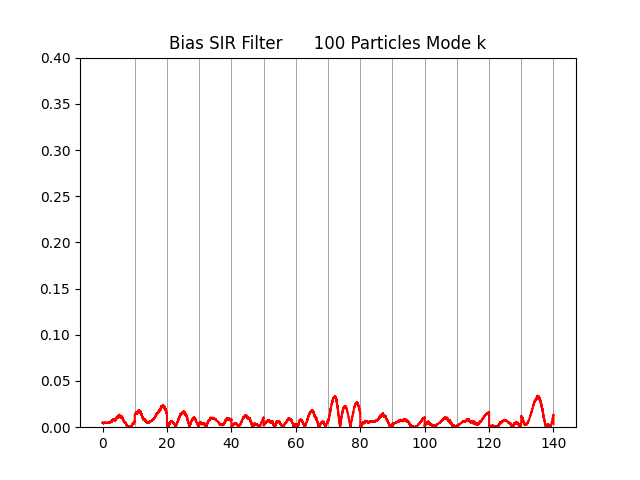}
    \caption{EST bias k.}
    \label{fig:LU_bias_k}
\end{subfigure}
\hfill
\begin{subfigure}{0.3\textwidth}
    \includegraphics[width=\textwidth]{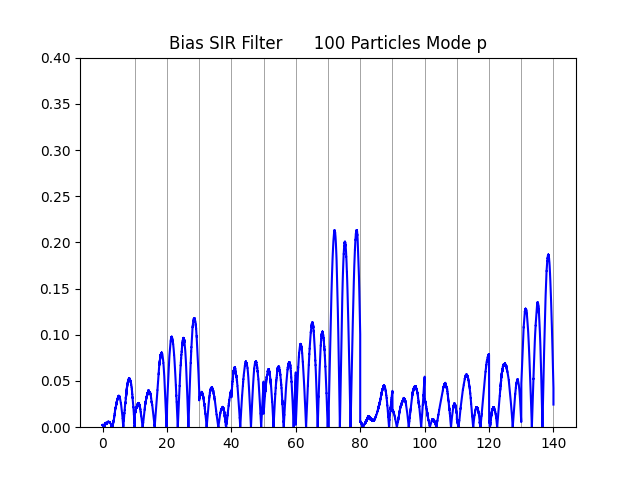}
    \caption{EST bias p.}
    \label{fig:LU_bias_p}
\end{subfigure}
\hfill
\begin{subfigure}{0.3\textwidth}
    \includegraphics[width=\textwidth]{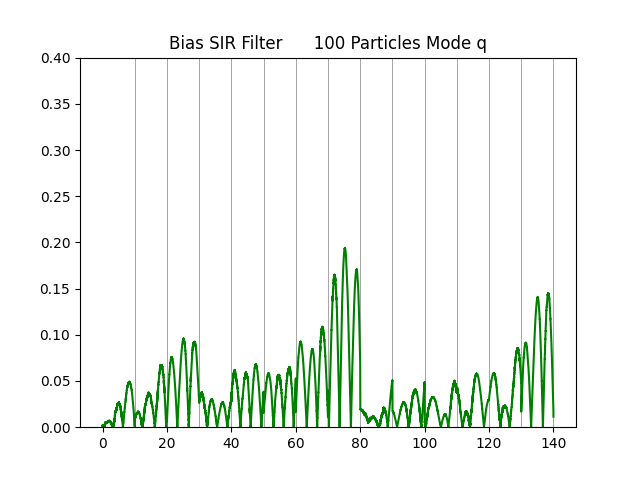}
    \caption{EST bias q.}
    \label{fig:LU_bias_q}
\end{subfigure}
\begin{subfigure}{0.3\textwidth}
    \includegraphics[width=\textwidth]{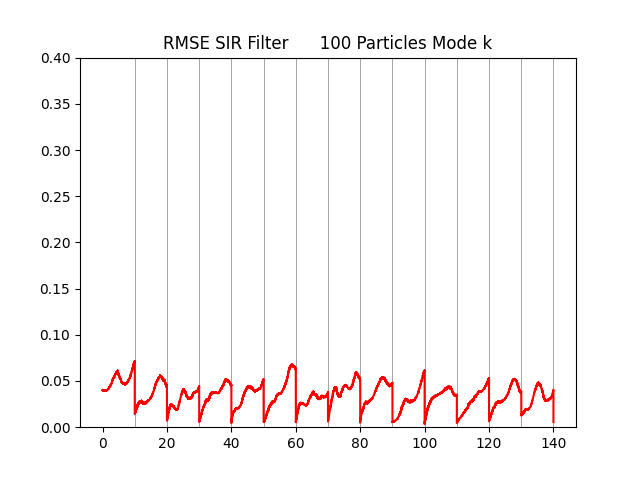}
    \caption{EST RMSE k.}
    \label{fig:LU_rmse_k}
\end{subfigure}
\hfill
\begin{subfigure}{0.3\textwidth}
    \includegraphics[width=\textwidth]{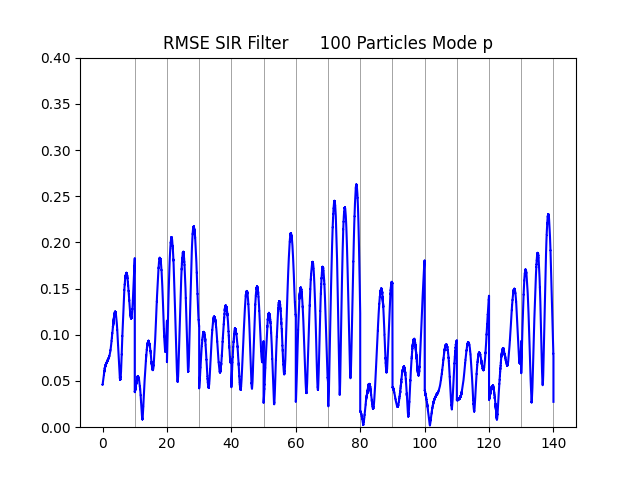}
    \caption{EST RMSE p.}
    \label{fig:LU_rmse_p}
\end{subfigure}
\hfill
\begin{subfigure}{0.3\textwidth}
    \includegraphics[width=\textwidth]{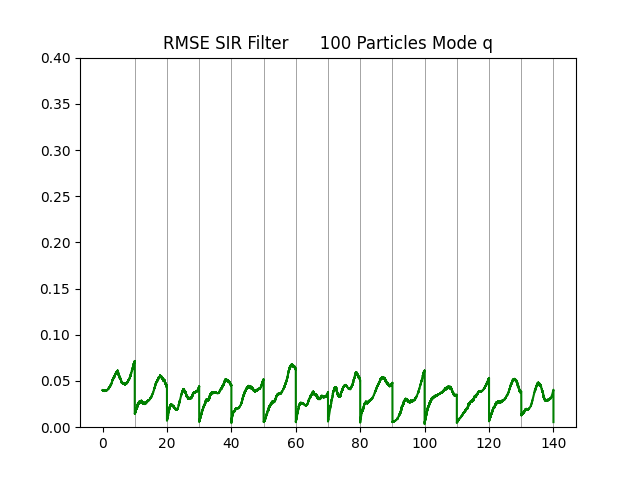}
    \caption{EST RMSE q.}
    \label{fig:LU_rmse_q}
\end{subfigure}

\caption{Filtering experiment for EST model using SIR particle filter. (a-c) Ensemble evolution for $100$ particles in mode k (a, red), p (b, blue), and q (c, green). The signal (grey) is the deterministic model and the observations (black stars) are noisy and taken and assimilated every $10$ time units. (d-f) Bias of the filtering ensemble. (g-h) RMSE of the filtering ensemble.}
\label{fig:filtering_LU}
\end{figure}

\begin{figure}
\centering
\begin{subfigure}{0.3\textwidth}
    \includegraphics[width=\textwidth]{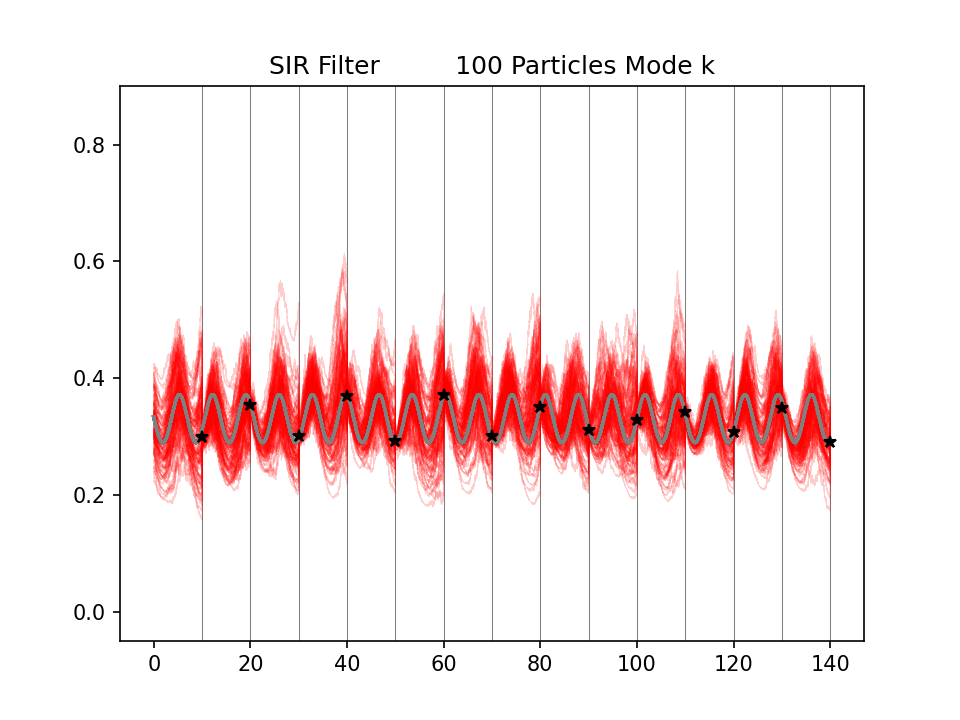}
    \caption{HST filter k.}
    \label{fig:SALT_filter_k}
\end{subfigure}
\hfill
\begin{subfigure}{0.3\textwidth}
    \includegraphics[width=\textwidth]{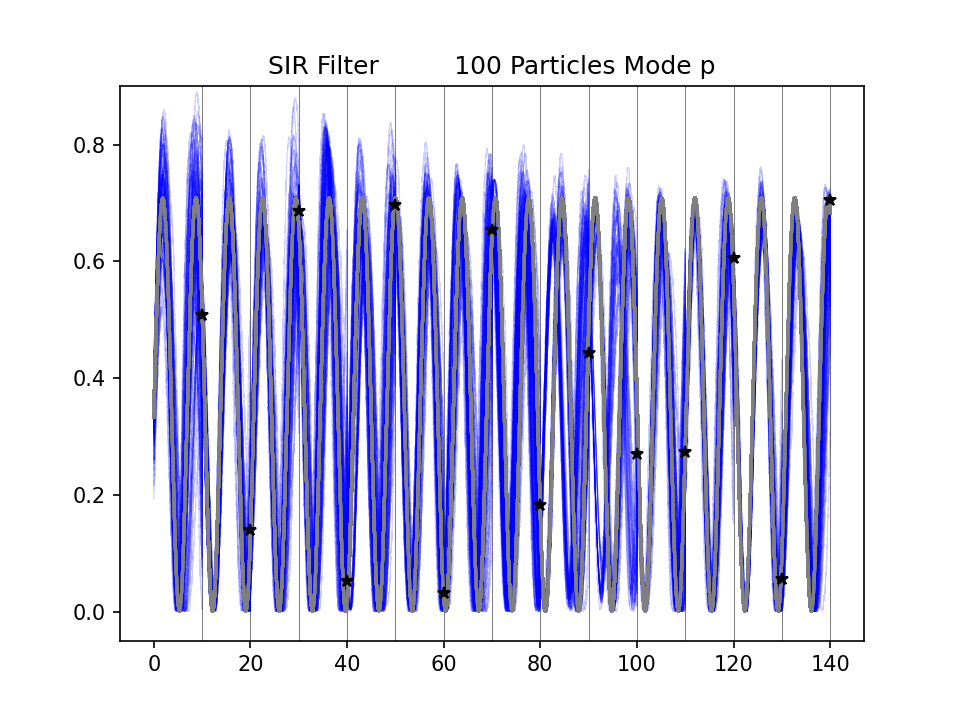}
    \caption{HST filter p.}
    \label{fig:SALT_filter_p}
\end{subfigure}
\hfill
\begin{subfigure}{0.3\textwidth}
    \includegraphics[width=\textwidth]{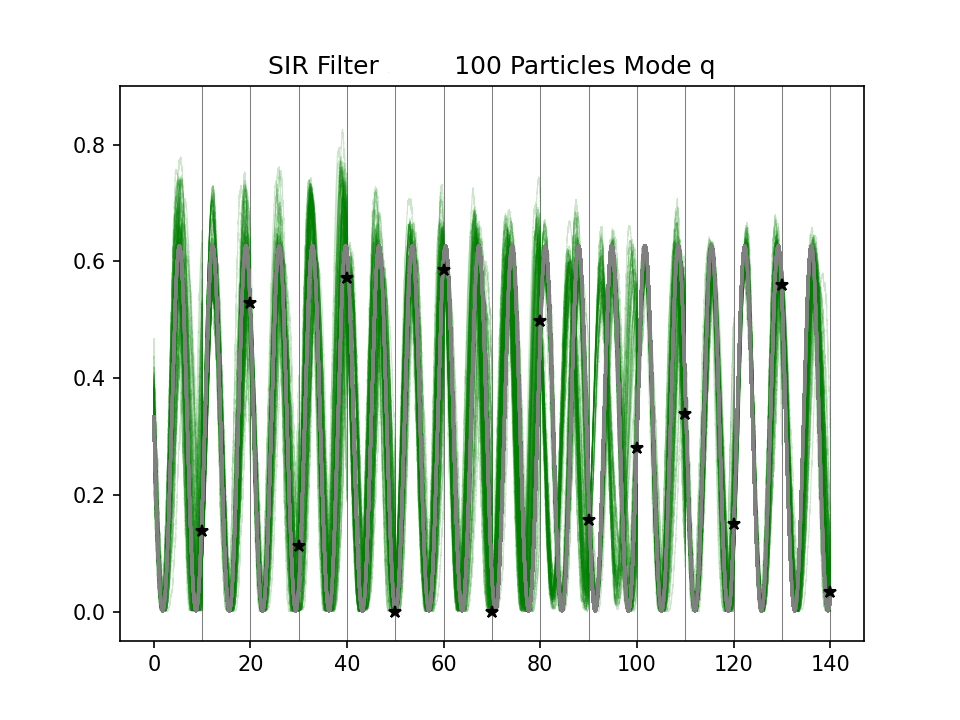}
    \caption{HST filter q.}
    \label{fig:SALT_filter_q}
\end{subfigure}

\begin{subfigure}{0.3\textwidth}
    \includegraphics[width=\textwidth]{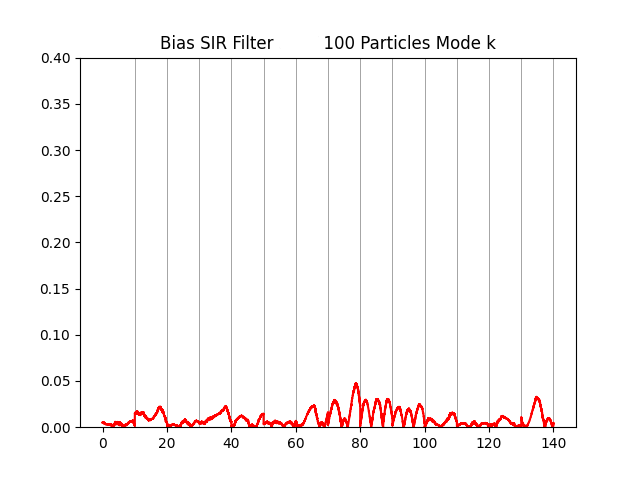}
    \caption{HST bias k.}
    \label{fig:SALT_bias_k}
\end{subfigure}
\hfill
\begin{subfigure}{0.3\textwidth}
    \includegraphics[width=\textwidth]{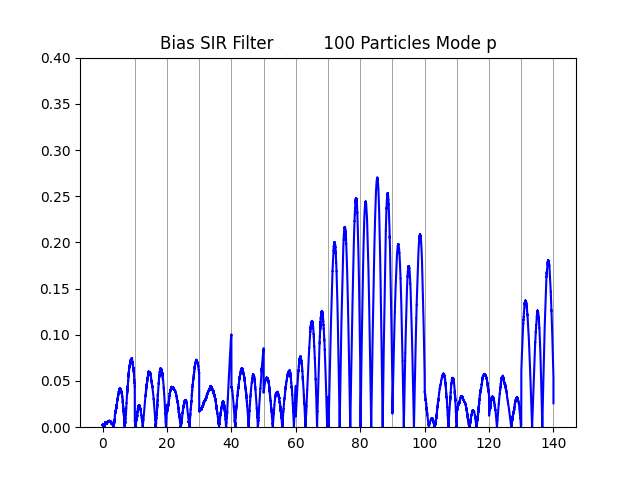}
    \caption{HST bias p.}
    \label{fig:SALT_bias_p}
\end{subfigure}
\hfill
\begin{subfigure}{0.3\textwidth}
    \includegraphics[width=\textwidth]{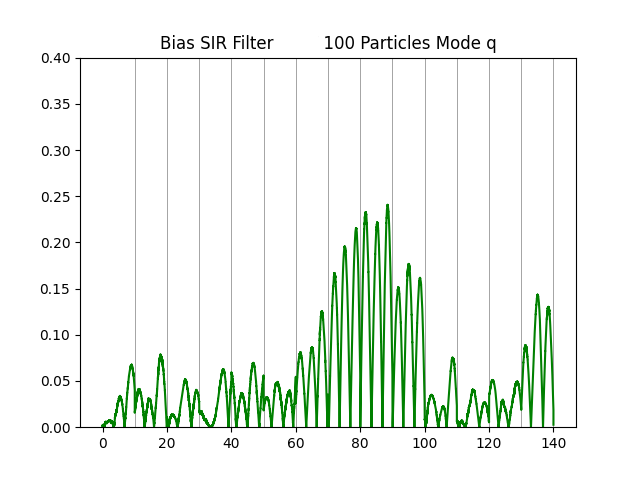}
    \caption{HST bias q.}
    \label{fig:SALT_bias_q}
\end{subfigure}
\begin{subfigure}{0.3\textwidth}
    \includegraphics[width=\textwidth]{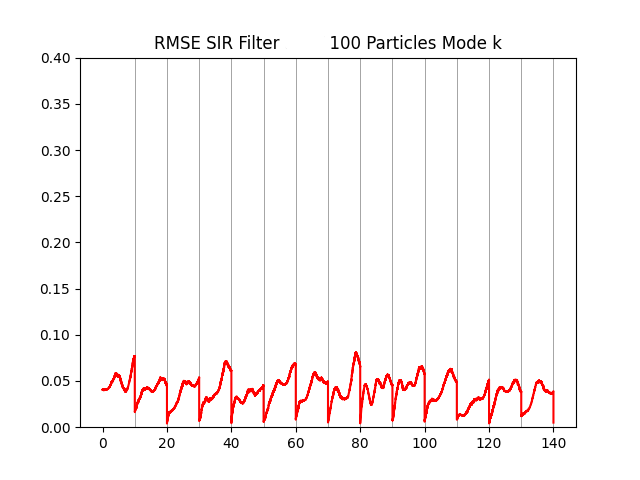}
    \caption{HST RMSE k.}
    \label{fig:SALT_rmse_k}
\end{subfigure}
\hfill
\begin{subfigure}{0.3\textwidth}
    \includegraphics[width=\textwidth]{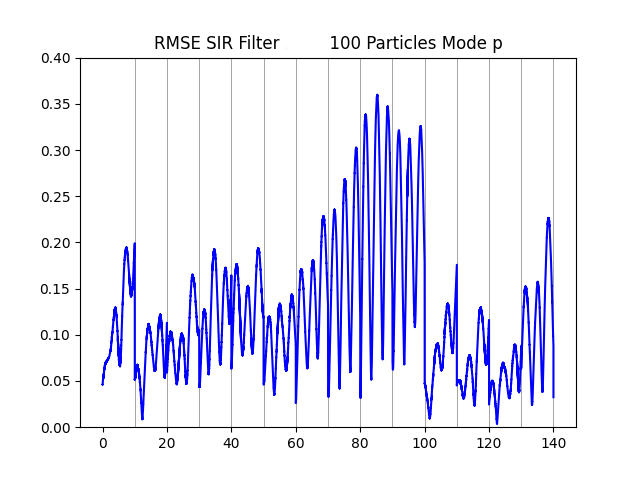}
    \caption{HST RMSE p.}
    \label{fig:SALT_rmse_p}
\end{subfigure}
\hfill
\begin{subfigure}{0.3\textwidth}
    \includegraphics[width=\textwidth]{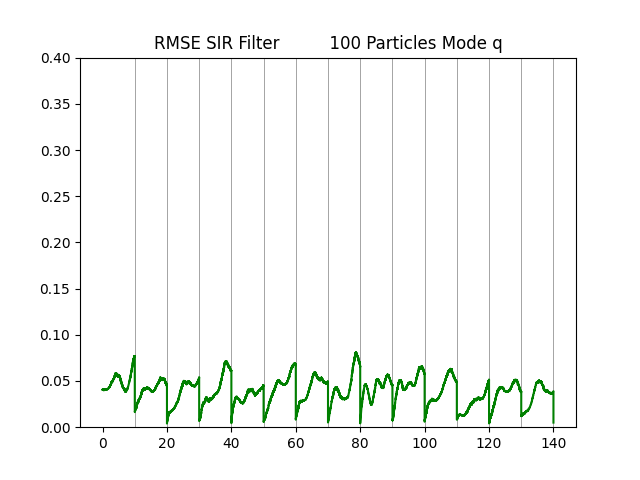}
    \caption{HST RMSE q.}
    \label{fig:SALT_rmse_q}
\end{subfigure}

\caption{Filtering experiment for HST model using SIR particle filter. (a-c) Ensemble evolution for $100$ particles in mode k (a, red), p (b, blue), and q (c, green). The signal (grey) is the deterministic model and the observations (black stars) are noisy and taken and assimilated every $10$ time units. (d-f) Bias of the filtering ensemble. (g-h) RMSE of the filtering ensemble.}
\label{fig:filtering_SALT}

\end{figure}

We use the sequential importance resampling (SIR) particle filter to assimilate the periodically observed signal process under the influence of observation noise. The particles evolve according to the stochastic triad models.
Figures~\ref{fig:filtering_LU}
 and~\ref{fig:filtering_SALT} show the results of filtering the ensemble of $n=100$ particles of the EST and HST triad models, respectively. The ensembles are assessed in terms of the bias and RMSE statistics. We analyse the comparison details below:

\paragraph{Mode k} This is the least energetic of all the modes (hence the reason why we observe it with the least amount of measurement noise). The cloud of particles is well placed around the truth even with the small sample. The bias remains small for both the HST and the EST versions and it reduces significantly when observations are assimilated more frequently (see Figure~\ref{fig:filtering_LU_small} in Appendix~\ref{sec:filter-verification}) as well as when we use a large number ($500$) of particles (see Figure~\ref{fig:filtering_LU_500} in Appendix~\ref{sec:filter-verification}). The RMSE remains small in all cases and decreases (though not substantially) when the DA step is small.           

\paragraph{Modes p and q} These are the two energetic modes of the system. We used here a measurement noise that is one order of magnitude larger. Despite this, the results remain equally good. The cloud of particles provide a good envelope for the truth at all times. This validates the choice of the stochasticity: the uncertainty is properly modelled. For both models the bias can become very large, reaching $30\%$ of the size of the oscillations for the HST model and $25\%$ of the size of the oscillations for the EST model. As expected, it is drastically reduced when observations are assimilated more frequently. The RMSE for mode p is also large but substantially smaller for mode q. The addition of more intermediate DA steps or more particles has a less pronounced effect for the q mode. 

\begin{remark}
\label{rem:ess}
We record the Effective Sample Size (ESS) for a typical run (for both EST and HST) in Figure \ref{fig:filtering_ess}. As usual, the ESS is computed just before the application of the resampling procedure. The ESS is seen to decay dramatically from 100 down to single digits numbers in most instances in time.     
\end{remark}

\begin{remark}
\label{rem:stats_EST}
We record the results over 10 independent runs of the filtering experiment for the EST model with 500 ensemble members in figure \ref{fig:EST_mean_stats}. More precisely, in graphs \ref{fig:EST_mean_k}, \ref{fig:EST_mean_p} and \ref{fig:EST_mean_q} each, we plot the mean across the 10 independent filtering runs together with the the evolution of the signal and the individual ensemble means for each mode.
The mean bias as well as the envelope obtained from the independent runs are shown in graphs \ref{fig:EST_bias_k_mean}, \ref{fig:EST_bias_p_mean} and \ref{fig:EST_bias_q_mean}.
The same is shown for the RMSE in graphs \ref{fig:EST_rmse_k_mean}, \ref{fig:EST_rmse_p_mean} and \ref{fig:EST_rmse_q_mean}.
Compared with a single run of the same experiment reported in Figure \ref{fig:filtering_LU_500} we observe the approximations are now near perfect (the statistical error has been drastically reduced). 

\end{remark}

\section{Conclusions}\label{sect:conclusions}

The introduction of stochasticity into the deterministic triad models leads to two new stochastic models. Stochasticity is introduced in a principled way (rather than ad-hoc). It starts with a full scale fluid dynamic model which is randomly perturbed. At the full scale, the stochastic parametrisation models the small-scale effects in fluid dynamics modelling. In particular, it efficiently captures the high-frequency small-scale dynamics and correctly correlates it with the slow, large-scale fluid motion. In addition, it is constrained to conserve either the helicity or the kinetic energy of the system. This inspires two different  stochastic triad models of Euler type which we  compare using data assimilation procedures based on particle filtering. The methodology we employ can be used as a benchmark when analysing new types of stochastic parametrisations: ours is the first study that assesses the efficiency of stochastic parametrisations from a data assimilation perspective. 

The introduction of stochasticity ensures that the
correct spread (one that preserves the physical properties of the system)  is introduced in the ensemble of particles. In its absence, the particle filter degenerates quite rapidly: after a few DA steps, a single particle survives the culling procedure which does offer a good approximation of the truth. A purely deterministic transition kernel does not work, generating a rapid degeneracy of the particle filter. 

The two stochastic systems (one which preserves helicity and the other one which preserves energy) are analysed using a standard particle filter. There is no need for additional procedures (such as tempering, nudging, or jittering). They perform equally well from the viewpoint of DA: both the RMSE and the bias are drastically reduced and stabilised when the noise is carefully calibrated. The two different stochastic kernels require different noise calibrations in order to perform well in similar data assimilation scenarios. This is somehow expected, given that the underlying stochastic parametrisations preserve different physical quantities.

\nocite{*}
\bibliographystyle{plainurl}
\bibliography{main}

\begin{appendices}
\numberwithin{equation}{section}

\section{Notation and Basic Identities}\label{appendix:notations}

\subsection{Notation}

In this work we use the following notation. We write $\Z$, $\R$ and $\C$ for the sets of integers, real numbers and complex numbers, respectively. Boldface letters denote three-dimensional complex vectors. For two complex vectors $\ba$ and $\bb$ with components $a_j$ and $b_j$, their dot product is denoted by
\begin{equation*}
    \ba\cdot\bb = \sum_{j=1}^3 a_j b_j = a_j b_j \in \C.
\end{equation*}
This paper uses the Einstein convention of summing over repeated indices. The norm of the complex vector $\ba$ is defined as
\begin{equation*}
    |\ba| = \sqrt{\ba \cdot \ba^*} = \sqrt{a_j a_j^*} \geq 0,
\end{equation*}
with the superscript symbol $^*$ denoting complex conjugation.
Further, the cross product of two vectors $\ba$ and $\bb$ is given by
\begin{equation*}
    \ba\times\bb = (
    a_2 b_3-a_3b_2,
    a_3b_1-a_1b_3, 
    a_1b_2-a_2b_1
    ) \in \C^3.
\end{equation*}
The gradient of a scalar field $\phi: \mathcal{D}\subseteq \R^3 \to \C$ at a point $\bx\in\mathcal{D}$ is denoted by
\begin{equation*}
    \nabla\phi(\bx) = 
    \big(
    \partial_1\phi(\bx), 
    \partial_2\phi(\bx),
    \partial_3\phi(\bx)
    \big) \in \C^3.
\end{equation*}
The divergence of a vector field $\bs{\psi}: \mathcal{D}\subseteq \R^3 \to \C^3$ with components $\psi_j$ at a point $\bx\in\mathcal{D}$ is defined as
\begin{equation*}
    \nabla\cdot\bs{\psi}(\bx) = \partial_1{\psi}_1(\bx)+\partial_2{\psi}_2(\bx)+\partial_3{\psi}_3(\bx) = \partial_j{\psi}_j(\bx) \in \C
\end{equation*}
and the curl of $\bs{\psi}$ at $\bx$ is given by
\begin{equation*}
    \nabla\times\bs{\psi}(\bx) = 
    \big(
    \partial_2\psi_3(\bx)-\partial_3\psi_2(\bx),
    \partial_3\psi_1(\bx)-\partial_1\psi_3(\bx), \partial_1\psi_2(\bx)-\partial_2\psi_1(\bx)
    \big) \in \C^3.
\end{equation*}

\subsection{Vector identities}

For three vectors $\ba$, $\bb$ and $\bc$, we have the following algebraic vector identities
\begin{align*}
    &\ba\cdot(\bb\times\bc)=\bb\cdot(\bc\times\ba)=\bc\cdot(\ba\times\bb),\\
    &\ba\times(\bb\times\bc) = (\ba\cdot\bc)\bb - (\ba\cdot\bb)\bc,\\
    &(\ba\times\bb)\cdot(\bc\times\bd) = (\ba\cdot\bc)(\bb\cdot\bd) - (\bb\cdot\bc)(\ba\cdot\bd),\\
    &\ba\times\ba=\bs{0}.
\end{align*}
Moreover, we have the vector calculus identity
\begin{equation}
    (\ba\cdot \nabla)\bb +  b_j\nabla a_j = -\ba\times(\nabla\times\bb) + \nabla(\ba\cdot\bb).
    \label{eq:vec_calc_curl_form}
\end{equation}

\section{Derivation of Triad Models}\label{appendix:derivTriademodels}

\label{sec:derivations}

\subsection{Deterministic Euler}\label{adx:detEuler}
We compute the projection of the terms corresponding to the deterministic Euler vorticity equation in curl form:
\begin{equation}
\label{eqn: Eulcurl}
    \p_t\bv  - \bv \times (\nabla \times \bv) + \tfrac12\nabla |\bv|^2  = -\nabla  p .
\end{equation}
onto the helical basis.
For the time-derivative we get
\begin{align*}
    \begin{split}
    \int_{\cal D} \partial_t\bv(\bx, t)\cdot \bh_{s_k}^*(\bk)e^{-i\bk\cdot\bx} \id\bx
    &= \sum_{\bp}\sum_{s_p} \partial_t a_{s_p}(\bp, t)\bh_{s_p}(\bp)\cdot \bh_{s_k}^*(\bk)\int_{\cal D}e^{i(\bp-\bk)\cdot\bx}\id\bx\\
    &=L^3 \sum_{s_p} \partial_t a_{s_p}(\bk, t)\bh_{s_p}(\bk)\cdot \bh_{s_k}^*(\bk)\\
    &= L^3 \partial_t a_{s_k}(\bk, t)\bh_{s_k}(\bk)\cdot \bh_{s_k}^*(\bk)\\
    &= 2L^3\partial_t a_{s_k}(\bk, t).
    \end{split}
\end{align*}
The vorticity term gives
\begin{align}
\begin{split}
    &\int_{\cal D} (\bv(\bx, t)\times \bs{\omega}(\bx,t)) \cdot \bh_{s_k}^*(\bk)e^{-i\bk\cdot\bx} \id\bx \\
    &= \sum_{\bp,\bq}\sum_{s_p, s_q} a_{s_p}(\bp, t)s_q|\bq|a_{s_q}(\bq, t)\bh_{s_p}(\bp)\times\bh_{s_q}(\bq)\cdot \bh_{s_k}^*(\bk)\int_{\cal D}e^{i(\bp+\bq-\bk)\cdot\bx}\id\bx\\
    &= L^3\sum_{\bp+\bq+\bk=0}\sum_{s_p, s_q} a_{s_p}^*(\bp, t)s_q|\bq|a_{s_q}^*(\bq, t)\bh_{s_p}^*(\bp)\times\bh_{s_q}^*(\bq)\cdot \bh_{s_k}^*(\bk). 
\end{split}
\label{eq:pre_sym_euler}
\end{align}
Note that we can write \eqref{eq:pre_sym_euler} in a form which is symmetric in $\bp$ and $\bq$ since, renaming the indices,
\begin{align*}
    &\sum_{\bp+\bq+\bk=0}\sum_{s_p, s_q} a_{s_p}^*(\bp, t)s_q|\bq|a_{s_q}^*(\bq, t)\bh_{s_p}^*(\bp)\times\bh_{s_q}^*(\bq)\cdot \bh_{s_k}^*(\bk)\\
    &=\sum_{\bp+\bq+\bk=0}\sum_{s_p, s_q} a_{s_q}^*(\bq, t)s_p|\bp|a_{s_p}^*(\bp, t)\bh_{s_q}^*(\bq)\times\bh_{s_p}^*(\bp)\cdot \bh_{s_k}^*(\bk)\\
    &= -\sum_{\bp+\bq+\bk=0}\sum_{s_p, s_q} a_{s_q}^*(\bq, t)s_p|\bp|a_{s_p}^*(\bp, t) \bh_{s_p}^*(\bp)\times\bh_{s_q}^*(\bq) \cdot \bh_{s_k}^*(\bk).
\end{align*}
Therefore,
\begin{align}
    &\int_{\cal D} (\bv(\bx, t)\times \bs{\omega}(\bx,t)) \cdot \bh_{s_k}^*(\bk)e^{-i\bk\cdot\bx} \id\bx\\
    &=\frac{L^3}{2} \sum_{\bp+\bq+\bk=0}\sum_{s_p, s_q} (s_q|\bq|-s_p|\bp|)a_{s_p}^*(\bp, t)a_{s_q}^*(\bq, t)\bh_{s_p}^*(\bp)\times\bh_{s_q}^*(\bq)\cdot \bh_{s_k}^*(\bk).
\end{align}
Moreover, the gradient terms in \eqref{eqn: Eulcurl} vanish upon expansion into helical modes.
Thus, the Euler equation \eqref{eqn: Eulcurl} in helical basis becomes 
\begin{equation*}
    \partial_t a_{s_k}(\bk, t) = -\frac{1}{4} \sum_{\bp+\bq+\bk=0}\sum_{s_p, s_q} (s_p|\bp|-s_q|\bq|)a_{s_p}^*(\bp, t)a_{s_q}^*(\bq, t)\bh_{s_p}^*(\bp)\times\bh_{s_q}^*(\bq)\cdot \bh_{s_k}^*(\bk).
\end{equation*}

\subsection{SALT Euler}\label{adx:saltEuler}

We expand the 3D SALT Navier-Stokes equation \eqref{eqn: curlform} using \eqref{eqn: defs} and \eqref{eqn: SALTamp}. Assume $b_s(-\bp)=b_s^*(\bp)$.
\begin{equation*}
    \int_{\mathcal{D}} \rmd\bv(\bx,t)\cdot \bh_{s_k}^*(\bk)e^{-i\bk\cdot \bx}\id\bx = \sum_\bp\sum_{s_p} \rmd a_{s_p}(\bp,t)\bh_{s_p}(\bp)\cdot\bh_{s_k}^*(\bk)\int_{\mathcal{D}}e^{i(\bp-\bk)\cdot\bx}\id\bx = 2L^3\rmd a_{s_k}(\bk,t)
\end{equation*}
And
\begin{align*}
    &\int_{\mathcal{D}} \rmd \bx_t(\bx,t) \times \bs{\omega}(\bx, t) \cdot \bh_{s_k}^*(\bk)e^{-i\bk\cdot\bx} \id\bx \\
    &= \sum_{\bp,\bq}\sum_{s_p,s_q} \big[ a_{s_p}(\bp,t)\rmd t + b_{s_p}(\bp)\circ \rmd W_t\big]s_q|\bq|a_{s_q}(\bq,t) \bh_{s_p}(\bp)\times\bh_{s_q}(\bq) \cdot\bh_{s_k}^*(\bk)\int_\mathcal{D} e^{i(\bp+\bq-\bk)\cdot\bx}\id\bx\\
    &=L^3\sum_{\bp+\bq+\bk=0}\sum_{s_p,s_q} \big[ a_{s_p}^*(\bp,t)\rmd t + b_{s_p}^*(\bp)\circ \rmd W_t\big]s_q|\bq|a_{s_q}^*(\bq,t) \bh_{s_p}^*(\bp)\times\bh_{s_q}^*(\bq) \cdot\bh_{s_k}^*(\bk).
\end{align*}
Renaming the indices, we can write
\begin{align*}
    &L^3\sum_{\bp+\bq+\bk=0}\sum_{s_p,s_q} \big[ a_{s_p}^*(\bp,t)\rmd t + b_{s_p}^*(\bp)\circ \rmd W_t\big]s_q|\bq|a_{s_q}^*(\bq,t) \bh_{s_p}^*(\bp)\times\bh_{s_q}^*(\bq) \cdot\bh_{s_k}^*(\bk)\\
    &= L^3\sum_{\bp+\bq+\bk=0}\sum_{s_p,s_q} \big[ a_{s_q}^*(\bq,t)\rmd t + b_{s_q}^*(\bq)\circ \rmd W_t\big]s_p|\bp|a_{s_p}^*(\bp,t) \bh_{s_q}^*(\bq)\times\bh_{s_p}^*(\bp) \cdot\bh_{s_k}^*(\bk)\\
    &= -L^3\sum_{\bp+\bq+\bk=0}\sum_{s_p,s_q} \big[ a_{s_q}^*(\bq,t)\rmd t + b_{s_q}^*(\bq)\circ \rmd W_t\big]s_p|\bp|a_{s_p}^*(\bp,t) \bh_{s_p}^*(\bp)\times\bh_{s_q}^*(\bq) \cdot\bh_{s_k}^*(\bk).
\end{align*}
Thus, we arrive at
\begin{align*}
&\int_{\mathcal{D}} \rmd \bx_t(\bx,t) \times \bs{\omega}(\bx, t) \cdot \bh_{s_k}^*(\bk)e^{-i\bk\cdot\bx} \id\bx\\
&= \frac{L^3}{2}
\sum_{\bp+\bq+\bk=0}\sum_{s_p,s_q}
\Bigg[
(s_q|\bq|b_{s_p}^*(\bp)\circ \rmd W_t a_{s_q}^*(\bq,t)-s_p|\bp|b_{s_q}^*(\bq)\circ \rmd W_t a_{s_p}^*(\bp,t))
\\&+
(
a_{s_p}^*(\bp,t)\rmd t s_q|\bq|a_{s_q}^*(\bq,t)
- a_{s_q}^*(\bq,t)\rmd t s_p|\bp|a_{s_p}^*(\bp,t)
) 
\Bigg]\bh_{s_p}^*(\bp)\times\bh_{s_q}^*(\bq) \cdot\bh_{s_k}^*(\bk).
\end{align*}
Therefore,
\begin{align*}
    \rmd a_{s_k}(\bk, t) &= \frac{1}{4}
\sum_{\bp+\bq+\bk=0}\sum_{s_p,s_q}
\Bigg[
(s_q|\bq|b_{s_p}^*(\bp)\circ \rmd W_t a_{s_q}^*(\bq,t)-s_p|\bp|b_{s_q}^*(\bq)\circ \rmd W_t a_{s_p}^*(\bp,t))
\\&+
(
a_{s_p}^*(\bp,t)\rmd t s_q|\bq|a_{s_q}^*(\bq,t)
- a_{s_q}^*(\bq,t)\rmd t s_p|\bp|a_{s_p}^*(\bp,t)
) 
\Bigg]\bh_{s_p}^*(\bp)\times\bh_{s_q}^*(\bq) \cdot\bh_{s_k}^*(\bk).
\end{align*}

\subsection{LU Euler}\label{adx:luEuler}
Written in terms of the SALT model, the LU model is
\begin{equation*}
    \id \bv + \id \bx_t\cdot\nabla \bv + \bv^j\nabla \id\bx_t^j - \bv^j\nabla (\bs{\xi} \circ\id W_t)^j 
    = -\nabla  \id p.
\end{equation*}
Expanding the additional term gives
\begin{equation*}
    \nabla (\bs{\xi} \circ\id W_t)^j = \nabla\sum_{\bq}(b_{\pm}(\bq)\bh_{\pm}(\bq)e^{i\bq\cdot\bx} \circ\id W_t)^j = \sum_{\bq} i \bq b_{\pm}(\bq)\bh^j_{\pm}(\bq) e^{i\bq\cdot\bx} \circ\id W_t.
\end{equation*}
Thus we get
\begin{equation*}
    \bv^j\nabla (\bs{\xi} \circ\id W_t)^j = \sum_{\bp, \bq}\sum_{s_p, s_q} a_{s_p}(\bp,t) \bh^j_{s_p}(\bp) i \bq b_{s_q}(\bq)\bh^j_{s_q}(\bq) e^{i(\bp+\bq)\cdot\bx} \circ\id W_t.
\end{equation*}
Now projecting and renaming, we have
\begin{align*}
    &\int_{\mathcal{D}} \bv^j\nabla (\bs{\xi} \circ\id W_t)^j \cdot a_{s_k}^*(\bk,t)\bh^*_{s_k}(\bk)e^{-i\bk\cdot \bx}\id\bx =\\
    &= \sum_{\bp, \bq}\sum_{s_p, s_q}a_{s_k}^*(\bk,t)a_{s_p}(\bp,t) [b_{s_q}(\bq)\circ\id W_t] (i\bq\cdot\bh^*_{s_k}(\bk))(\bh_{s_p}(\bp)\cdot\bh_{s_q}(\bq)) \int_{\mathcal{D}} e^{(\bp+\bq-\bk)\cdot\bx} \id\bx\\
    &=L^3 \sum_{\bp+\bq+\bk=0}\sum_{s_p, s_q}a_{s_k}^*(\bk,t)a^*_{s_p}(\bp,t) [b^*_{s_q}(\bq)\circ\id W_t] (-i\bq\cdot\bh^*_{s_k}(\bk))(\bh^*_{s_p}(\bp)\cdot\bh^*_{s_q}(\bq))\\
    &= L^3 \sum_{\bp+\bq+\bk=0}\sum_{s_p, s_q} (f_k^{pq})^* a_{s_k}^*(\bk,t)(a^*_{s_p}(\bp,t) [b^*_{s_q}(\bq)\circ\id W_t]+a^*_{s_q}(\bq,t) [b^*_{s_p}(\bp)\circ\id W_t])
\end{align*}
with
\begin{equation*}
    f_k^{pq} = (-i(\bp+\bq)\cdot\bh_{s_k}(\bk))(\bh_{s_p}(\bp)\cdot\bh_{s_q}(\bq)).
\end{equation*}
Note that, due to the triad condition $\bp+\bq+\bk=0$,
\begin{equation*}
    f_k^{pq} = f_p^{kq} = f_q^{kp} = 0
\end{equation*}
so that the difference term between SALT and LU vanishes in the helical projection and the two projected models coincide.

\section{Supplementary Numerics}
\label{sec:suppl_num}

\subsection{Calibration of the Noise Amplitude}
\label{sec:noise_calib}

To calibrate the noise amplitude for the data assimilation experiments, we rely on two forecast verification metrics. The rank histogram (or Talagrand histogram) and the continuous ranked probability score (CRPS). We evaluated these for both models on $64$ different noise amplitude vectors. The metrics were recorded by running the data assimilation/particle filtering experiment described in the main text for the noise amplitude vectors $\bb=[b_k, b_p, b_q]$ resulting from all possible combinations of $b_k \in \{0.05, 0.1, 0.2, 0.5\}$, $b_p\in\{0.025, 0.05, 0.1, 0.2\}$ and $b_q\in \{0.01, 0.02, 0.04, 0.1\}$.
We present the mean CRPS scores for the top $5$ tested noise amplitude vectors in Table~\ref{tab:crps}. The mean is taken across both models. Based on this calibration, we chose the case $\bb = [0.10, 0.05, 0.01]$ for the data assimilation experiment, as the other ones have inferior rank histograms, so we achieve a good balance between the visual judgment of rank histograms and the CRPS score.

\begin{table}[h]
\centering
\begin{tabular}{ l c c c l}
 $\bb$ &                    EST &           HST &           Mean    & Histograms\\ \hline
 $[0.05, 0.025, 0.01]$ &    $0.0282$ &      $0.0394$ &    $0.0338$  & Figure~\ref{fig:LU_5_2_1} +~\ref{fig:SALT_5_2_1}\\  
 $[0.05, 0.025, 0.02]$ &    $0.0338$ &      $0.0392$ &    $0.0365$  & Figure~\ref{fig:LU_5_2_2} +~\ref{fig:SALT_5_2_2}\\
 $[0.10, 0.025, 0.02]$ &    $0.0352$ &      $0.0385$ &    $0.0368$  & Figure~\ref{fig:LU_10_2_2} +~\ref{fig:SALT_10_2_2}\\
 $[0.10, 0.05, 0.02]$ &     $0.0393$ &      $0.0356$ &    $0.0375$  & Figure~\ref{fig:LU_10_5_2} +~\ref{fig:SALT_10_5_2}\\
 $[0.10, 0.05, 0.01]$ &     $0.0363$ &      $0.0409$ &    $0.0386$  & Figure~\ref{fig:LU_10_5_1} +~\ref{fig:SALT_10_5_1}
\end{tabular}
\caption{CRPS scores (lower is better) for the five best -- in terms of overall mean CRPS score -- tested noise amplitude vectors $\bb$ (first column). The mean CRPS score over the three modal energies is shown for the EST (second column) and HST (third column) models. The overall mean CRPS score for the respective noise amplitude is shown in the \emph{Mean} column (fourth column). Finally, we provide references to the figures showing the associated rank histograms (last column).}
\label{tab:crps}
\end{table}

\begin{figure}[h]

\centering
\begin{subfigure}{0.18\textwidth}
    \includegraphics[width=\textwidth]{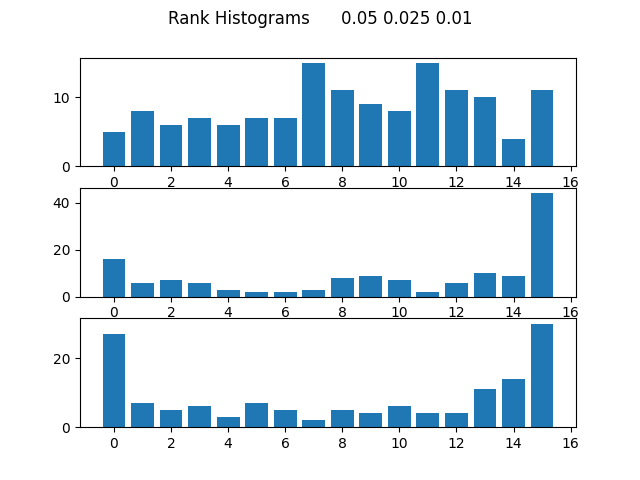}
    \caption{EST}
    \label{fig:LU_5_2_1}
\end{subfigure}
\hfill
\begin{subfigure}{0.18\textwidth}
    \includegraphics[width=\textwidth]{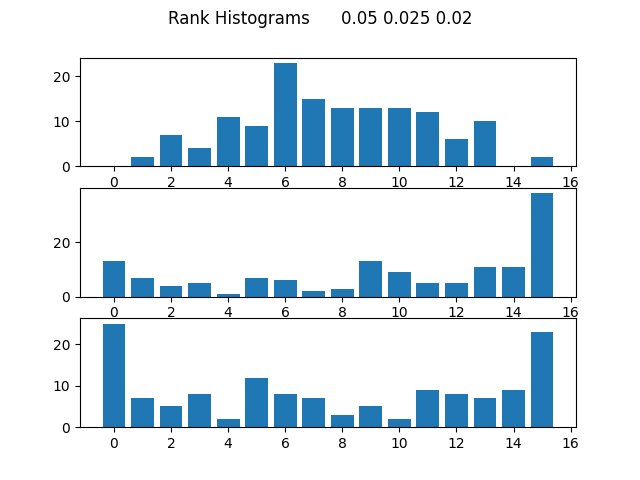}
    \caption{EST}
    \label{fig:LU_5_2_2}
\end{subfigure}
\hfill
\begin{subfigure}{0.18\textwidth}
    \includegraphics[width=\textwidth]{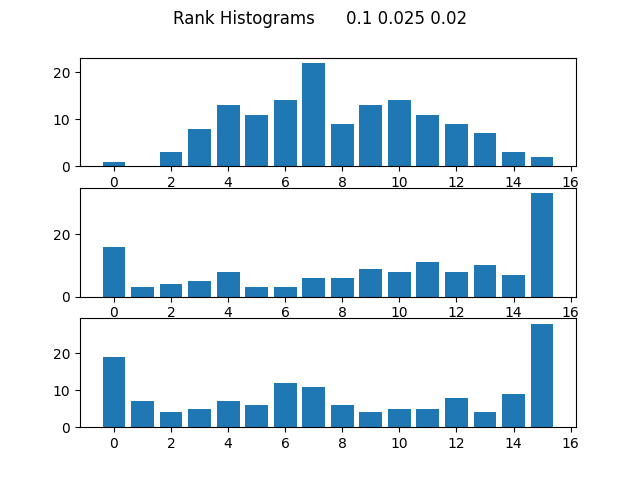}
    \caption{EST}
    \label{fig:LU_10_2_2}
\end{subfigure}
\hfill
\begin{subfigure}{0.18\textwidth}
    \includegraphics[width=\textwidth]{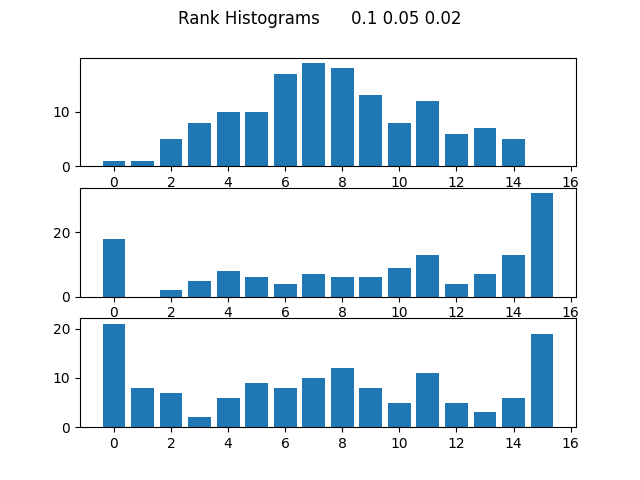}
    \caption{EST}
    \label{fig:LU_10_5_2}
\end{subfigure}
\hfill
\begin{subfigure}{0.18\textwidth}
    \includegraphics[width=\textwidth]{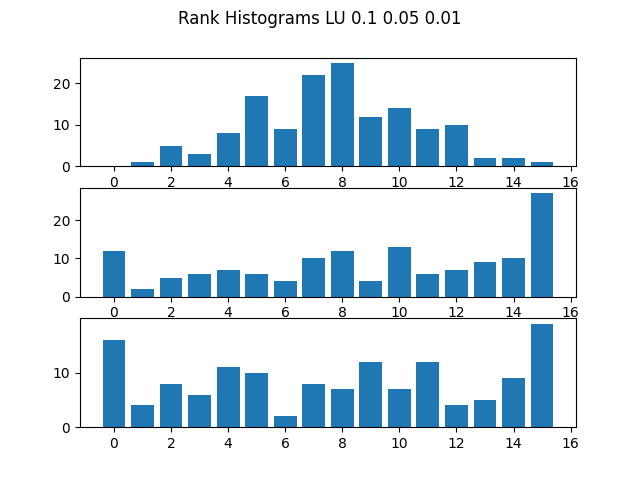}
    \caption{EST}
    \label{fig:LU_10_5_1}
\end{subfigure}

\begin{subfigure}{0.18\textwidth}
    \includegraphics[width=\textwidth]{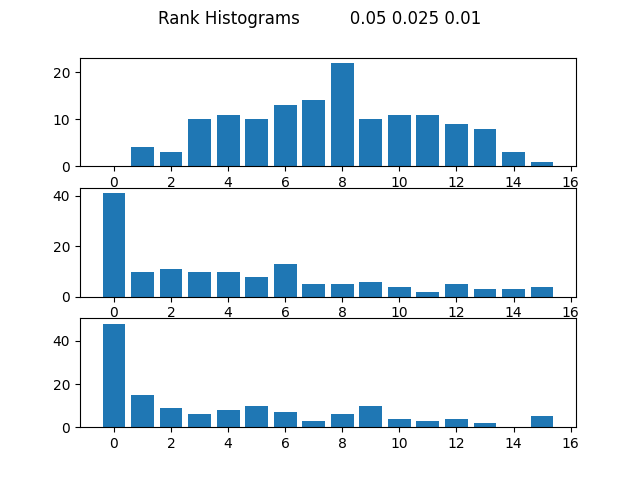}
    \caption{HST}
    \label{fig:SALT_5_2_1}
\end{subfigure}
\hfill
\begin{subfigure}{0.18\textwidth}
    \includegraphics[width=\textwidth]{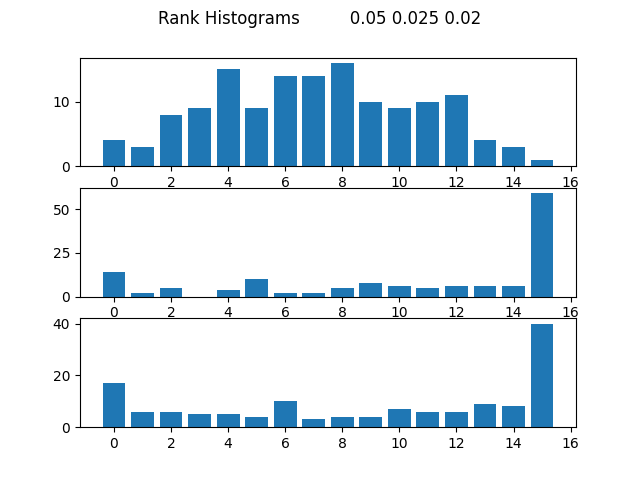}
    \caption{HST}
    \label{fig:SALT_5_2_2}
\end{subfigure}
\hfill
\begin{subfigure}{0.18\textwidth}
    \includegraphics[width=\textwidth]{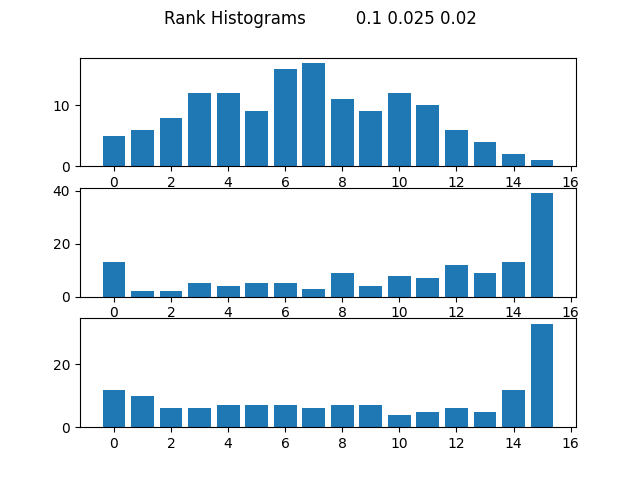}
    \caption{HST}
    \label{fig:SALT_10_2_2}
\end{subfigure}
\hfill
\begin{subfigure}{0.18\textwidth}
    \includegraphics[width=\textwidth]{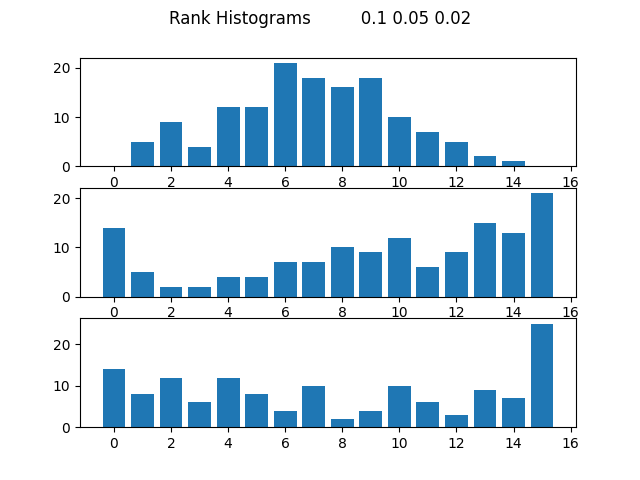}
    \caption{HST}
    \label{fig:SALT_10_5_2}
\end{subfigure}
\hfill
\begin{subfigure}{0.18\textwidth}
    \includegraphics[width=\textwidth]{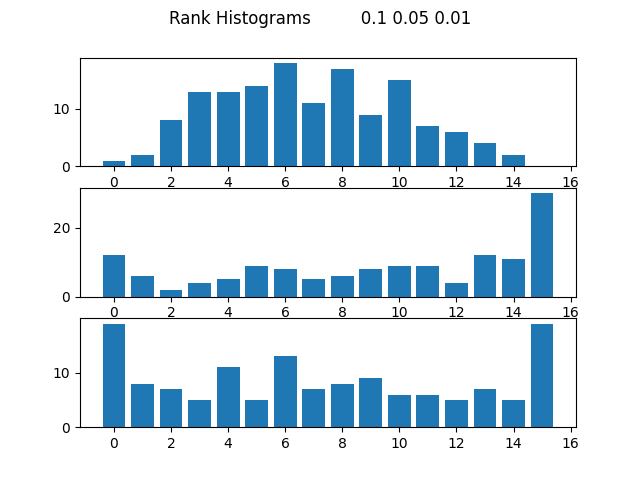}
    \caption{HST}
    \label{fig:SALT_10_5_1}
\end{subfigure}

\caption{Rank histograms for the $5$ best noise amplitude vectors in terms of CRPS score (see Table~\ref{tab:crps}) for the EST and HST models. Each individual graph shows the rank histogram of an ensemble of $15$ particles run up to a final time of $1400$, with data assimilation performed every $10$ time units. The top histogram in each subfigure represents the ensemble for mode $\bk$, the middle histogram represents the ensemble for mode $\bp$ and the bottom one for mode $\bq$.}
\label{fig:hist}

\end{figure}

\subsection{Data Assimilation Verification}
\label{sec:filter-verification}

\begin{figure}[h]

\centering
\begin{subfigure}{0.3\textwidth}
    \includegraphics[width=\textwidth]{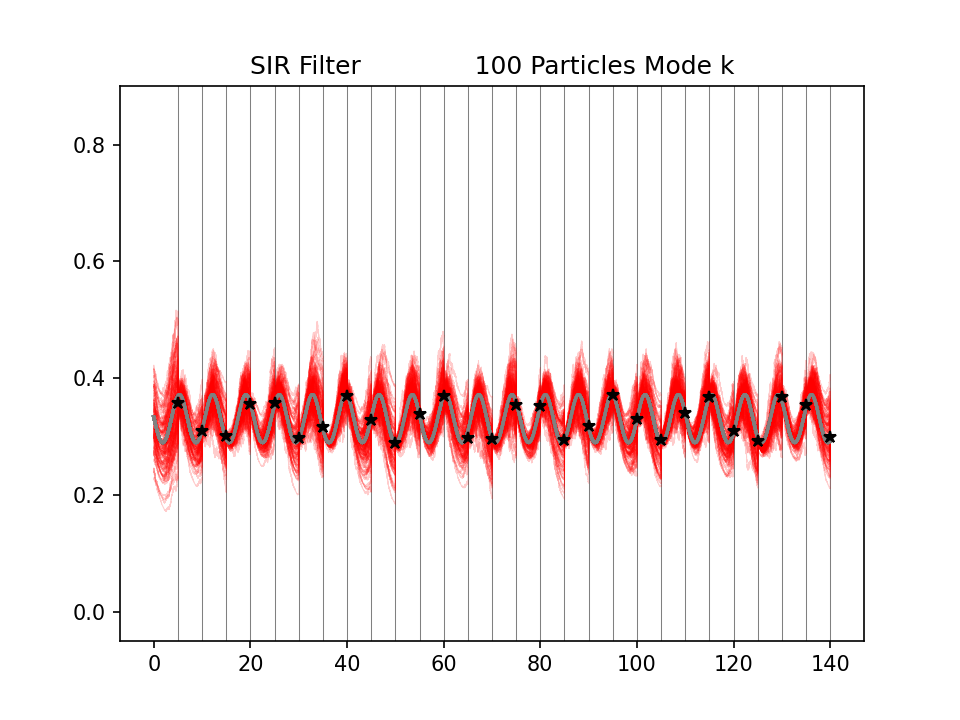}
    \caption{EST filter k.}
    \label{fig:LU_small_filter_k}
\end{subfigure}
\hfill
\begin{subfigure}{0.3\textwidth}
    \includegraphics[width=\textwidth]{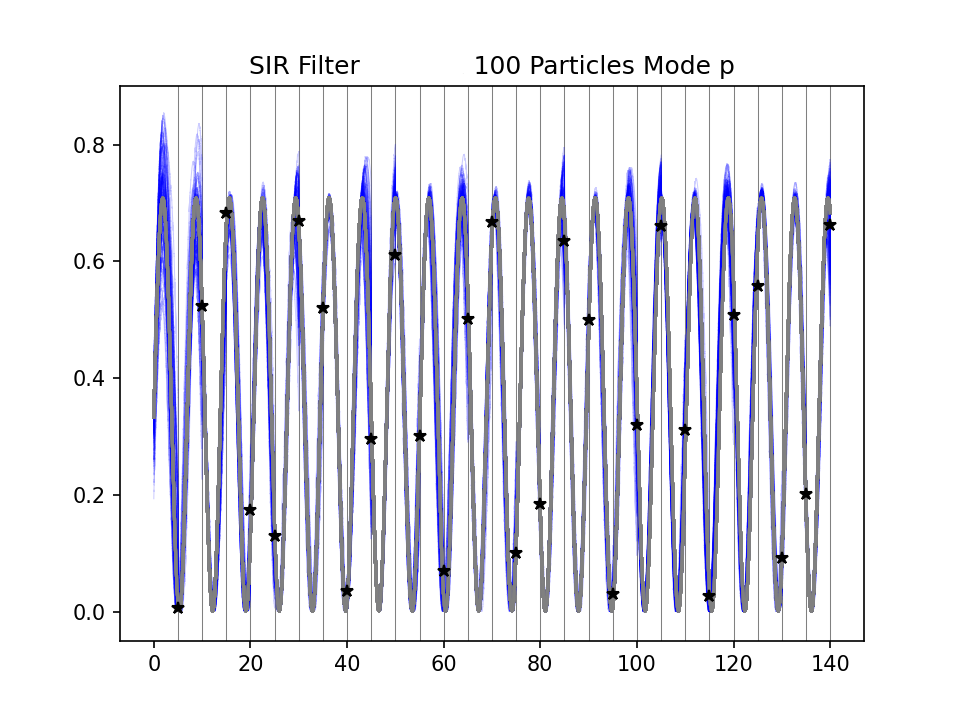}
    \caption{EST filter p.}
    \label{fig:LU_small_filter_p}
\end{subfigure}
\hfill
\begin{subfigure}{0.3\textwidth}
    \includegraphics[width=\textwidth]{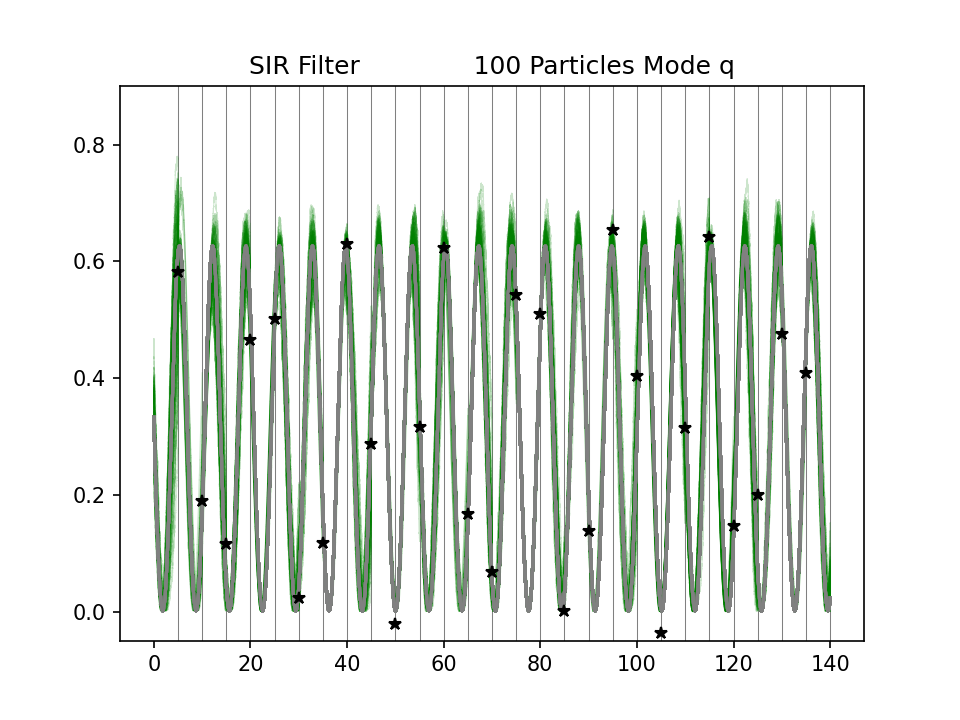}
    \caption{EST filter q.}
    \label{fig:LU_small_filter_q}
\end{subfigure}

\begin{subfigure}{0.3\textwidth}
    \includegraphics[width=\textwidth]{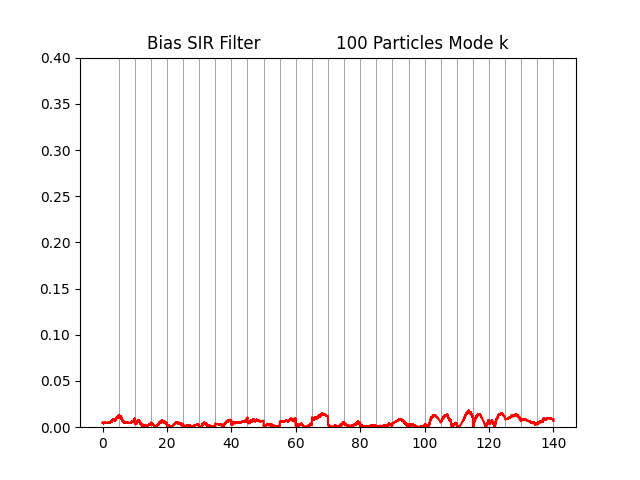}
    \caption{EST bias k.}
    \label{fig:LU_small_bias_k}
\end{subfigure}
\hfill
\begin{subfigure}{0.3\textwidth}
    \includegraphics[width=\textwidth]{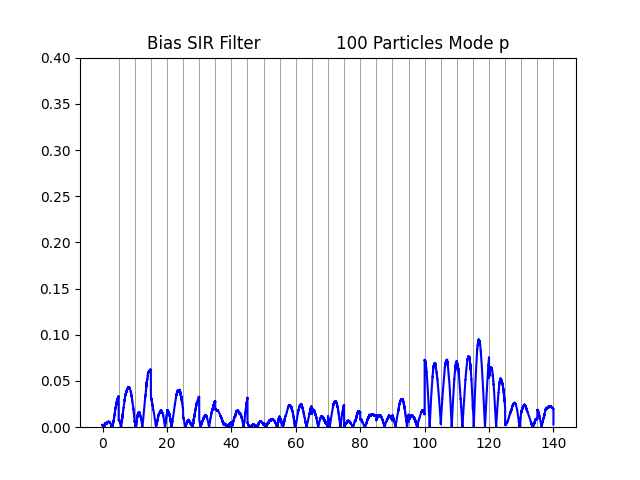}
    \caption{EST bias p.}
    \label{fig:LU_small_bias_p}
\end{subfigure}
\hfill
\begin{subfigure}{0.3\textwidth}
    \includegraphics[width=\textwidth]{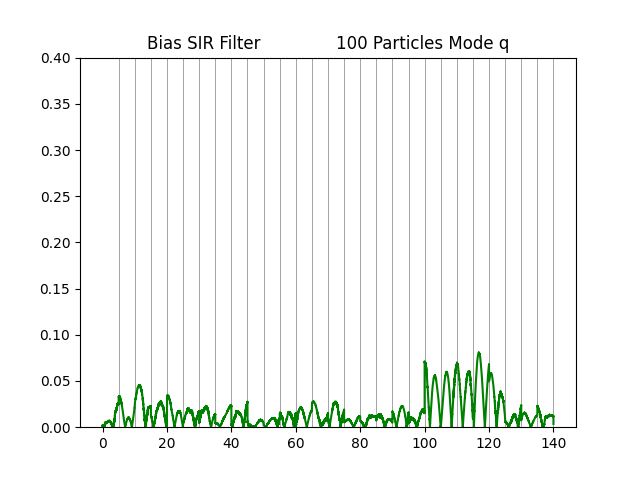}
    \caption{EST bias q.}
    \label{fig:LU_small_bias_q}
\end{subfigure}
\begin{subfigure}{0.3\textwidth}
    \includegraphics[width=\textwidth]{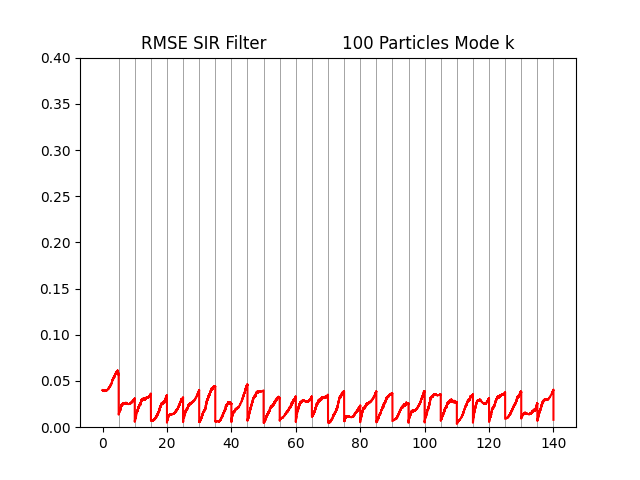}
    \caption{EST RMSE k.}
    \label{fig:LU_small_rmse_k}
\end{subfigure}
\hfill
\begin{subfigure}{0.3\textwidth}
    \includegraphics[width=\textwidth]{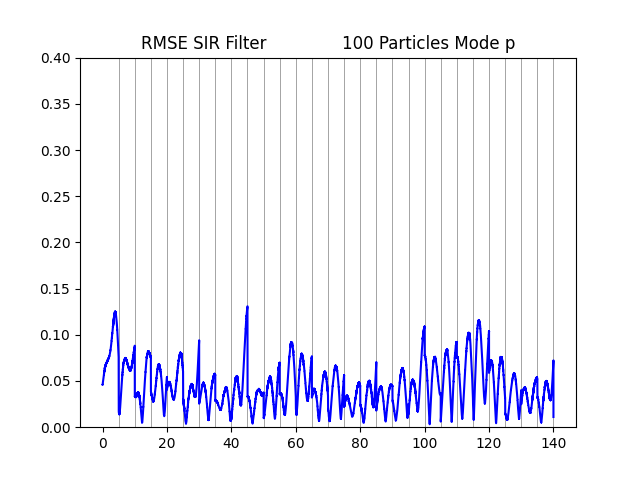}
    \caption{EST RMSE p.}
    \label{fig:LU_small_rmse_p}
\end{subfigure}
\hfill
\begin{subfigure}{0.3\textwidth}
    \includegraphics[width=\textwidth]{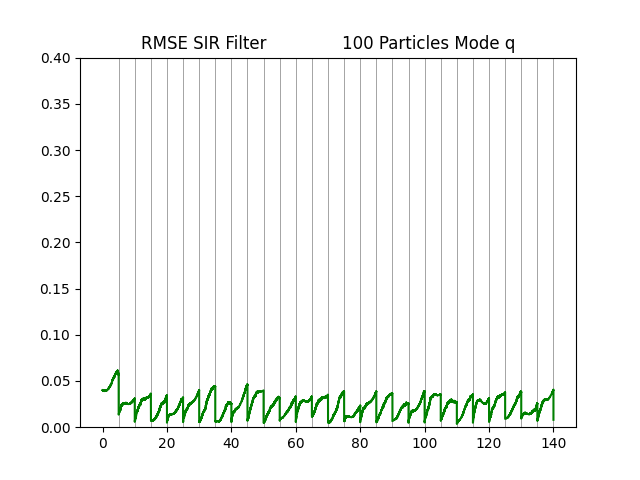}
    \caption{EST RMSE q.}
    \label{fig:LU_small_rmse_q}
\end{subfigure}

\caption{Filtering experiment for EST model using SIR particle filter with a small data assimilation interval of $5$ time units. (a-c) Ensemble evolution for $100$ particles in mode k (a, red), p (b, blue), and q (c, green). The signal (grey) is the deterministic model and the observations (black stars) are noisy and taken and assimilated every $5$ time units. (d-f) Bias of the filtering ensemble. (g-h) RMSE of the filtering ensemble.}
\label{fig:filtering_LU_small}

\end{figure}

\begin{figure}[h]

\centering
\begin{subfigure}{0.3\textwidth}
    \includegraphics[width=\textwidth]{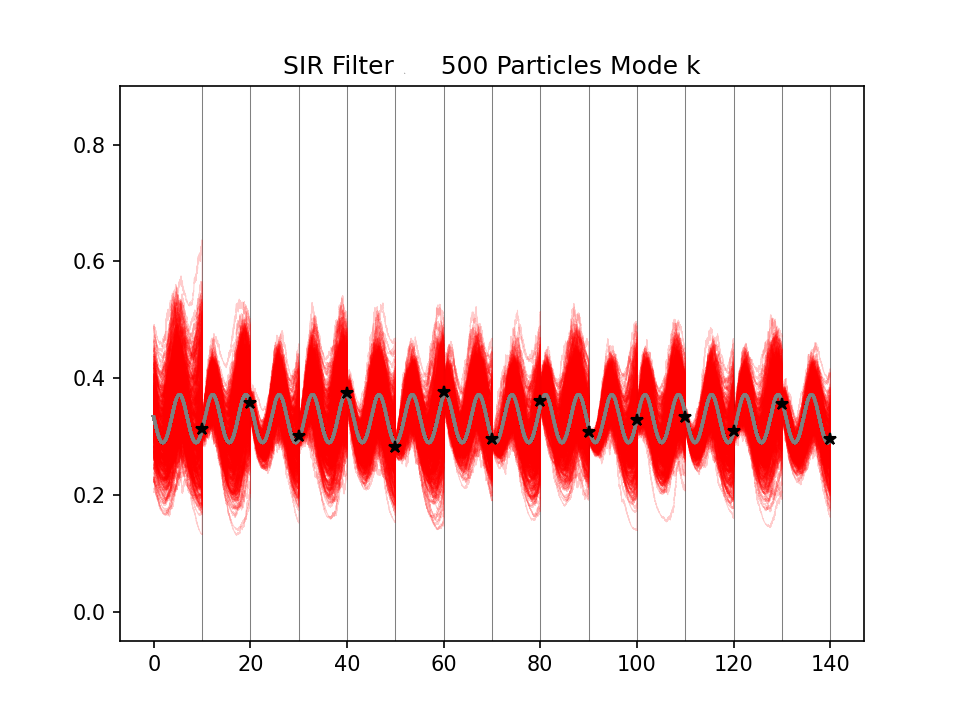}
    \caption{EST filter k.}
    \label{fig:LU_500_filter_k}
\end{subfigure}
\hfill
\begin{subfigure}{0.3\textwidth}
    \includegraphics[width=\textwidth]{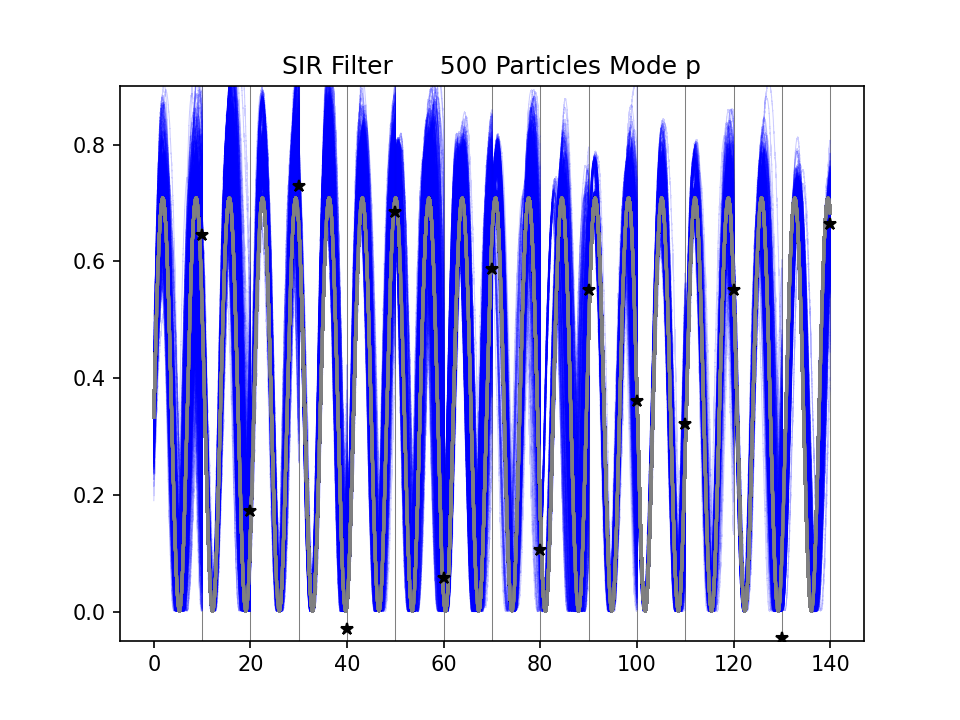}
    \caption{EST filter p.}
    \label{fig:LU_500_filter_p}
\end{subfigure}
\hfill
\begin{subfigure}{0.3\textwidth}
    \includegraphics[width=\textwidth]{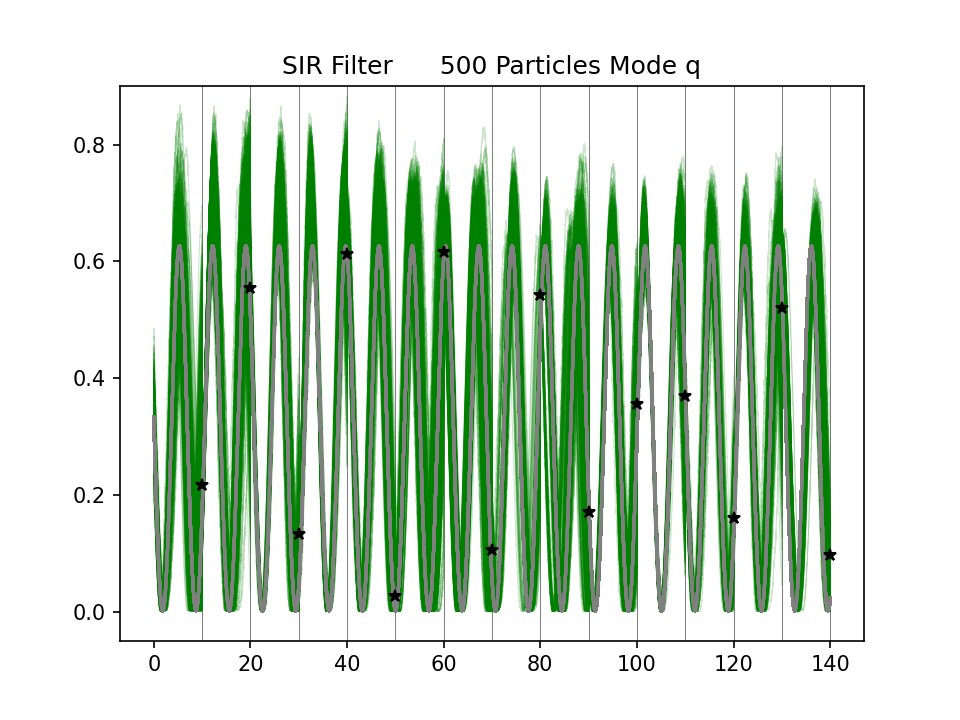}
    \caption{EST filter q.}
    \label{fig:LU_500_filter_q}
\end{subfigure}

\begin{subfigure}{0.3\textwidth}
    \includegraphics[width=\textwidth]{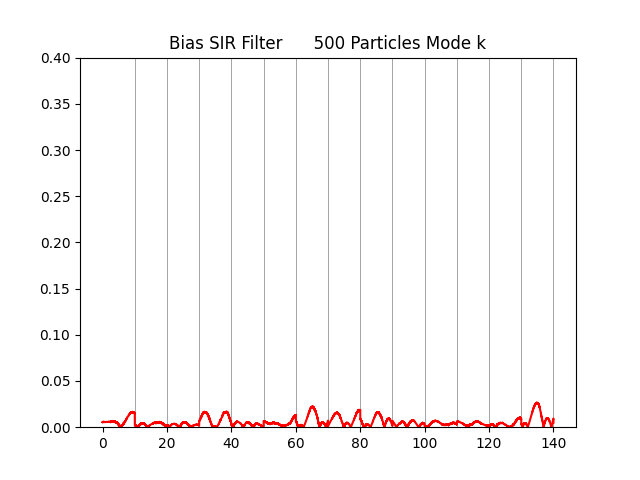}
    \caption{EST bias k.}
    \label{fig:LU_500_bias_k}
\end{subfigure}
\hfill
\begin{subfigure}{0.3\textwidth}
    \includegraphics[width=\textwidth]{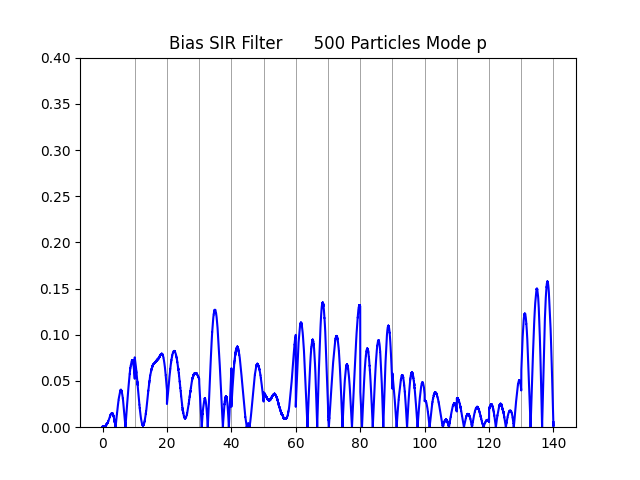}
    \caption{EST bias p.}
    \label{fig:LU_500_bias_p}
\end{subfigure}
\hfill
\begin{subfigure}{0.3\textwidth}
    \includegraphics[width=\textwidth]{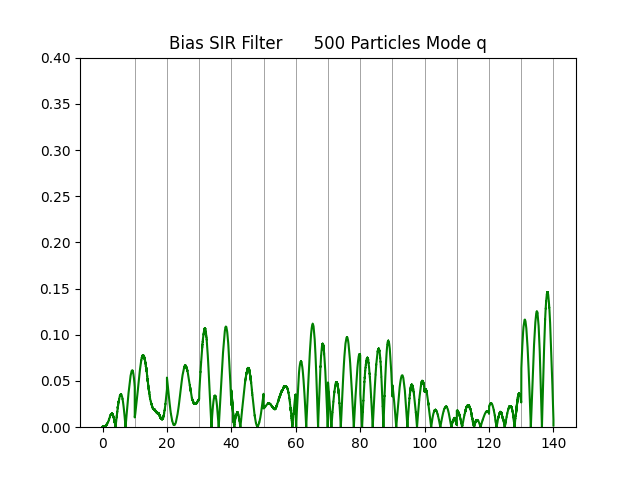}
    \caption{EST bias q.}
    \label{fig:LU_500_bias_q}
\end{subfigure}
\begin{subfigure}{0.3\textwidth}
    \includegraphics[width=\textwidth]{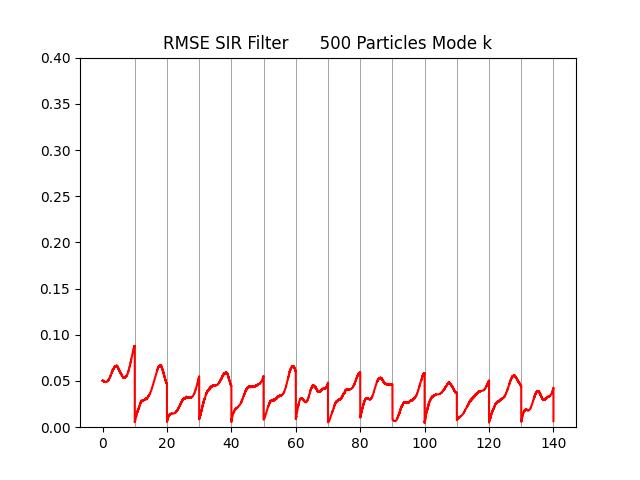}
    \caption{EST RMSE k.}
    \label{fig:LU_500_rmse_k}
\end{subfigure}
\hfill
\begin{subfigure}{0.3\textwidth}
    \includegraphics[width=\textwidth]{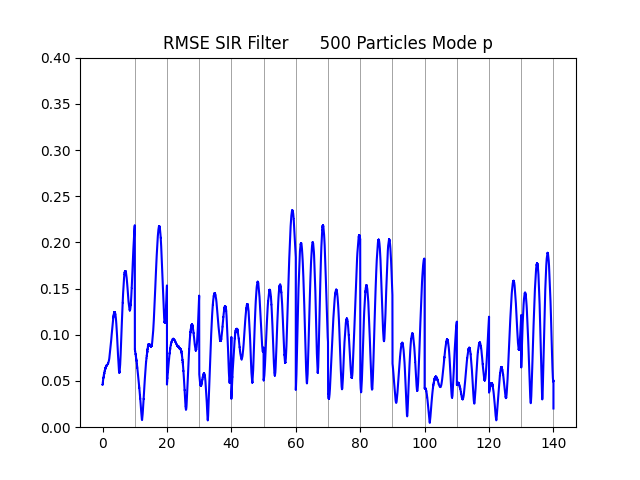}
    \caption{EST RMSE p.}
    \label{fig:LU_500_rmse_p}
\end{subfigure}
\hfill
\begin{subfigure}{0.3\textwidth}
    \includegraphics[width=\textwidth]{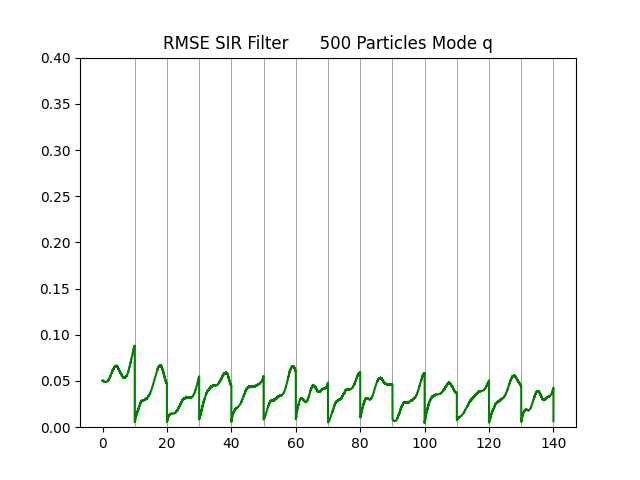}
    \caption{EST RMSE q.}
    \label{fig:LU_500_rmse_q}
\end{subfigure}

\caption{Filtering experiment for EST model using SIR particle filter with a large particle ensemble of $500$ members. (a-c) Ensemble evolution for $500$ particles in mode k (a, red), p (b, blue), and q (c, green). The signal (grey) is the deterministic model and the observations (black stars) are noisy and taken and assimilated every $5$ time units. (d-f) Bias of the filtering ensemble. (g-h) RMSE of the filtering ensemble.}
\label{fig:filtering_LU_500}

\end{figure}

\begin{figure}[h]

\centering
\begin{subfigure}{0.48\textwidth}
    \includegraphics[width=\textwidth]{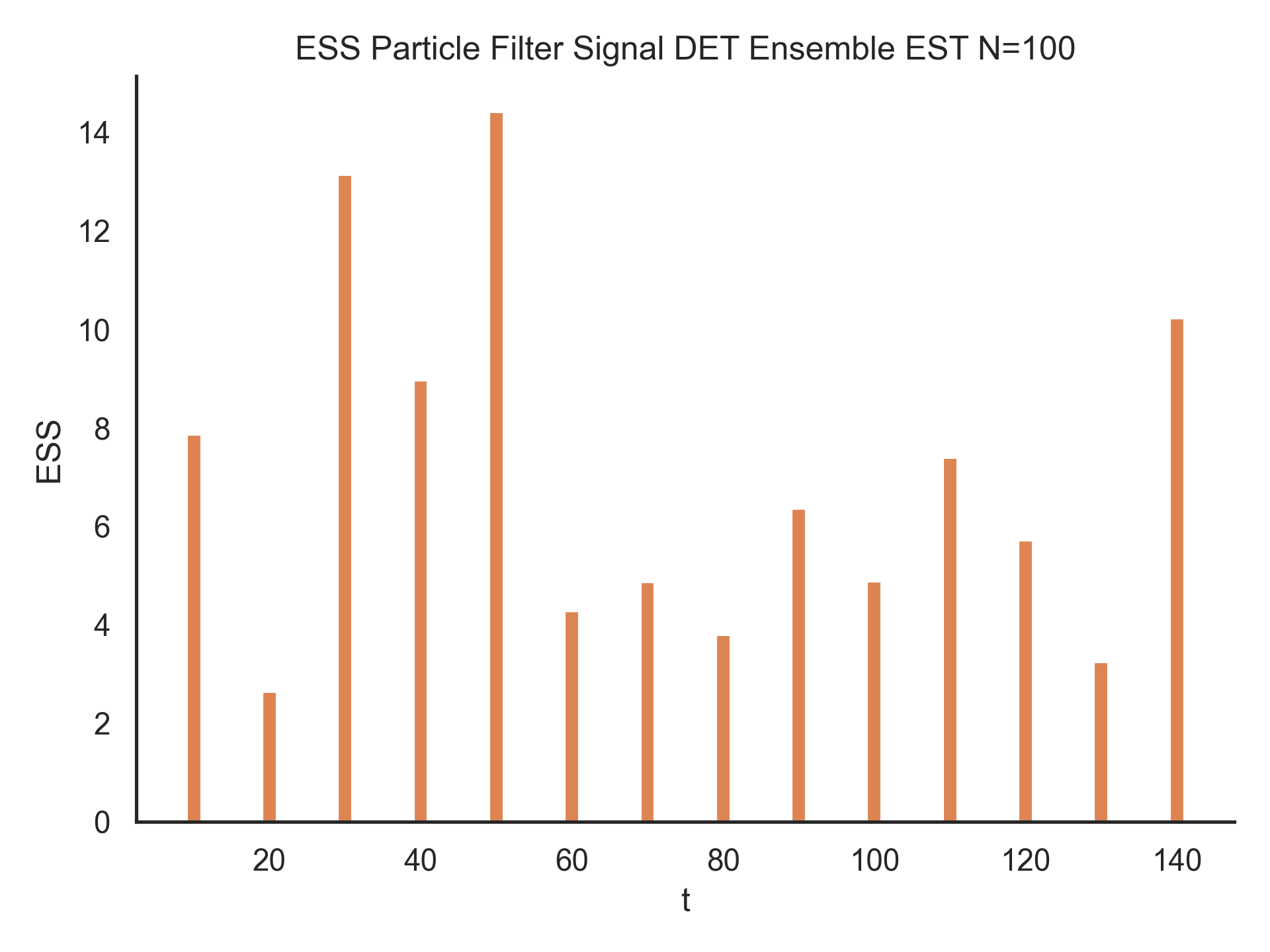}
    \caption{ESS for EST experiment.}
    \label{fig:ess_est}
\end{subfigure}
\hfill
\begin{subfigure}{0.48\textwidth}
    \includegraphics[width=\textwidth]{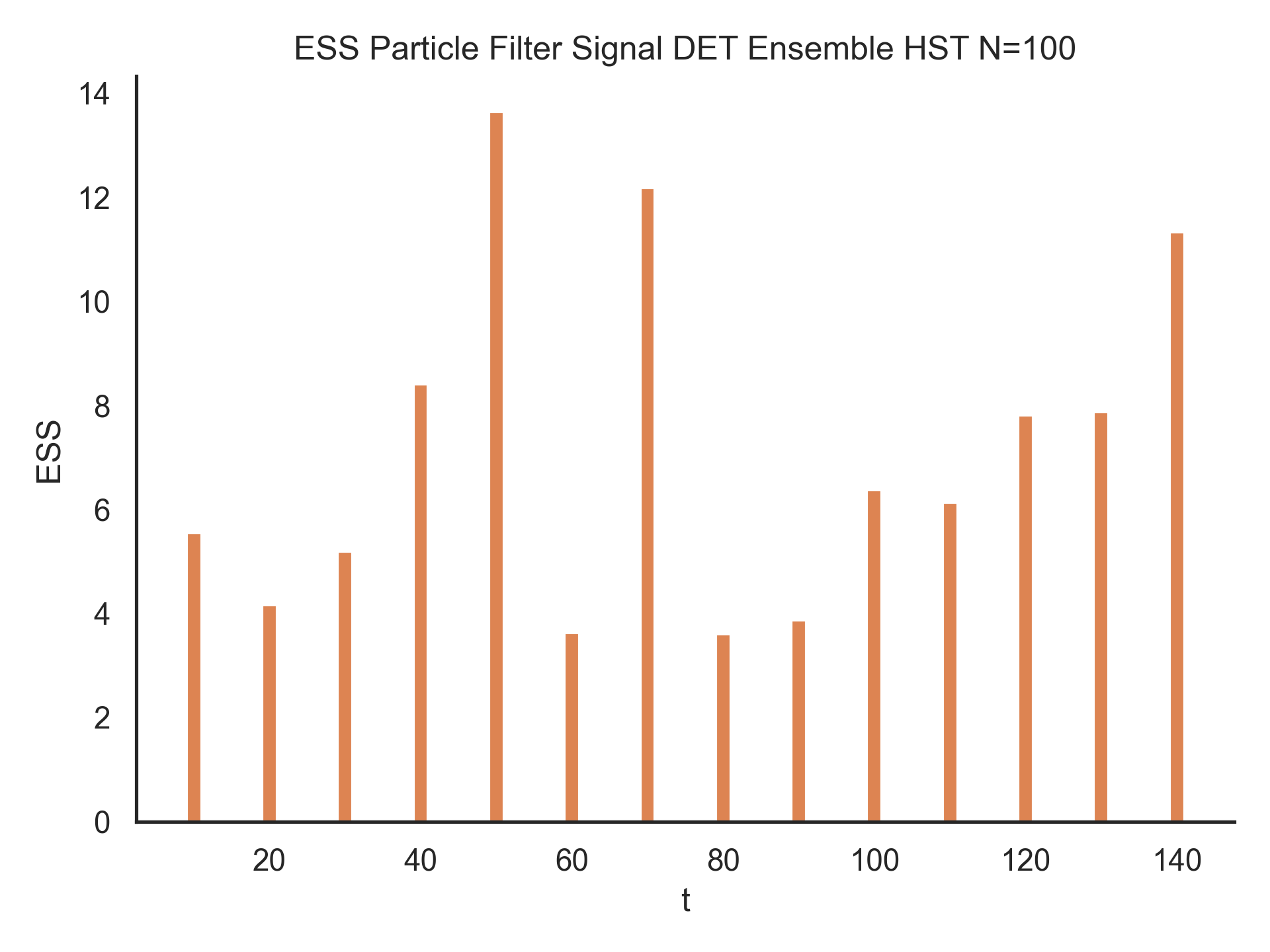}
    \caption{ESS for HST experiment.}
    \label{fig:ess_hst}
\end{subfigure}
\hfill

\caption{Typical ESS for the filtering experiments. See Remark~\ref{rem:ess} in the main text.}
\label{fig:filtering_ess}

\end{figure}

\begin{figure}[h]

\centering
\begin{subfigure}{0.3\textwidth}
    \includegraphics[width=\textwidth]{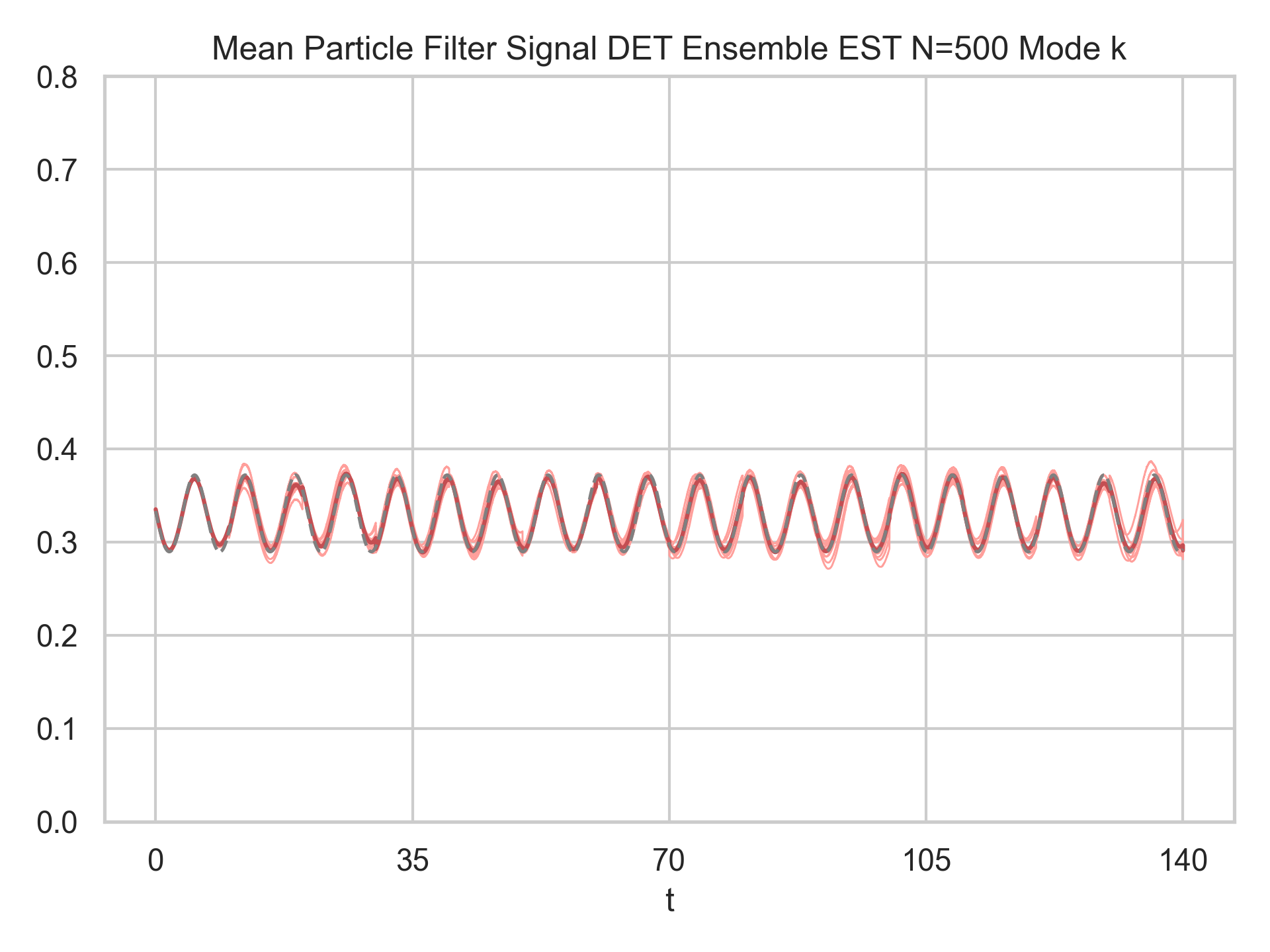}
    \caption{EST filter k.}
    \label{fig:EST_mean_k}
\end{subfigure}
\hfill
\begin{subfigure}{0.3\textwidth}
    \includegraphics[width=\textwidth]{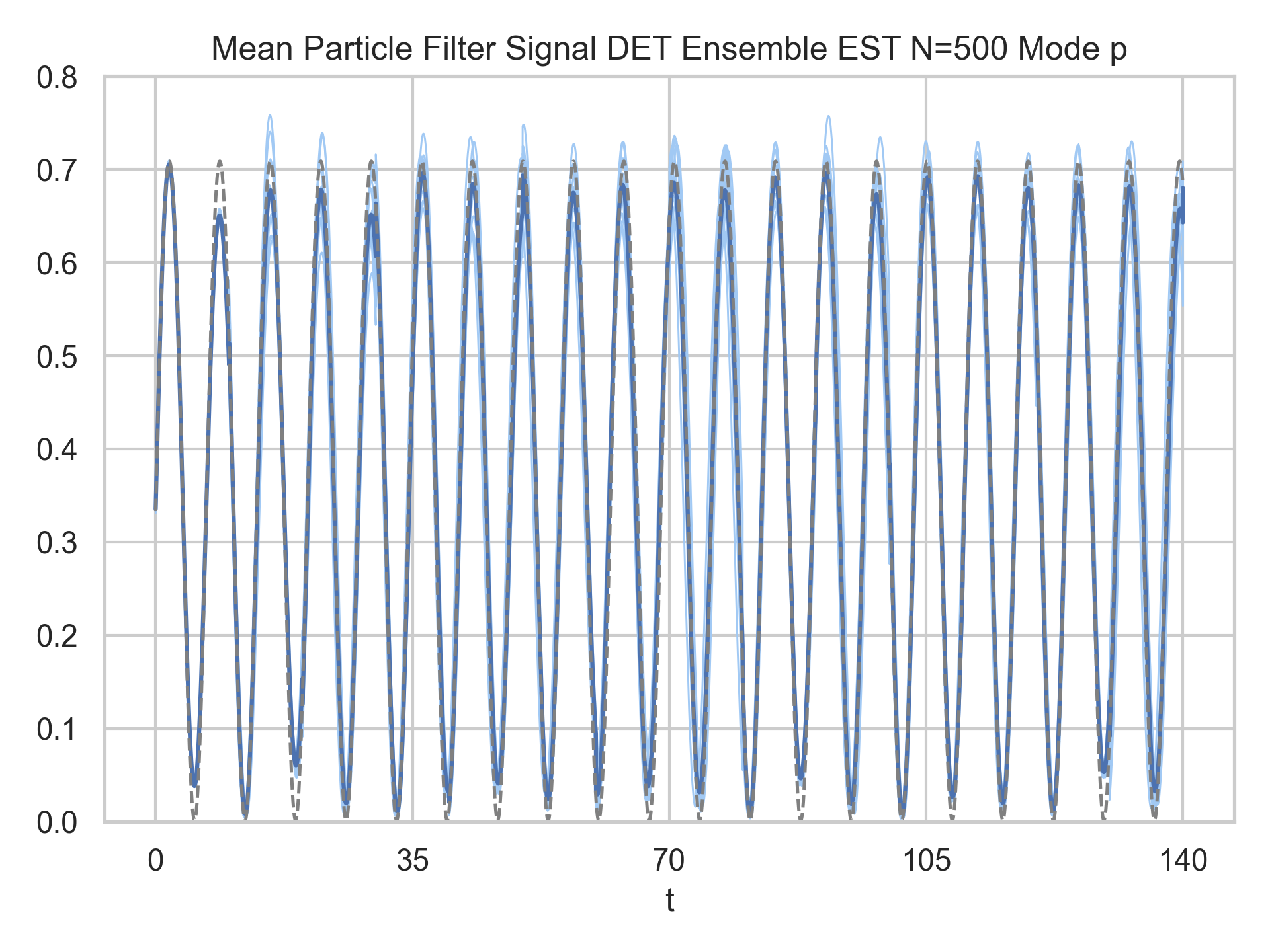}
    \caption{EST filter p.}
    \label{fig:EST_mean_p}
\end{subfigure}
\hfill
\begin{subfigure}{0.3\textwidth}
    \includegraphics[width=\textwidth]{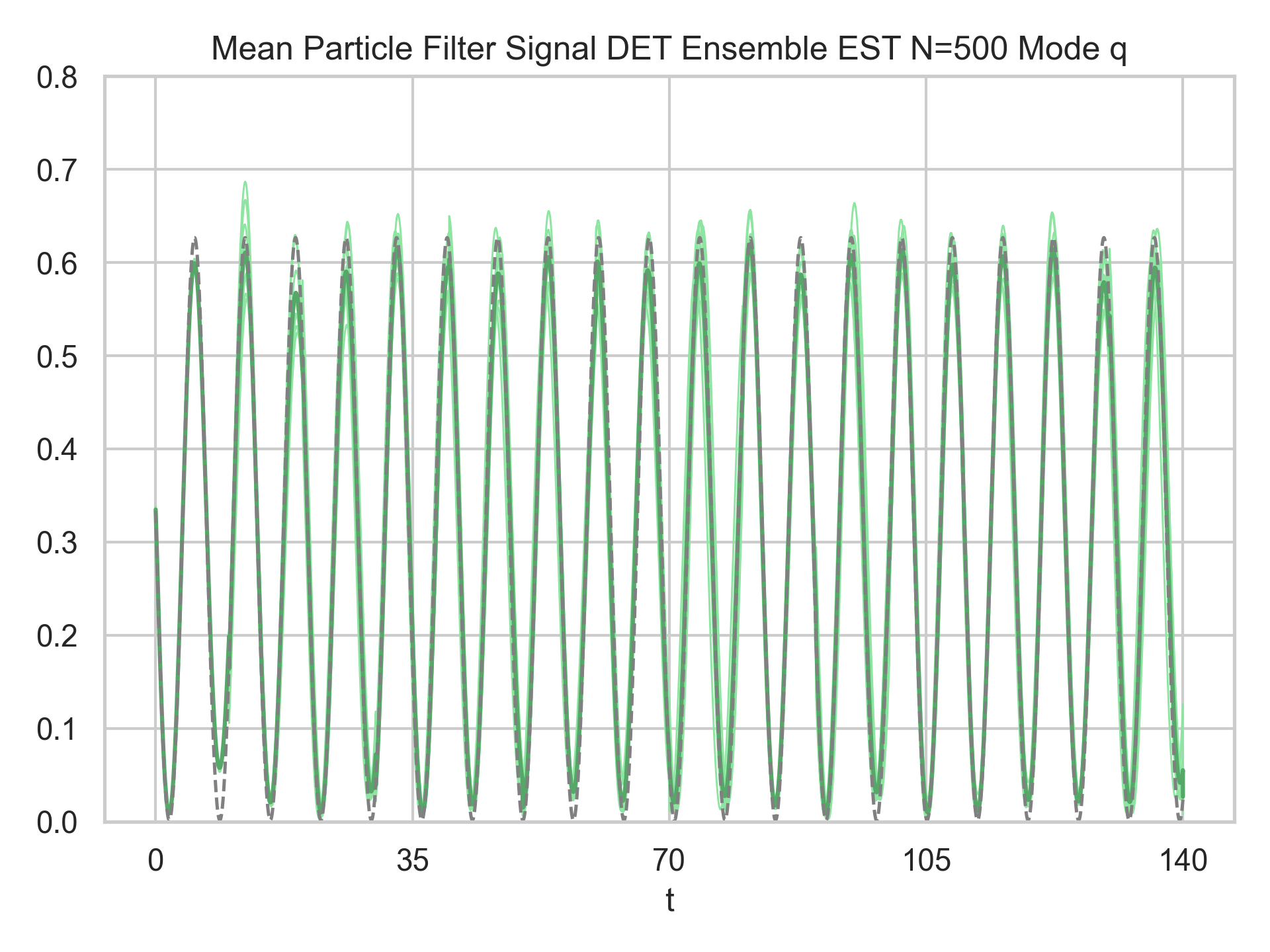}
    \caption{EST filter q.}
    \label{fig:EST_mean_q}
\end{subfigure}

\begin{subfigure}{0.3\textwidth}
    \includegraphics[width=\textwidth]{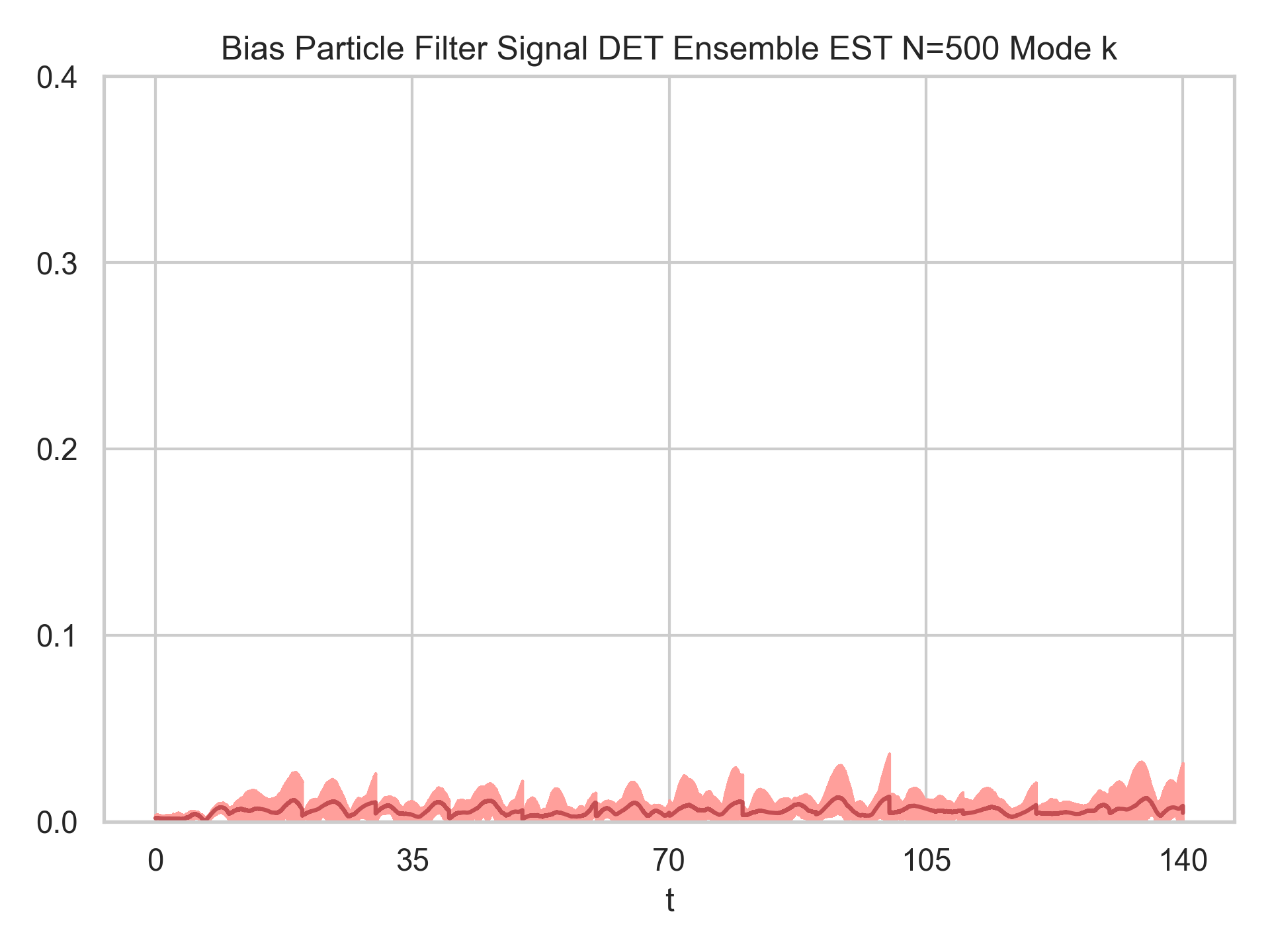}
    \caption{EST bias k.}
    \label{fig:EST_bias_k_mean}
\end{subfigure}
\hfill
\begin{subfigure}{0.3\textwidth}
    \includegraphics[width=\textwidth]{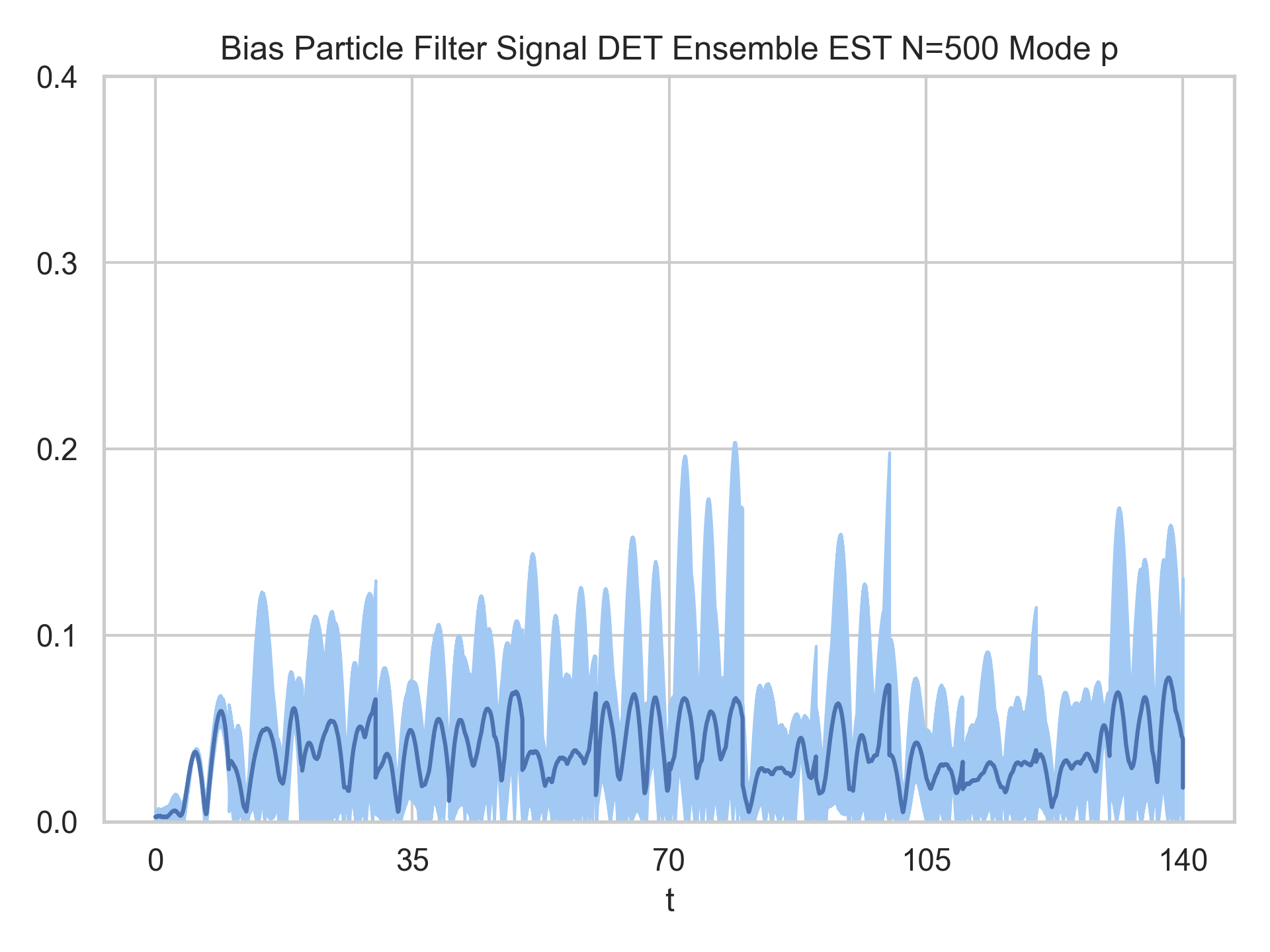}
    \caption{EST bias p.}
    \label{fig:EST_bias_p_mean}
\end{subfigure}
\hfill
\begin{subfigure}{0.3\textwidth}
    \includegraphics[width=\textwidth]{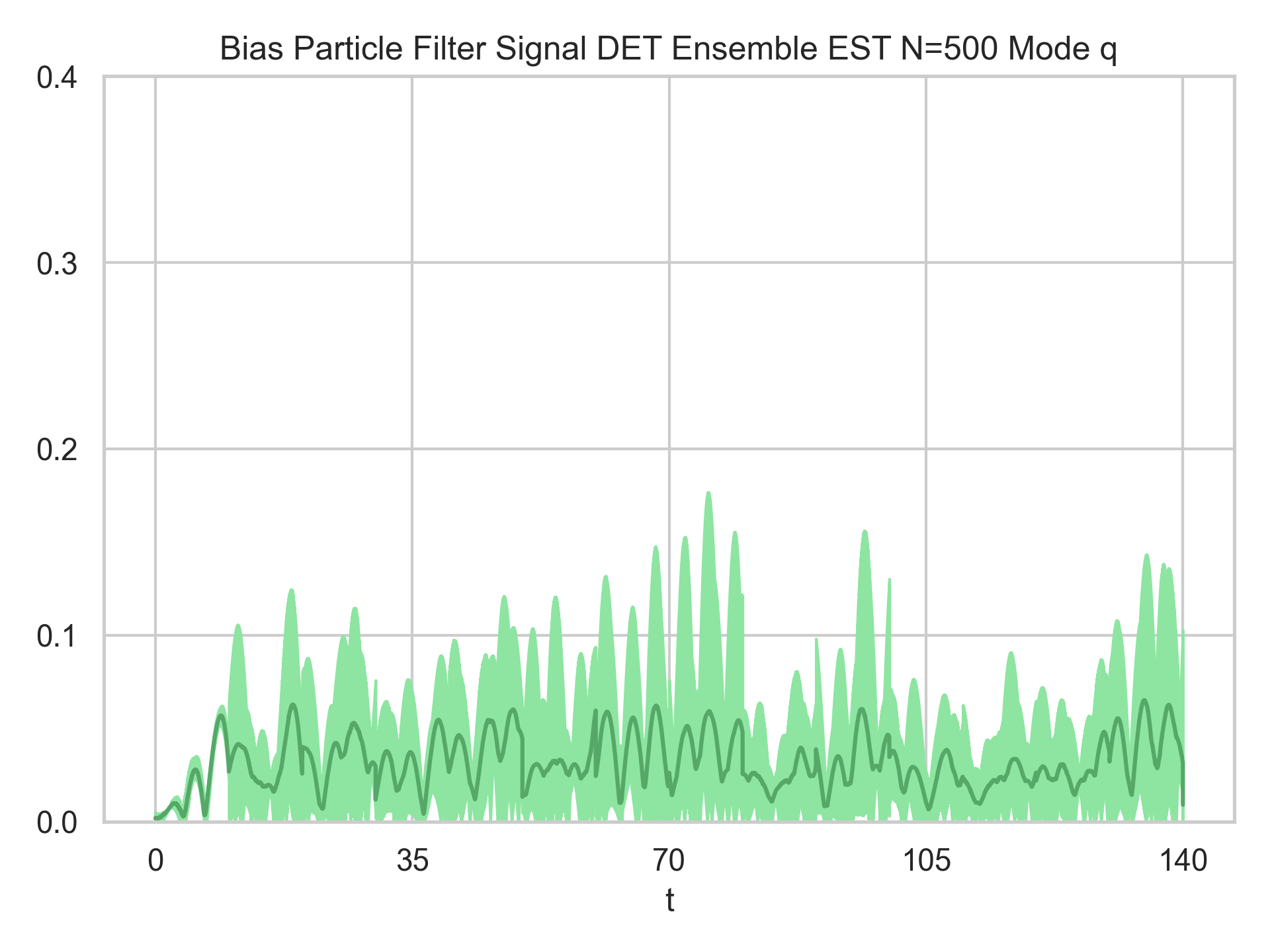}
    \caption{EST bias q.}
    \label{fig:EST_bias_q_mean}
\end{subfigure}
\begin{subfigure}{0.3\textwidth}
    \includegraphics[width=\textwidth]{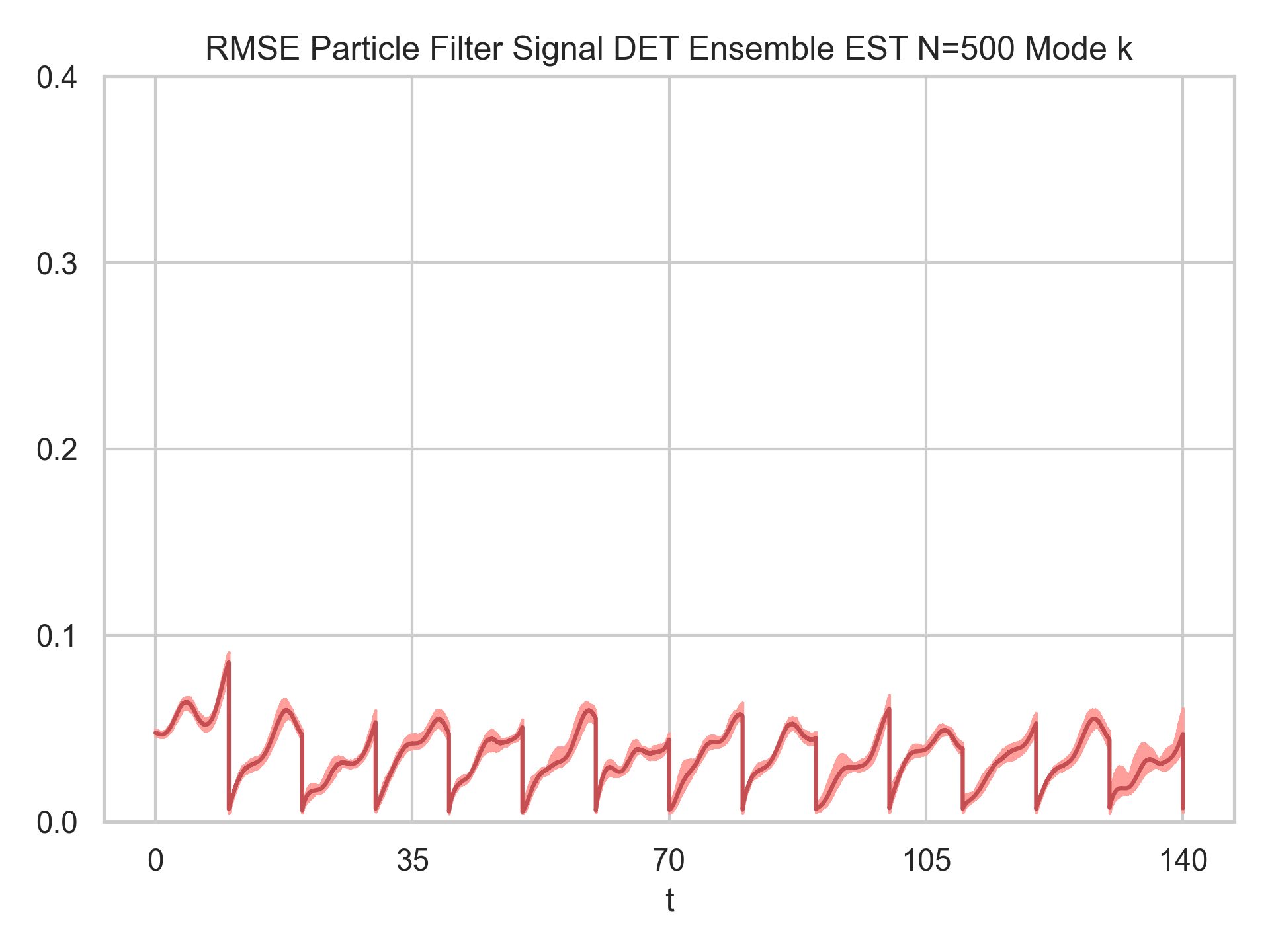}
    \caption{EST RMSE k.}
    \label{fig:EST_rmse_k_mean}
\end{subfigure}
\hfill
\begin{subfigure}{0.3\textwidth}
    \includegraphics[width=\textwidth]{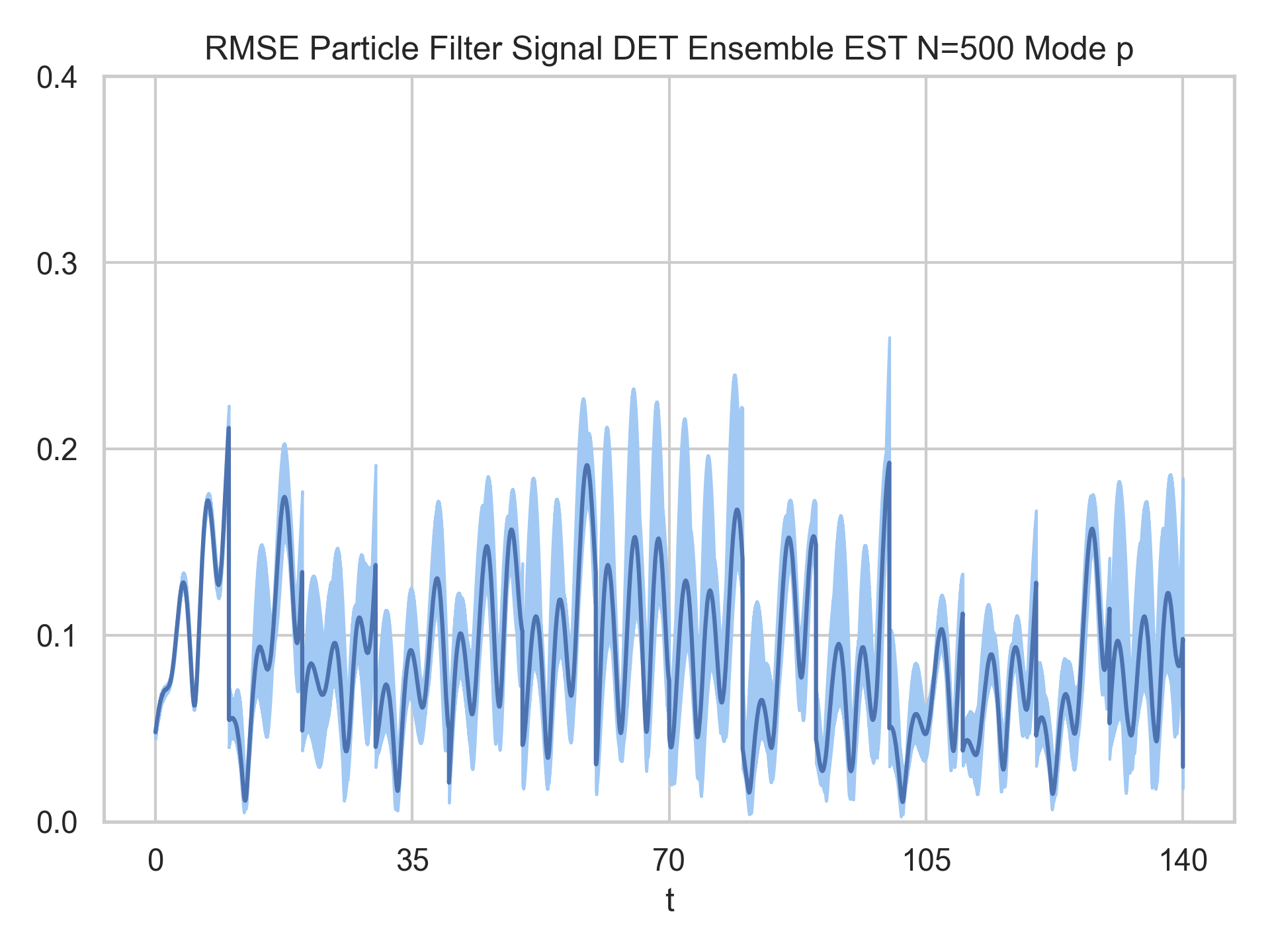}
    \caption{EST RMSE p.}
    \label{fig:EST_rmse_p_mean}
\end{subfigure}
\hfill
\begin{subfigure}{0.3\textwidth}
    \includegraphics[width=\textwidth]{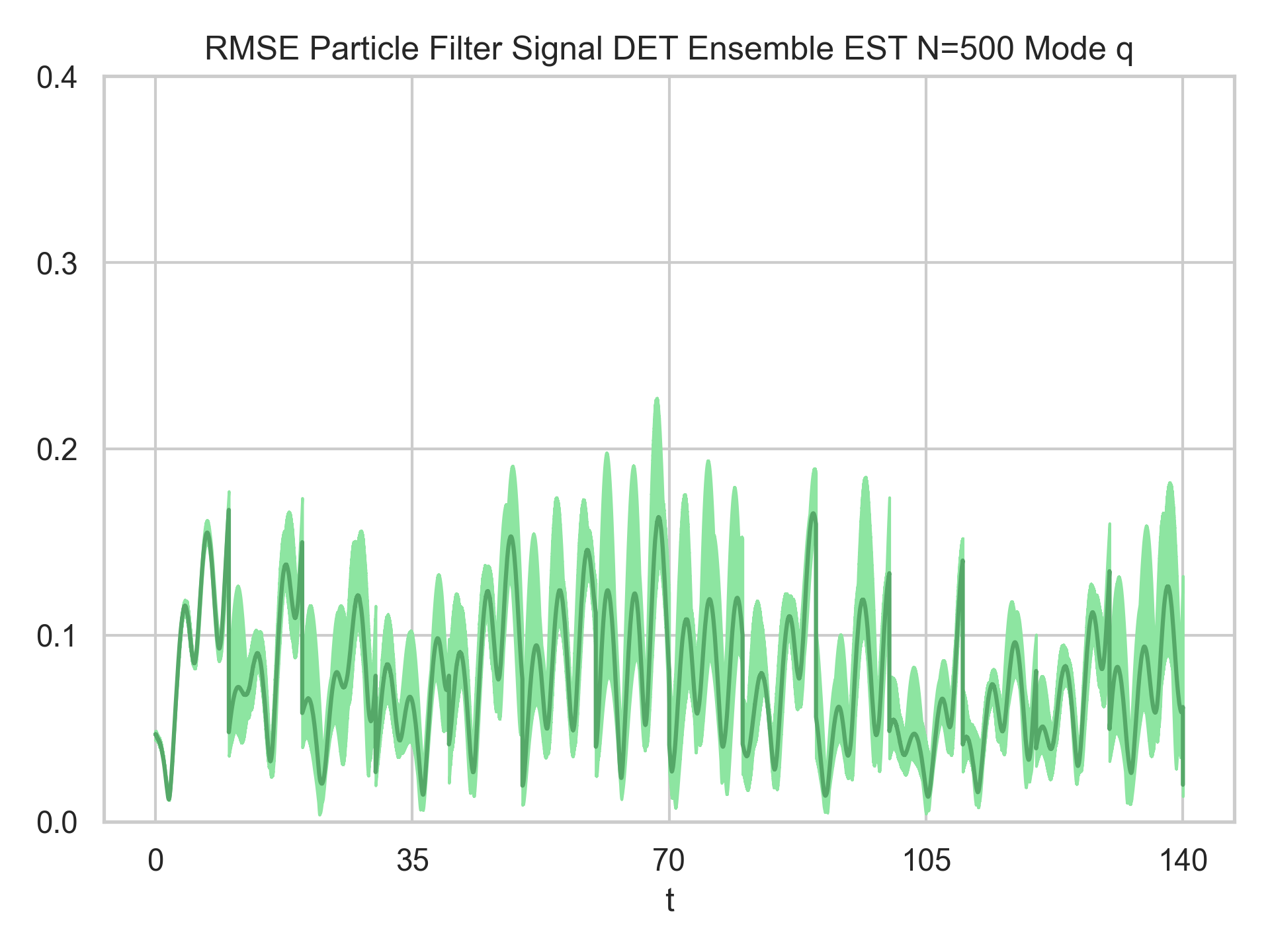}
    \caption{EST RMSE q.}
    \label{fig:EST_rmse_q_mean}
\end{subfigure}

\caption{Statistics of the mean ensemble over 10 independent runs of the particle filter with $500$ ensemble members for the EST model. See Remark~\ref{rem:stats_EST} in the main text.}
\label{fig:EST_mean_stats}
\end{figure}

\end{appendices}

\end{document}